\documentclass[amsfonts,12pt]{article}%
\usepackage{amssymb,amsmath,latexsym}
\usepackage{graphicx}
\usepackage{amsmath}
\usepackage{amsfonts}
\usepackage{amssymb}%
\providecommand{\U}[1]{\protect\rule{.1in}{.1in}}
\def \ed{\end{document}}
\numberwithin{equation}{section}
\def \n1{\newpage}
\def \L1{\frak L}

\def \U{{\mathcal U}}

\def \bqn{\begin{equation}}
\def \9{\end{equation}}
\def \3{\begin{eqnarray*}}
\def \4{\end{eqnarray*}}
\def \1{\begin{eqnarray}}
\def \2{\end{eqnarray}}

\def \big1{\bigcap}

\def \fr{\frac}

\newcounter{corollary}

\newcounter{proposition}
\newcounter{definition}

\def \brem{\begin{remark}}
\def \erem{\end{remark}}
\def \bth{\begin{theorem}}
\def \eth{\end{theorem}}
\def \bpr{\begin{proposition}}
\def \epr{\end{proposition}}
\def \bprf{\begin{proof}}
\def \eprf{\end{proof}}
\def \bex{\begin{example}}
\def \eex{\end{example}}
\def \bprf{\begin{proof}}
\def \eprf{\end{proof}}
\def \blem{\begin{lemma}}
\def \elem{\end{lemma}}
\newcounter{theorem}

\def \big{\bigcap}

\def \bl{\begin{lemma}}
\def \bcor{\begin{corollary}}
\def \ecor{\end{corollary}}
\def \el{\end{lemma}}
\def \beq*{\begin{eqnarray*}}
\def \eeq*{\end{eqnarray*}}
\def \6{\vspace*{7mm}}

\def \s1{\sqrt}

\def \mb{\mbox}

\def \bt{\begin{tabular}}
\def \et{\end{tabular}}

\def \l{\left}

\def \r{\right}
\def \hs1{\hspace*{3mm}}
\def \q2{\hspace*{9mm}}

\def \un1{\underline}
\def \mb{\mbox}
\def \vs1{\vspace*{4mm}}

\def \ba{\begin{array}}
\def \ea{\end{array}}

\newcommand{\ec}{\end{center}}
\newcommand{\bc}{\begin{center}}
\newcommand{\be}{\begin{equation}}

\newcommand{\ee}{\end{equation}}
\newcommand{\bn}{\begin{enumerate}}
\newcommand{\en}{\end{enumerate}}
\newcommand{\bi}{\begin{itemize}}
\newcommand{\ei}{\end{itemize}}
\newtheorem{theorem}{Theorem}

\newtheorem{corollary}{Corollary}

\newtheorem{example}{Example}

\newtheorem{lemma}{Lemma}

\newtheorem{proposition}{Proposition}
\newtheorem{remark}{Remark}

\newenvironment{proof}[1][Proof]{\textbf{#1.} }{\
\rule{0.5em}{0.5em}}
\begin{document}

\title{Higher order self-adjoint operators with polynomial coefficients}
\author{H. Azad$^{*}$, A. Laradji$^{*}$ and M. T. Mustafa$^{**}$ \\
$^{*}$Department of Mathematics and Statistics, King Fahd University of \\ Petroleum and Minerals, Dhahran 31261, Saudi Arabia\\
$^{**}$Department of Mathematics, Statistics and Physics, \\Qatar University, Doha 2713, Qatar\\
hassanaz@kfupm.edu.sa, alaradji@kfupm.edu.sa, tahir.mustafa@qu.edu.qa\\
}
\date{}
\maketitle

\begin{abstract}
Algebraic and analytic aspects of self-adjoint operators of order four or more with polynomial coefficients  are investigated. As a consequence, a systematic way
of constructing such operators is given. The procedure is applied to obtain
many examples up to order 8; similar examples can be constructed for all even
order operators. In particular, a complete classification of all order 4
operators is given.

\end{abstract}

\vspace{-0.5in}

\begin{center}

\end{center}

\baselineskip=18pt


\renewcommand{\theequation}{\thesection.\arabic{equation}}

\section{Introduction}
The classification of self-adjoint second order operators with polynomial coefficients is a classical subject going back to Brenke \cite{brenke}.
This paper is a contribution to certain algebraic and analytic aspects of
higher order self-adjoint operators with polynomial coefficients. Its main aim
is to construct such operators. This involves determining the explicit
differential equations for the polynomial coefficients of the operators and
the boundary conditions which ensure self-adjointness. These operators are not
in general iterates of second order classical operators - as stated in
\cite{littlejohnpams1994}; (cf. \cite{littlejohn-krall89}). As the weight
function which makes these operators self-adjoint depends only on the first
two leading terms of the operator, therefore, if one can find a second
operator with the same weight function, the eigenpolynomials for both
operators would be the same; see Section 4.

We should point out that some of the most important recent contributions to this subject are due to Kwon, Littlejohn and Yoon \cite{kwon2001}, Bavnick \cite{bavinck95, bavinck98, bavinck2000}, Koekoek \cite{koekoek93} and Kockock-Koekock \cite{kockock2000}; see also the references therein.

A classical reference for higher order Sturm-Liouville theory is the book of Ince \cite[Chap.IX]{ince}. This theory was revived by Everitt in \cite{everitt57}; the paper \cite{everitt2001} by Everitt et al. deals with the same subject.

Classical references that deal with various aspects of polynomial solutions of differential equations are the references \cite{Rou, brenke, Boc, lesky62, Kra40}. More recent papers that deal with the same subject are \cite{koornwinder84, krall81, littlejohn84, Turbiner, azad1, azad2}. A reference for related topic of orthogonal polynomials is \cite{Sze}. A more recent reference for this topic which also has an extensive bibliography is the book \cite{Ism}.

Here is a more detailed description of the results of this paper. We consider linear differential operators with polynomial coefficients that map polynomials of degree $k$ into itself- for all $k$. Proposition 2.1 gives necessary and sufficient conditions for such an $n^{\rm th}$ order linear operator to be self-adjoint. The integrability and asymptotic properties of the weight function and its derivatives near the zeroes of the leading term given in Propositions 3.2 and 3.3 carry enough information to determine the form of the first two terms of the operator in specific cases. This is used further to determine the full operator, using the boundary conditions and the determining equations of Proposition 2.1. Although similar determining equations are known (cf. \cite{kwon2001}), boundary conditions involving all polynomial coefficients of the linear operator do not seem to have been considered earlier and these are equally crucial to the construction of operators given here.

The classical second order operators are completely determined by the integrability of
the associated weight function.

In general, if the operator is of the form $L(y)= a y^{(n)}+b y^{(n-1)}+ \cdots$ then, as shown in Proposition 3.4, for $n>2$,  the multiplicity of each root of $a$ is at least $2$ and its multiplicity in $b$ is less than its multiplicity in $a$ and it is of multiplicity at least $1$ in $b$. In particular, if the operator is of fourth order and the leading term has distinct roots, then every root occurs with multiplicity 2 and therefore the leading term must have exactly two distinct roots and their multiplicity in the next term is $1$ - as shown in Proposition 3.4.

In \cite{littlejohn-krall89}, fourth order Sturm-Liouville systems were given for ordinary weights and for weights involving distributions. In particular, for ordinary weights, the authors  found  operators that are iterates of second order ones. In this paper,  a systematic way of producing Sturm-Liouville systems with ordinary weights for all even orders is given.
We recover the fourth order operators obtained in \cite{littlejohn-krall89} as well as new classes that are not, in general, iterates of second order operators - as stated in \cite{littlejohnpams1994}; cf. \cite{littlejohn-krall89}.

The examples II.6, II.7 of fourth order operators in \cite{littlejohn-krall89} with weights involving delta and Heaviside functions  and the sixth order operator in Example 4.1 in \cite{kwon2001} with weight involving delta functions are solutions of the determining equations given in Section 4. For these 4th order operators of \cite{littlejohn-krall89}, the 3rd boundary condition in the sense of  (\ref{boundary-n-4}) fails at one or both the boundary points. For the sixth order operator of Example 4.1 in \cite{kwon2001} the last three boundary conditions in (\ref{boundary-n-6}) do not hold. Many such examples can be constructed but the solutions of the determining equations are too many to be listed efficiently (cf. Example 4.1, Section 4.3).

We also give examples of sixth and eighth order operators where all the boundary conditions hold. In fact similar examples can be constructed for any even order; such constructions involve increased computational complexity.

All the solutions presented in the examples in Section 4 have been first verified using Mathematica and then directly generated (from the same file) by Mathematica as LaTeX output for the paper.

\section{Algebraic aspects of higher order Sturm-Liouville theory}
Consider, on the space $C^{\infty }$, the $n$th-order linear operator $%
L=\sum\limits_{k=1}^{n}a_{k}(x)D^{k},$ where $D$ is the usual differential
operator and each $a_{k}:=a_{k}(x)$ is a polynomial of degree at most $k$.
In this way, for each natural number $N,$ the vector space $\mathbb{P}_{N}$
of all polynomials of degree at most $N$ is $L$-invariant. Our first
objective is to obtain conditions on the polynomials $a_{k}$ for the
existence of an inner product $\left\langle u,v\right\rangle =\int_{I}puvdx$
on $C^{\infty }$ for which $L$ is self-adjoint, and where the weight $p$ is
a piecewise smooth function on some real interval $I$. This smoothness
assumption is reasonable since it is satisfied by all weights of classical
orthogonal polynomials. For a function $f$ and an interval $J$ with
(possibly infinite) endpoints $\alpha <\beta ,$ $\partial J$ denotes the
boundary $\{\alpha ,\beta \}$ of $J,$ and $\left[ f\right] _{J}$ means $%
\lim\limits_{x\rightarrow \beta ^{-}}f(x)-\lim\limits_{x\rightarrow \alpha
^{+}}f(x),$ where both limits are finite. For notational convenience, let $%
b_{j}:=pa_{j}$ ($1\leq j\leq n),$ and $b_{0}$ be the zero function.
\bigskip

\noindent Our main result in this section is the following.

\bigskip

\noindent\textbf{Proposition 2.1.} \textit{With the above notation, for }$L$\textit{%
\ to be self-adjoint with respect to the inner product }$\left\langle
u,v\right\rangle =\int_{I}p(x)u(x)v(x)dx,$ \textit{it is necessary that }$%
p(x)=\dfrac{e^{\frac{2}{n}\int \frac{a_{n-1}(x)}{a_{n}(x)}dx}}{\left\vert
a_{n}(x)\right\vert }$ \textit{\ on each subinterval of }$I$ \textit{where
it is smooth} \textit{and that $n$ be even. Conversely if  $p(x)=\dfrac{e^{\frac{2}{n}\int\frac{a_{n-1}(x)}{a_{n}(x)}dx}%
}{\left\vert a_{n}(x)\right\vert }$ and }$b_{j}=pa_{j}$,
\textit{\ it is necessary and
sufficient for $L$ to be self adjoint that the following conditions hold for }$1\leq j\leq n$\textit{.}

\textit{\bigskip}

(i)\textit{\ }$\binom{n}{j}b_{n}^{(n-j)}-\binom{n-1}{j}b_{n-1}^{(n-j-1)}%
+\cdots+(-1)^{n-j}b_{j}=b_{j}$\textit{\ on }$I,$

\textit{\bigskip}

(ii)\textit{\ }$%
\begin{bmatrix}
\binom{n-1}{j-1} & -\binom{n-2}{j-1} & \cdots & (-1)^{n-j}\\
&  &  & \\
\binom{n-2}{j-2} & -\binom{n-3}{j-2} & \cdots & (-1)^{n-j}\\
&  &  & \\
\cdots & \cdots & \cdots & \cdots\\
&  &  & \\
\binom{n-j}{0} & -\binom{n-j-1}{0} & \cdots & (-1)^{n-j}%
\end{bmatrix}%
\begin{bmatrix}
b_{n}^{(n-j)}\\
\\
b_{n-1}^{(n-j-1)}\\
\cdots\\
\\
b_{j}%
\end{bmatrix}
=0$\textit{\ on }$\partial I.$

\bigskip

\noindent\textit{In
particular, }$n$ \textit{must be even}.\textit{ }

\bigskip

\bigskip

\noindent To prove Proposition 2.1 we shall need the following lemmas. The first
is a formula for repeated integration by parts, while the last one may be of
independent interest.

\bigskip

\noindent\textbf{Lemma 2.2} \textit{(See \cite{naimark}) Let }$f$\textit{\ and }%
$y$\textit{\ be functions }$k$\textit{\ times differentiable on some interval
}$I.$\textit{\ Then}$\int_{I}fy^{(k)}dx=(-1)^{k}\int_{I}f^{(k)}ydx+\left[
\sum\limits_{j=0}^{k-1}(-1)^{j}f^{(j)}y^{(k-1-j)}\right]  _{I}.$\bigskip

\noindent\textbf{Lemma 2.3} \textit{Let }$f_{j}$\textit{\ }$(1\leq j\leq
r)$\textit{\ be functions continuous on some interval }$(a,b),$\textit{\ where
}$b$\textit{\ may be infinite. If there exists a non-singular square matrix
}$A=[a_{ij}]_{1\leq i,j\leq r}$\textit{\ such that }$\lim\limits_{x\rightarrow
b^{-}}\sum\limits_{j=1}^{r}a_{ij}f_{j}(x)=0$\textit{\ for }$1\leq i\leq
r,$\textit{\ then }$\lim\limits_{x\rightarrow b^{-}}f_{j}(x)=0$\textit{\ for
}$1\leq i\leq r.$

\bigskip

\noindent\textbf{Proof.} Put $g_{i}(x)=\sum\limits_{j=1}^{r}a_{ij}f_{j}(x)$
$(1\leq i\leq r),$ so that for all $x\in(a,b)$
\[
A(f_{1}(x),f_{2}(x),...,f_{r}(x))^{T}=(g_{1}(x),g_{2}(x),...,g_{r}(x))^{T}.
\]

\noindent The conclusion follows from the fact that $\lim\limits_{x\rightarrow
b^{-}}g_{i}(x)=0$ and that
\[
(f_{1}(x),f_{2}(x),...,f_{r}(x))^{T}=A^{-1}(g_{1}(x),g_{2}(x),...,g_{r}%
(x))^{T}.\text{ }\square
\]

\bigskip

\noindent\textbf{Lemma 2.4} \textit{Let }$v_{i}$\textit{\ }$(0\leq i\leq
r)$\textit{\ be functions continuous on a real interval }$I$\textit{\ such
that for all polynomials }$u,$\textit{\ }$\left[  \sum\limits_{i=0}^{r}%
v_{i}u^{(i)}\right]  _{I}=0.$\textit{\ Then }$v_{i}=0$\textit{\ }$(0\leq i\leq
r)$\textit{\ at each endpoint of }$I.$

\bigskip

\noindent\textbf{Proof.} Suppose first that the endpoints $\alpha,\beta$
($\alpha<\beta)$ of $I$ are finite. We need only show that $v_{r}=0$ at each
endpoint of $I$, the statement for the remaining $v_{i}$ would then follow by
straightforward (reverse) induction. From $u=(x-\alpha)^{r}(x-\beta)^{r}z$
$(z$ a polynomial$)$ we get $\left[  r!v_{r}z\right]  _{I}=0$ i.e. $\left[
v_{r}z\right]  _{I}=0.$ Then, from $z=1$ and $z=x$ respectively, we get
$v_{r}(\beta)-v_{r}(\alpha)=0$ and $\beta v_{r}(\beta)-\alpha v_{r}%
(\alpha)=0.$ These equations imply $v_{r}(\beta)=v_{r}(\alpha)=0,$ as required.

\noindent Suppose now that $\beta=\infty$ (with $\alpha$ possibly infinite).
Put $u_{j}=x^{r+j}$ $(1\leq j\leq r+1).$ Then, since $\lim
\limits_{x\rightarrow\infty}x^{j}=\infty$ and $\lim\limits_{x\rightarrow
\infty}\sum_{i=0}^{r}\dfrac{(r+j)!}{(r+j-i)!}x^{r+j-i}v_{i}(x)$ is finite for
each $j$ ($1\leq j\leq r+1),$ we obtain $\sum\limits_{i=0}^{r}\dfrac
{(r+j)!x^{r-i}v_{i}(x)}{(r+j-i)!}\rightarrow0$ as $x\rightarrow\infty.$ If we
now put $f_{i}(x)=x^{r-i}v_{i}(x)$, $0\leq i\leq r,$ we get $\lim
\limits_{x\rightarrow\infty}\sum\limits_{i=0}^{r}a_{ij}f_{i}(x)=0$ for $1\leq
j\leq r+1$, where $a_{ij}=\dfrac{(r+j)!}{(r+j-i)!}.$ In view of Lemma 2.3, we
need only show that the matrix $A=[a_{ij}]_{1\leq i,j\leq r}$ is non-singular.
We have $a_{ij}=i!\binom{r+j}{i}$ hence $A$ is non-singular if and only if the
matrix $B=\left[  \binom{r+j}{i}\right]  _{0\leq i\leq r,1\leq j\leq r+1}$ is
non-singular. If, more generally, we let $D(m,n):=\det\left[  \binom
{m+i-1}{j-1}\right]  _{1\leq i,j\leq n+1}$ for $m\geq n\geq2,$ then it is easy
to show that $D(m,n)=D(m,n-1).$ This gives $D(m,n)=D(m,1)=1,$ so that $\det
B=D(r+1,r)\neq0,$ i.e. $A$ is non-singular. The case when $\alpha$ is infinite
is similarly dealt with. $\square$\bigskip

\bigskip

\noindent\textbf{Lemma 2.5} \textit{Let }$c$\textit{\ and }$c_{j}$%
\textit{\ }$(0\leq j\leq n)$\textit{\ be functions continuous on an interval
}$I.$\textit{\ If }$\left[  \sum\limits_{j=0}^{n}c_{j}y^{(j)}\right]
_{I}=\int_{I}cydx$\textit{\ for all }$n$\textit{-times differentiable
functions }$y,$\textit{\ then }$c=0$\textit{\ on }$I$\textit{\ and each
}$c_{j}$\textit{\ vanishes at each endpoint of }$I.$

\bigskip

\noindent\textbf{Proof.} We first prove that $c=0$ on $I.$ Suppose on the
contrary that $c\left(  \gamma\right)  \neq0$ for some $\gamma$ in $I.$ We can
assume that $c>0$ on some subinterval $\left[  \delta,\varepsilon\right]  $ of
$I$ containing $\gamma.$ Let
\[
\phi\left(  x\right)  =\left\{
\begin{array}
[c]{cc}%
e^{\frac{1}{\left(  x-\delta\right)  \left(  x-\varepsilon\right)  }} &
\text{if }\delta\leq x\leq\varepsilon\\
0 & \text{otherwise}%
\end{array}
\right.
\]

\noindent Then, putting $y=(x-\delta)^{2n}(x-\varepsilon)^{2n}\phi(x),$ we get
$\left[  \sum\limits_{j=0}^{n}c_{j}y^{(j)}\right]  _{I}=0$ and so $\int%
_{I}cydx=\int_{\delta}^{\varepsilon}c(x-\delta)^{2n}(x-\varepsilon)^{2n}%
\phi(x)dx=0.$ This is impossible since the integrand in this last integral is
positive. This shows that $c=0$ on $I.$ By Lemma 2.4, each $c_{j}$ equals zero
at each endpoint of $I,$ and the proof is complete. $\square$

\bigskip

\bigskip

\noindent\textbf{Proof of Proposition 2.1} By Lemma 2.2, we have for any
functions $y$ and $u$ in $C^{\infty }$%
\begin{align}
\left\langle Ly,u\right\rangle  &  =\int_{I}p(Ly)udx=\sum\limits_{k=1}^{n}%
\int_{I}pa_{k}y^{(k)}udx=\sum\limits_{k=1}^{n}\int_{I}\left(  pa_{k}u\right)
y^{(k)}dx\nonumber\\
&  =\sum\limits_{k=1}^{n}\left(  \left[  \sum\limits_{j=0}^{k-1}%
(-1)^{j}(pa_{k}u)^{(j)}y^{(k-1-j)}\right]  _{I}+(-1)^{k}\int_{I}(pa_{k}%
u)^{(k)}ydx\right)  .
\end{align}

\noindent Suppose first that $L$ is self-adjoint. Then for all polynomials $y$
and $u$, $\left\langle Ly,u\right\rangle =\left\langle y,Lu\right\rangle ,$
i.e.%
\begin{equation}
\sum\limits_{k=1}^{n}\left(  \left[  \sum\limits_{j=0}^{k-1}(-1)^{j}%
(pa_{k}u)^{(j)}y^{(k-1-j)}\right]  _{I}+(-1)^{k}\int_{I}(pa_{k}u)^{(k)}%
ydx\right)  =\sum\limits_{k=1}^{n}\int_{I}pa_{k}u^{(k)}ydx.
\end{equation}

\noindent Fix $u$ and put $c=\sum\limits_{k=1}^{n}(-1)^{k}(b_{k}u)^{(k)}%
-\sum\limits_{k=1}^{n}b_{k}u^{(k)},$ $c_{j}=\sum\limits_{k=j+1}^{n}%
(-1)^{k-1-j}\left(  b_{k}u\right)  ^{(k-1-j)}$ $(0\leq j\leq n-1).$ Then (2.2)
gives $\left[  \sum\limits_{j=0}^{n-1}c_{j}y^{(j)}\right]  _{I}=\int_{I}cydx,$
which by Lemma 2.5 implies $c=0$ on $I$ and $c_{j}=0$ $(0\leq j\leq n-1)$ at
each endpoint of $I.$ Applying Leibniz rule to the terms $\left( b_{k}u\right) ^{(k)}$ of $c$ we
obtain

\begin{equation}
\sum_{k=1}^{n}b_{k}u^{(k)}=\displaystyle\sum_{k=1}^{n}\left( -1\right)
^{k}\sum\limits_{j=0}^{k}\binom{k}{j}b_{k}^{(k-j)}u^{(j)}.  \tag{2.3}
\end{equation}

Since this is true for all $u$ in $C^{\infty }$, we can equate coefficients
of $u^{(k)}$ and get from $k=n$ that $n$ is even and that for $0\leq k\leq
n-1$

\begin{equation}
\left.
\begin{array}{c}
2b_{k}=\binom{k+1}{k}b_{k+1}^{\prime }-\binom{k+2}{k}b_{k+2}^{^{\prime
\prime }}+\cdots +\binom{n}{k}b_{n}^{\left( n-k\right) }\ \text{if }k\text{
is odd} \\
0=-\dbinom{k+1}{k}b_{k+1}^{\prime }+\binom{k+2}{k}b_{k+2}^{^{\prime \prime
}}-\cdots +\binom{n}{k}b_{n}^{\left( n-k\right) }\text{ if }k\text{ is even}%
\end{array}%
\right\}   \tag{2.4}
\end{equation}

From $k=n-1,$ we obtain the equation $2pa_{n-1}=n(pa_{n})^{\prime },$ whose
solution is $p=\dfrac{e^{\frac{2}{n}\int \frac{a_{n-1}}{a_{n}}dx}}{%
\left\vert a_{n}\right\vert }."$

\bigskip

\noindent Applying Leibniz rule again to $c_{j}=\sum\limits_{k=j+1}%
^{n}(-1)^{k-1-j}\left(  b_{k}u\right)  ^{(k-1-j)}$ $(0\leq j\leq n-1),$ we
obtain in a similar manner, but this time on $\partial I$ (i.e. at the
endpoints of $I),$ the following equations for $1\leq j\leq n$

\begin{center}%
\begin{gather*}
\binom{n-1}{j-1}b_{n}^{(n-j)}-\binom{n-2}{j-1}b_{n-1}^{(n-j-1)}+\cdots
+(-1)^{n-j}b_{j}=0\\
\binom{n-2}{j-2}b_{n}^{(n-j)}-\binom{n-3}{j-2}b_{n-1}^{(n-j-1)}+\cdots
+(-1)^{n-j}b_{j}=0\\
...........................................................\\
\binom{n-j}{0}b_{n}^{(n-j)}-\binom{n-j-1}{0}b_{n-1}^{(n-j-1)}+\cdots
+(-1)^{n-j}b_{j}=0
\end{gather*}

\end{center}

\noindent This can be put in matrix form $A_{j}%
\begin{bmatrix}
b_{n}^{(n-j)} & b_{n-1}^{(n-j-1)} & ... & b_{j}%
\end{bmatrix}
^{T}=0$ $(1\leq j\leq n),$ where $A_{j}$ is the $j\times(n-j+1)$ matrix%

\[%
\begin{bmatrix}
\binom{n-1}{j-1} & -\binom{n-2}{j-1} & \cdots & (-1)^{n-j}\\
&  &  & \\
\binom{n-2}{j-2} & -\binom{n-3}{j-2} & \cdots & (-1)^{n-j}\\
&  &  & \\
\cdots & \cdots & \cdots & \cdots\\
&  &  & \\
\binom{n-j}{0} & -\binom{n-j-1}{0} & \cdots & (-1)^{n-j}%
\end{bmatrix}
\]

\noindent We thus have for $1\leq j\leq n$%

\[
\bigskip%
\begin{bmatrix}
\binom{n}{j} & -\binom{n-1}{j} & \cdots & (-1)^{n-j}\\
&  &  & \\
\binom{n-1}{j-1} & -\binom{n-2}{j-1} & \cdots & (-1)^{n-j}\\
&  &  & \\
\cdots & \cdots & \cdots & \cdots\\
&  &  & \\
\binom{n-j}{0} & -\binom{n-j-1}{0} & \cdots & (-1)^{n-j}%
\end{bmatrix}%
\begin{bmatrix}
b_{n}^{(n-j)}\\
\\
b_{n-1}^{(n-j-1)}\\
\cdots\\
\\
b_{j}%
\end{bmatrix}
=%
\begin{bmatrix}
b_{j}\\
\\
0\\
\cdots\\
\\
0
\end{bmatrix}
,
\]
\bigskip

\noindent where the first equation is on $I$ and the remaining ones are on
$\partial I.$

\noindent Conversely, it is clear that if for $1\leq j\leq n,$ $\binom{n}%
{j}b_{n}^{(n-j)}-\binom{n-1}{j}b_{n-1}^{(n-j-1)}+\cdots+(-1)^{n-j}b_{j}=b_{j}$
on $I$ and $%
\begin{bmatrix}
\binom{n-1}{j-1} & -\binom{n-2}{j-1} & \cdots & (-1)^{n-j}\\
&  &  & \\
\binom{n-2}{j-2} & -\binom{n-3}{j-2} & \cdots & (-1)^{n-j}\\
&  &  & \\
\cdots & \cdots & \cdots & \cdots\\
&  &  & \\
\binom{n-j}{0} & -\binom{n-j-1}{0} & \cdots & (-1)^{n-j}%
\end{bmatrix}%
\begin{bmatrix}
b_{n}^{(n-j)}\\
\\
b_{n-1}^{(n-j-1)}\\
\cdots\\
\\
b_{j}%
\end{bmatrix}
=0$ on $\partial I,$ then equation (2.2) holds for all functions $y,$ $u$ in
$C^{\infty }$, and therefore $L$ is self-adjoint. $\square$

\bigskip

\bigskip

\noindent The following observation is particularly useful. When, in the above
proof, $j\geq n-j+1$ i.e. $j\geq1+n/2,$ we get more equations than "unknowns"
$b_{n}^{(n-j)},$ $b_{n-1}^{(n-j-1)},...,$ $b_{j}.$ Thus, deleting the first
$2j-n-1$ rows of $A_{j}$ and putting $k=n-j$ $(0\leq k\leq n/2-1),$ we obtain
the equations $B_{k}%
\begin{bmatrix}
b_{n}^{(k)} & b_{n-1}^{(k-1)} & ... & b_{n-k}%
\end{bmatrix}
^{T}=0$ where $B_{k}$ is the $(k+1)\times(k+1)$ matrix
\[%
\begin{bmatrix}
\binom{2k}{k} & -\binom{2k-1}{k-1} & \cdots & (-1)^{k}\\
&  &  & \\
\binom{2k-1}{k} & -\binom{2k-2}{k-1} & \cdots & (-1)^{k}\\
&  &  & \\
\cdots & \cdots & \cdots & \cdots\\
&  &  & \\
\binom{k}{k} & -\binom{k-1}{k-1} & \cdots & (-1)^{k}%
\end{bmatrix}
.
\]
Clearly, $\det B_{k}=\pm\det E_{k}$ where $E_{k}=%
\begin{bmatrix}
\binom{0}{0} & \binom{1}{0} & \cdots & \binom{k}{0}\\
&  &  & \\
\binom{1}{1} & \binom{2}{1} & \cdots & \binom{k+1}{1}\\
&  &  & \\
\cdots & \cdots & \cdots & \cdots\\
&  &  & \\
\binom{k}{k} & \binom{k+1}{k} & \cdots & \binom{2k}{k}%
\end{bmatrix}
.$ As in the proof of Lemma 2.4, if we let $E(m,k)=%
\begin{vmatrix}
\binom{m-k}{0} & \binom{m-k+1}{0} & \cdots & \binom{m}{0}\\
&  &  & \\
\binom{m-k+1}{1} & \binom{m-k+2}{1} & \cdots & \binom{m+1}{1}\\
&  &  & \\
\cdots & \cdots & \cdots & \cdots\\
&  &  & \\
\binom{m}{k} & \binom{m+1}{k} & \cdots & \binom{m+k}{k}%
\end{vmatrix}
,$ then elementary column operations give $E(m,k)=E(m,k-1)=\cdots=E(m,1)=1.$
Therefore $\det E_{k}=E(k,k)\neq0,$ i.e. $B_{k}$ is non-singular, and we
obtain%
\begin{equation}
b_{n}^{(k)}=b_{n-1}^{(k-1)}=\cdots=b_{n-k}=0\text{ on }\partial I\text{ for
}0\leq k\leq n/2-1.
\end{equation}

\noindent An interesting consequence of this is that if the weight $p=1$, then
$a_{n}^{(k)}=a_{n-1}^{(k-1)}=0$ for $0\leq k\leq n/2-1$ on $\partial I,$ and
so, if $a_{n}$ is not constant, $I$ must be a finite interval $[\alpha,\beta]$
with $a_{n}=A(x-\alpha)^{n/2}(x-\beta)^{n/2}$ for some non-zero constant $A$
and $a_{n-1}=\dfrac{An^{2}}{2}(x-(\alpha+\beta)/2)(x-\alpha)^{n/2-1}%
(x-\beta)^{n/2-1}$ (recall that the degree of $a_{k}$ is at most $k $ and that
$2pa_{n-1}=n(pa_{n})^{\prime}$)$.$ We thus obtain the form of the two leading
polynomial coefficients of what may be considered as the $n$-th order Legendre
differential equation.

\bigskip

\noindent Proposition 2.1 gives a necessary and sufficient set of conditions under which the operator $L$ is self-adjoint with respect to an inner product of the form $\left\langle u,v\right\rangle =\int_{I}puvdx.$ This was achieved under the assumption that $p$ is an admissible weight, that is $\int_{I}pudx$ is integrable for all polynomials $u$, and that $p$ satisfies certain differentiability conditions. It is therefore highly desirable that we obtain conditions under which this assumption holds. The next section is devoted to such an analysis

\section{Analytic aspects of Sturm-Liouville theory}
Keeping the same notation as before, let $\displaystyle
p=\frac{1}{|a|}e^{\frac{2}{n}\int \frac{b}{a}dx}$ be the weight
function, where, for brevity, $a=a_{n}$ is a polynomial of degree
$\leq n$ and $b=a_{n-1}$ is a polynomial degree $\leq n-1.$ Without
loss of generality, we may assume $a$ to be monic. The weight function is actually defined piecewise. If $I$ is an interval that has no zeros of the leading polynomial $a=a_{n}$ then choosing a base point $p_0$ in $I$ the weight function is given by $p(x)=\frac{e^{\frac{2}{n}\int_{p_0}^x \frac{b(t)}{a(t)}dt}}{|a(x)|}$. In this section we discuss the integrability of the weight function over an interval $I$ whose endpoints are consecutive zeros of $a(x)$. The integrability is basically a consequence of the following Lemma. \\

\noindent
\textbf{Lemma 3.1}. \textit{For $\epsilon >0$ and $d, \alpha$ integers with $\alpha >0$,  $\int_{0}^{\epsilon} \frac{e^{kx^d}}{x^\alpha}dx$ exists if and only if $k<0$ and $d<0$.   }\\

\noindent
We say that a function $f(x)$ defined and continuous on an open interval $I$ containing $0$ as a left end point is left integrable at $0$  if for any $\eta\in I$, $\displaystyle \lim_{\epsilon \to 0^+} \int_{\epsilon}^{\eta}f(x)dx$ exists. \\
Similarly if  $f(x)$ is a function defined and continuous on an
  an open interval $I$ containing $0$ as a right end point is right integrable at $0$  if for any $\eta\in I$, $\displaystyle \lim_{\epsilon \to 0^-} \int_{\eta}^{\epsilon}f(x)dx$ exists. Clearly this is equivalent to saying that the function $g(x)=f(-x)$ is left integrable at $0$. By suitable translations, one can define the concept of left and right integrability at the end points of an interval $I$ on which the given function is defined and continuous.\\
  Let $r$ be a zero of $a(x)=a_n(x)$ and let $m_a(r)=\alpha$ and $m_b(r)=\beta$ be the multiplicities of $r$ as a root of $a(x)$ and $b(x)=b_{n-1}(x)$. Thus $\frac{b(x)}{a(x)}=(x-r)^{\beta-\alpha}\phi(x)$ where $\phi(x)$ is a rational function with $\phi(r)\ne 0$. Hence $\frac{b(x)}{a(x)}=\phi(r) (x-r)^{\beta-\alpha}\psi(x)$ where $\psi(r)=1$.\\

\noindent
\textbf{Definition}. We say that a root $r$ of $a(x)$ is an ordinary root if $m_b(r)-n_a(r)+1\ne 0$, and it is a logarithmic root if $m_b(r)-m_a(r)+1 = 0$.\\

  Using Lemma 3.1 we have the following result which is one of the main tools for explicit determination of self-adjoint operators.  \\

\noindent\textbf{Proposition 3.2}. \textit{Let }$\beta=m_{b}(r),\alpha
=m_{a}(r),$ so that $\frac{b(x)}{a(x)}=(x-r)^{\beta-\alpha}\phi(x),$\textit{
where }$\phi(x)$\textit{ is a rational function with }$\displaystyle  \phi(r)=\lim
_{x\rightarrow r}(x-r)^{\alpha-\beta}\frac{b(x)}{a(x)}$ $\neq0$\textit{. Then}

\begin{enumerate}
\item[(i)] \textit{ For an ordinary root }$r$\textit{ of }$a$\textit{, the
weight function }$p(x)$\textit{ is integrable from the right at }$r$\textit{
if and only if }$\alpha-\beta\geq2$\textit{ and }$\phi(r)>0$\textit{. \newline
It is integrable from the left at }$r$\textit{ if and only if }$\alpha
-\beta\geq2$\textit{ and }$(-1)^{\alpha-\beta}\phi(r)<0$\textit{.\newline In
this case, the weight function }$p(x)$\textit{ is respectively right/left
}$C^{\infty}$\textit{ differentiable at }$r$\textit{ and }$p$\textit{ and all
its (one sided) derivatives vanish at }$r$\textit{. }

\item[(ii)] \textit{ For a logarithmic root }$r$\textit{ of }$a$\textit{, the
weight function }$p(x)$\textit{ is right/left integrable near }$r$\textit{ if
and only if }$\frac{|x-r|^{\frac{2}{n}\phi(r)}}{|x-r|^{\alpha}}$\textit{ is
integrable near }$r$\textit{ if and only if }$\dfrac{2}{n}\phi(r)-\alpha+1>0$.
\end{enumerate}

\noindent Using lower and upper bounds on the asymptotic form of the weight
function (as $x\rightarrow\infty$), or partial fraction decomposition of
$\dfrac{b}{a},$ one can prove the following result.

\bigskip

\noindent\textbf{Proposition 3.3}.
\begin{enumerate}
\item[(i)] \textit{If }$a$\textit{ has no real roots and } $p(x)=\frac{e^{\frac{2}{n}\int_{p_0}^x \frac{b(t)}{a(t)}dt}}{|a(x)|}$ \textit{ then the weight
function }$p(x)$\textit{ gives finite norm for all polynomials if and only if
}$\deg b-\deg a$ \textit{is an odd positive integer and the leading term of
}$b$\textit{ is negative. }

\item[(ii)] \textit{If }$a$\textit{ has only one root, say }$0$\textit{, then
} $p(x)=\frac{e^{\frac{2}{n}\int_{p_0}^x \frac{b(t)}{a(t)}dt}}{|a(x)|}$ \textit{
gives a finite norm for all polynomials if and only if }

\begin{description}
\item[(a)] $\deg b-\deg a\geq0$ \textit{and the leading term of} $b$
\textit{is negative.}

\item[(b)] \textit{If }$a=x^{\alpha}\left(  A_{0}+A_{1}x+\cdots\right)
$\textit{ and }$b=x^{\beta}\left(  B_{0}+B_{1}x+\cdots\right)  ,$\textit{
where }$A_{0}$\textit{ and }$B_{0}$\textit{ are nonzero constants, then
}$\alpha-\beta\geq1$\textit{, and }$\dfrac{B_{0}}{A_{0}}>0$\textit{ for
}$\alpha-\beta\geq2$\textit{ whereas }$\dfrac{2B_{0}}{nA_{0}}-\alpha
+1>0$\textit{ for }$\alpha-\beta=1.$
\end{description}
\end{enumerate}

The differentiability properties of ordinary roots have already been discussed. We now assume that the multiplicity of a root $r$ of $a(x)=a_n(x)$ is $\alpha$ and its multiplicity in $b(x)=a_{n-1}(x)$ is $\beta$. For convenience of notation we assume that $r$ is zero.

As above, we have $a_n(x)=x^\alpha (A_{0} + A_{1} x+ \cdots)$, $a_{n-1}(x)=x^\beta (B_{0} + B_{1} x+ \cdots)$ with $(\beta-\alpha)=-1$. Thus near $x=0$,
the weight is of the form $p(x)=\frac{1}{|A_{0}|}|x|^{\left(\frac{2}{n}\frac{B_{0}}{A_{0}}-\alpha\right)} \frac{e^{\phi(x)}}{1+\psi(x)}$ where $\phi$, $\psi$ are analytic near zero and $\psi(0)=0$.
When there is no danger of confusion we will write $p(x)\sim |x|^{\left(\frac{2}{n}\frac{B_{0}}{A_{0}}-\alpha\right)}$.
Now $p'=p\left(\frac{2}{n}\frac{b}{a}-\frac{a'}{a} \right)$. Therefore, all higher derivatives of $p$ are of the form $p \rho$ where $\rho$ is a rational function and all higher derivatives of $p \rho$ are also multiples of $p$ by rational functions.
For later use we record the asymptotic behavior of $p'$ near a zero of $a_n$. \\ $$p'(x)=\frac{1}{|A_{0}|}|x|^{\left(\frac{2}{n}\frac{B_{0}}{A_{0}}-\alpha\right)} \frac{e^{\phi(x)}}{1+\psi(x)} \left(\frac{2}{n}\frac{b}{a}-\frac{a'}{a} \right)$$
 The weight $p$ is integrable near zero if and only if
$\left(\frac{2}{n}\frac{B_{0}}{A_{0}}-\alpha +1\right)>0$.
Moreover $\displaystyle \lim_{x\to 0}p(x) a_n(x) =0$ if and only if $\frac{2}{n}\frac{B_{0}}{A_{0}} >0$. By the integrability of the weight $\frac{2}{n}\frac{B_{0}}{A_{0}} >\alpha -1 \geq 0$. Thus
the boundary condition $\displaystyle \lim_{x\to 0}p(x) a_n(x) =0$ is a consequence of the integrability of the weight near zero.
Similarly $p(x)a_{n-1}(x)=\frac{1}{|A_{0}|}|x|^{\left(\frac{2}{n}\frac{B_{0}}{A_{0}}-1\right)} \frac{e^{\phi(x)}}{1+\psi(x)} (B_{0} + B_{1} x + \cdots)$, keeping in mind that $\alpha-\beta = 1$. Hence $\displaystyle \lim_{x\to 0}p(x) a_{n-1}(x)=0$ if and only if
$\left(\frac{2}{n}\frac{B_{0}}{A_{0}}-1\right)>0$.
\subsection{Higher order operators}

The principal aim of this section is to prove the following result.\\

\noindent \textbf{Proposition 3.4}. \textit{Let $L=a_n(x) y^{(n)}+a_{n-1}(x) y^{(n)} +\cdots + a_2(x) y''+a_1(x) y'$ be a self-adjoint operator of order $n$ with $n>2$. If $a_n$ has a real root then the multiplicity of the root is at least $2$ and the multiplicity of the same root in $a_{n-1}$ is positive and less than its multiplicity in $a_n$. }

\noindent \textbf{Proof}.
Let $r$ be a real root of $a_n$ and assume that it is a simple root. It is then a logarithmic root. Therefore, near $r$, we have
$$
p(x)\sim |x-r|^{\left(\frac{2}{n}\frac{B_{0}}{A_{0}}-1\right)},
$$
where $a_n(x)=(x-r) (A_{0} + A_{1} (x-r)+ \cdots)$, $a_{n-1}(x)= (B_{0} + B_{1} (x-r)+ \cdots)$ and $A_{0}$, $B_{0}$ are  not zero. \\
Now $p'=p\left(\frac{2}{n}\frac{a_{n-1}}{a_n}-\frac{a_{n}'}{a_n} \right)$. Therefore
$$
p'(x)\sim|x-r|^{\left(\frac{2}{n}\frac{B_{0}}{A_{0}}-1\right)}
\left(\frac{2}{n}\frac{a_{n-1}(x)}{a_n(x)}-\frac{a_{n}'(x)}{a_n(x)} \right).
$$
Similarly
$$
p(x) a_{n-1}(x)\sim |x-r|^{\left(\frac{2}{n}\frac{B_{0}}{A_{0}}-1\right)} (B_{0} + B_{1} (x-r)+ \cdots)
$$
The boundary conditions in Proposition 2.1 imply that $(a_n p)$, $(a_{n-1} p)$ and $(a_{n-1}p)'$ vanish on the boundary. \\
Now $\displaystyle \lim_{x\to r}a_n(x)p(x)=0$ is a consequence of the integrability of the weight near $r$. Similarly
$\displaystyle \lim_{x\to r}a_{n-1}(x)p(x)=0$ if and only if $\left(\frac{2}{n}\frac{B_{0}}{A_{0}}-1\right)>0$.\\
Let
$\displaystyle l_r = \lim_{x\to r} (x-r)^{\alpha-\beta}  \frac{a_{n-1}(x)}{a_n(x)}$. Clearly $l_r=\frac{B_{0}}{A_{0}}=\frac{a_{n-1}(r)}{a_{n}'(r)}$, as $\alpha-\beta=1$.\\
Since $\displaystyle \lim_{x\to r}pa_{n-1}=0$ and $a_{n-1}(r)\ne 0$ we see that $p$ must vanish at $r$ in the sense that its limit at $r$ is zero. The boundary condition $\displaystyle \lim_{x\to r}(a_{n-1}p)'=0$ now implies that $\displaystyle \lim_{x\to r}p'(x)=0$.
Now $p'=p\left(\frac{2}{n}\frac{a_{n-1}}{a_n}-\frac{a'_n}{a_n}\right)$. Thus near the root $r$,
$$
p'\sim |x-r|^{\left(\frac{2}{n}l_r-2\alpha\right)}\left(\frac{2}{n}a_{n-1}-a'_n  \right).$$\\
If $a_{n-1}-\frac{n}{2}a'_n\equiv 0$ then in particular $\left(\frac{2}{n}\frac{a_{n-1}(r)}{a'_n(r)}-1\right)=0$.\\
This means that
$\displaystyle \lim_{x\to r}(x-r)\frac{2}{n}\frac{a_{n-1}(x)}{a_n(x)}-1=0$ i.e. $\left(\frac{2}{n}\frac{B_{0}}{A_{0}}-1\right)=0$. This contradicts the boundary condition $\displaystyle \lim_{x\to r}a_{n-1}(x)p(x)=0$.\\

Let $(a_{n-1} - \frac{n}{2}a'_n)=(x-r)^{\lambda} H(x)$ where $\lambda \geq 0$ and $H(x)\ne 0$. If $(a_{n-1} - \frac{n}{2}a'_n)=(x-r)^{\lambda} H(x)$ with $\lambda > 0$ then $\displaystyle \lim_{x\to r}(x-r)\frac{2}{n}\frac{a_{n-1}(x)}{a_n(x)}-1=0$ which again contradicts the boundary condition $\displaystyle \lim_{x\to r}a_{n-1}(x)p(x)=0$.\\
Hence $p'\sim |x-r|^{\left(\frac{2}{n}l_r-2\alpha\right)}H(x)$ so $p' \to 0$ at $r$ if and only if $\left(\frac{2}{n}l_r-2\alpha\right)>0$.

By Proposition 2.1, the operator must satisfy - beside other equations - the determining equations
\begin{equation}\label{det1-gen}
n(a_{n}p)^{\prime}=2(a_{n-1}p)
\end{equation}
\begin{equation}\label{det2-gen}
\dfrac{\left(  n-1\right)  (n-2)}{6}(a_{n-1}p)^{\prime\prime}-\left(
n-2\right)  (a_{n-2}p)^{\prime}+2(a_{n-3}p)=0
\end{equation}
Equation (\ref{det2-gen}) is equivalent to
\begin{equation}\label{rational-gen}
C_1\left(\frac{a_{n-1}}{a_n}\right)^3+C_2\left(\frac{a_{n-1}}{a_n}\right)\left(\frac{a_{n-1}}{a_n}\right)'+
C_3\left(\frac{a_{n-1}}{a_n}\right)''+C_4 \left(\frac{a_{n-2}}{a_n}\right)' + C_5 \frac{a_{n-1}}{a_n}\frac{a_{n-2}}{a_n}+ C_6\frac{a_{n-3}}{a_n}=0
\end{equation}
where
$$
C_1=\frac{2(n-1)(n-2)}{3n^2},\ C_2=\frac{(n-1)(n-2)}{n},
$$
$$
C_3=\frac{(n-1)(n-2)}{6},\ C_4= -(n-2),\ C_5=-\frac{2(n-2)}{n}, \
C_6=2.
$$
This implies the identity
\begin{equation}\label{divide-general}
a_{n-1}\left(a_{n-1}-2\frac{n}{2}a'_n\right)\left(a_{n-1}-\frac{n}{2}a'_n\right)\equiv 0 ({\rm mod }a_n)
\end{equation}
Using this identity, as $(x-r)$ divides $a_n$ but it does not divide $a_{n-1}$ nor $\left(a_{n-1}-\frac{n}{2}a'_n\right)$, it must
divide $\left(a_{n-1}-2\frac{n}{2}a'_n\right)$. But then
$\displaystyle \lim_{x\to r}\left(a_{n-1}-2\frac{n}{2}a'_n\right)=0$. This means that
$\displaystyle \lim_{x\to r}\frac{2}{n}\frac{a_{n-1}}{a'_n}-2=0$ i.e. $\frac{2}{n}\frac{B_{0}}{A_{0}}-2=0$. As
seen above $p' \to 0$ at $r$ if and only if $\left(\frac{2}{n}l_r-2\alpha\right)>0$. Since $\alpha =1$ we have a contradiction.

Therefore $a_n$ cannot have a simple root and its multiplicity $\alpha$ in $a_n$ is at least $2$. Suppose that the multiplicity $\beta$ of $r$ in $a_{n-1}$ is zero. By considering the order of poles of $a_n$ in (\ref{rational-gen}) we see that $\beta$ cannot be zero.  This completes the proof of the proposition. $\square$\\

This result has an important consequence for fourth order self-adjoint operators.
\textbf{Corollary 3.5.} \textit{Let $L$ be a self-adjoint operator of order $4$ and $a_4$ be its leading term. If $a_4$ has more than one real root then it has exactly two real roots with multiplicity $2$. Moreover the multiplicity of each real root of $a_4$ in $a_{3}$ is $1$.
}

\bigskip
We also have the following result.\\
\noindent\textbf{Proposition 3.6}. \textit{Let }$n>2$\textit{ and suppose that
}$a=a_{n}$\textit{ has at most one real root. Then }$2\deg b-\deg a\leq
n-2\ \ $\textit{or \ }$3\deg b-2\deg a\leq n-3,$ where $b=a_{n-1}.$

\textit{\noindent}(i)\textit{ If }$a$\textit{ has no real root then }%
\[
\deg a<\deg b\leq n-3
\]

\textit{\noindent}(ii)\textit{ Suppose }$a$\textit{ has only one real root
}$r$\textit{ with multiplicity }$\alpha,$\textit{ let }$\beta$\textit{ be the
multiplicity of }$r$\textit{ as a root of }$b,$\textit{ and let }$a=\left(
x-r\right)  ^{\alpha}u,$\textit{ }$b=\left(  x-r\right)  ^{\beta}v.$\textit{
Then }%
\[
2\leq\deg a\leq\deg b\leq n-2\mathit{\ \ }\text{and}\mathit{\ \ }1+\deg
u\leq\deg v\leq n-3
\]

\noindent \textbf{Proof}. First, in all cases, $\deg b\geq1.$ This is because if $a$
has no real root then $\deg b$ is odd and if $a$ has (at least) one real root
then this will also be a root for $b$ (by Proposition 3.4). If we multiply by
$a^{3}$ both sides of (3.3), then the six terms on the left-hand side will be
polynomials with respective degrees
\begin{align*}
& 3\deg b,\ \ 2\deg b+\deg a-1,\ \ 2\deg a+\deg b-2,\ \ \deg a+\deg
a_{n-2}-1,\\
& \deg a+\deg b+\deg a_{n-2},\ \ 2\deg a+\deg a_{n-3}%
\end{align*}
\noindent A comparison of degrees shows that $a_{n-2}$ and $a_{n-3}$ cannot be
both zero and that%
\begin{align*}
2\deg b  & \leq\deg a+\deg a_{n-2},\ \text{or}\\
3\deg b  & \leq2\deg a+\deg a_{n-3}%
\end{align*}
\noindent Using the fact that $\deg a_{j}\leq j,$ we obtain
\[
2\deg b-\deg a\leq n-2\ \ \text{or }\ 3\deg b-2\deg a\leq n-3
\]

\noindent If $a$ has no real roots then, by Proposition 3.3 (i), $\deg b-\deg
a\geq1$ and hence%
\[
\deg a<\deg b\leq n-3
\]
\noindent If $a$ has only one real root $r$ with multiplicity $\alpha,$ then
$\alpha\geq2$ and $b$ has $r$ as a root with multiplicity $\beta,$ where
$1\leq\beta<\alpha$ (by Proposition 3.4). Since $\deg b\geq\deg a,$ we obtain
that
\[
2\leq\deg a\leq\deg b\leq n-2
\]
\noindent Let $a=\left(  x-r\right)  ^{\alpha}u,$ $b=\left(  x-r\right)
^{\beta}v.$ Then $\deg a=\alpha+\deg u$ and $\deg b=\beta+\deg v,$ and we
obtain $\deg v-\deg u\geq\alpha-\beta\geq1.$ Now $\beta+\deg v\leq n-2,$ so
$1+\deg u\leq\deg v\leq n-3$ and thus $\deg u\leq n-4.$ $\square$

\bigskip


%
%

\section{Examples of higher order operators, their eigenvalues and orthogonal eigenfunctions}

 Let$L$ be an operator of the form $L(y)=\sum\limits_{k=0}^{n}a_{k}(x)y^{(k)},$ where $\deg a_k\leq k$; then the eigenvalues of $L$ are the coefficients of $x^n$ in $L(x^n)$, $n=0,1,2,\cdots$.\\

\noindent \textbf{Proposition 4.1.} \textit{Let $L$ be a linear operator that maps the space $\mathbb{P}_{n}$ of all polynomials of degree at most $n$ into to itself for all $n\leq N$.
If the eigenvalues of $L$ are distinct or if $L$ is a self-adjoint operator then there is an eigenpolynomial of $L$ in every degree $\leq N$.  }\\

The proof is left to the reader. This means that if two operators leave the space of polynomials of degree at most $n$ invariant for all $n$ and the weight function which makes the two operators self-adjoint is the same, then they have the same eigenfunctions. The eigenvalues in general are not simple.

Let $\lambda$ be an eigenvalue of $L$ and $\mathbb{P}_{n}(\lambda)$ the corresponding eigenspace in the space $\mathbb{P}_{n}$ of all polynomials of degree$\leq n$.

If $n_0$ is the minimal degree in $\mathbb{P}_{n}(\lambda)$, then there is, up to a scalar only one polynomial in $\mathbb{P}_{n}(\lambda)$  of degree $n_0$. Choose a monic  polynomial $Q_1$ in $\mathbb{P}_{n_0}(\lambda)$.

Let $n_1$ be the smallest degree, if any, greater than $n_0$ in  $\mathbb{P}_{n}(\lambda)$. The codimension of $\mathbb{P}_{n_0}(\lambda)$ in $\mathbb{P}_{n}(\lambda)$ is $1$. Therefore, in the orthogonal complement of $\mathbb{P}_{n_0}(\lambda)$ in $\mathbb{P}_{n_1}(\lambda)$, choosing a monic polynomial $Q_2$,which will necessarily be of degree $n_1$, the polynomials $Q_1,\ Q_2$   give an orthogonal basis of $\mathbb{P}_{n_1}(\lambda)$. Continuing this process, we eventually get  an orthogonal basis of $\mathbb{P}_{n}(\lambda)$ consisting of monic polynomials.

We illustrate this by an example of a fourth order self-adjoint operator that has repeated eigenvalues but which is not an iterate of a second order operator. \\
Consider the  operator
\begin{equation}\label{op-special-eg}
L=(1-x^2)^2 y^{(4)}-8x(1-x^2) y''' + 8 y''-24x y'.
\end{equation}
Its eigenvalues are $\lambda_n = n[ (n-1)(n-2)(n+5)-24]$. The eigenvalue $\lambda = -24$ is repeated in degrees $n=1$ and $n=3$.   \\
The weight function for which this operator is self-adjoint is  $p(x)=1.$  The eigenpolynomials of degree at most $3$ are
$$
y_0(x)=1,\ \ y_1(x)=x,\ \ y_2(x)=x^2-\frac{1}{3}, \ \ y_3(x)=x^3.
$$
This gives the set of orthogonal polynomials $\{1,x,x^2-\frac{1}{3},x^3-\frac{3}{5}x\}$. Since the weight function is the same as that for the classical Legendre polynomials, this family up to scalars is the same as the corresponding classical Legendre polynomials.

We now return to examples of higher order operators. The restrictions on the parameters appearing in all the examples come from integrability of the weight and boundary conditions.
Before giving examples of higher order operators,  it is instructive to consider the classical case of second order operators in the frame work of section 3.

\subsection{Self-adjoint operators of order $2$}
Assume $n=2$  and that $a_2  (x)$ has distinct roots, which we may assume to be $-1$ and $1$. If $\alpha$ is the multiplicity of a root $r$ of  $a_2  (x)$ and $\beta$ is its multiplicity in $a_1 (x)$ then the integrability of the associated weight   gives the equation $\alpha=\beta + 1 + \delta$, with $\delta \geq 0$. As $\alpha=1$, we must have $\beta=0$, $\delta=0$. Thus, only the logarithmic case can occur.

Let $a_1(x)=  cx+d$. The integrability condition at a root $r$ reads $\displaystyle \lim_{x\to r}(x-r)\frac{a_1(x)}{a_2(x)}>0$.
As we are assuming that $a_2  (x)=x^2-1$, the integrability conditions at both the roots gives $c+d>0$, $-c+d<0$, so $-c<d<c$.

If  $a_2  (x)$ has no real roots, then, by Proposition 3.3, $a_2  (x)$ must be a constant, so taking $a_2  (x)=1$, we have $a_1(x) = cx+d$, with $c<0$.
Finally if  $a_2  (x)$ has only one real root, we may take it to be $0$. By Proposition 3.3, only the logarithmic case can occur, so taking $a_2  (x)=x$, we have $a_1(x) = cx+d$, with $c<0$, and $d>0$.

\subsection{Self-adjoint operators of order $4$}
In this section we determine all self-adjoint operators
\begin{equation}\label{main-op-n-4}
L=a_4(x) y^{(4)}+a_3(x) y''' + a_2(x) y''+a_1(x) y',
\end{equation}
with an admissible weight $p(x)=\dfrac{e^{\frac{1}{2}\int \frac{a_3(x)}{a_4(x)}dx}%
}{\left\vert a_4(x)\right\vert }$,
satisfying the differential equation
\begin{equation}\label{main-eq-n-4}
L(y)=\lambda y.
\end{equation}
By Proposition 2.1, the operator $L$ must satisfy the determining equations
\begin{eqnarray}
& (a_4p)' = \fr{1}{2} (a_3 p) \label{det1-n-4} \\
& (a_3p)'' - 2(a_2p)' + 2(a_1p)=0 \label{det2-n-4}
\end{eqnarray}
on $I$, subject to the vanishing of
\begin{equation}\label{boundary-n-4}
(a_4 p),\quad (a_3 p) \quad {\rm and} \quad \l((a_2p) -
\fr{(a_3p)'}{2}\r)
\end{equation}
 on the boundary $\partial I$.
\subsubsection{Case: $a_4(x)$ has no real roots}
If $a_4(x)$ is a monic polynomial having no real roots then $I = (-\infty,\infty)$ and, by Proposition 3.6,  we have $a_4(x)=1$, $a_3(x)$ linear.

Considering $a_4(x)=1$, $a_3(x)$ linear and solving determining equations (\ref{det1-n-4}), (\ref{det2-n-4}) subject to the constraints (\ref{boundary-n-4}) determines the fourth order self-adjoint operators (\ref{main-op-n-4}) and the differential equations (\ref{main-eq-n-4}) as
\[
\begin{split}
& a_4(x)=1, \quad a_3(x) = 2 (m_1-2 m_2^2 x), \\
& a_2(x) = 4  m_2^4 x^2-4 m_1  m_2^2 x+{A}, \quad
a_1(x) = (-m_1^2+2 m_2^2+{A}) (m_1-2 {m_2}^2 x) \end{split}\]
with the weight function
$$
p(x)= e^{-m_2^2 x^2 +{m_1}x +  {m_0}} \qquad (m_2 \ne 0)
$$
and the eigenvalues
$$
\lambda_n = 2 m_2^2 \left(m_1^2-{A}+2 m_2^2 (n-2)\right) n .
$$
Depending on the choice of $A, m_1, m_2$ one can get repeated eigenvalues. \\
The specific case $m_1=0, m_2=1, A=-4$ gives the standard fourth order Hermite operator \cite{littlejohn-krall89}
$$
a_4(x)=1, \quad a_3(x) =-4x, \quad a_2(x)=4 (x^2-1), \quad a_1(x)=4x
$$
with the weight $p(x)= e^{-x^2}$ and non-repeated eigenvalues $\lambda_n = 4n^2$.
It is worth noticing that this specific case is the only operator in the class of fourth order operators with $a_4(x)=1$ that is an iterate of the second order Hermite operator, and in general, this class is not obtained as iteration of the second order case.
\subsubsection{Case: $a_4(x)$ has only one real root}
In this case we have, by Proposition 3.6, $a_4(x)=x^2$ and $a_3(x)=x(a+bx)$ with $b\ne 0$. The weight is determined as $p(x) = |x|^{\fr{a}{2}-2} e^{\fr{b}{2}x}$ with $b<0$ and $I = (0,\infty)$. \\
Solving the determining equation  (\ref{det2-n-4}) subject to the constraints (\ref{boundary-n-4}) determines the fourth order self-adjoint operators (\ref{main-op-n-4}) and the differential equations (\ref{main-eq-n-4}) as
\[
\begin{split}
& p(x) = |x|^{\fr{a}{2}-2} e^{\fr{b}{2}x} \mb{ with } a>2,
b<0 \\
& a_4(x) = x^2, \quad a_3(x) = x(a+bx),\\
& a_2(x)  = \fr{1}{4}[-2a+a^2+x(4A+ b^2 x)], \quad a_1(x)  = \fr{1}{4}( 2A-ab)(-2+a+bx) \end{split}\]
with the eigenvalues
$$
\lambda_n= \frac{b}{4}  n (2 {A}-a b-b+b n),
$$
which in general are not iterates of second order operator.\\
The special case $a=4,b=-2, A=-5$ gives the classical fourth order Laguerre operator \cite{littlejohn-krall89} with the weight $p(x)=e^{-x}$ as
\[
\begin{split}
& a_4(x) = x^2, \quad a_3(x) = -2(-2+x)x, \quad a_2(x) = x^2-5 x+2, \quad  a_1(x) = -1+x
\end{split}
\]
and the eigenvalues $\lambda_n= n^2$. This coincides with the second iteration of the classical second order Laguerre operator

\subsubsection{Case: $a_4(x)$ has more that one real root}
By Corollary 3.5, $a_4(x)$ must have exactly two real roots with multiplicity $2$ and the multiplicity of each real root of $a_4$ in $a_{3}$ is $1$. By a linear change of variables and scaling, we may assume that the roots are $-1$ and $1$, and take $a_4(x)=(1-x^2)^2$.\\
Now  $a_3(x)  = K(1-x^2)$ is ruled out by Proposition 3.2 (ii).\\
So $a_3(x)$ must have the form $(k_1+k_2 x)(1-x^2)$ with $k_2\ne 0$. Without loss of generality, we consider the form $a_3(x) = -2(b+(-2+a)x)(-1+x^2)$ which determines the weight function.
$$
p(x) = \fr{(1+x)^{\fr{b-a-2}{2}}}{(1-x)^{\fr{b+a+2}{2}}} \mb{ with } b-a > 0 \mb{ and } b+a < 0.
$$
\begin{enumerate}
\item If $b=0$, then $a_3(x)=2 (-2+a) x(1-x^2)$ with $a < 0$.\\
Solving the determining equation  (\ref{det2-n-4}) subject to the constraints (\ref{boundary-n-4}) determines the fourth order self-adjoint operators (\ref{main-op-n-4}) and the differential equations (\ref{main-eq-n-4}) as
\[
\begin{split}
 & p(x)  = (1-x^2)^{-1-\fr{a}{2}} \qquad (a < 0)\\
 & a_4(x)  = (1-x^2)^2, \quad  a_3(x)  = -2 (-2+a) x(-1+x^2) \\
 & a_2(x) = -2a+a^2+A(-1+x^2), \quad a_1(x)  = a(2-3a+a^2-A) x
\end{split}\]
with the eigenvalues
$$
\lambda_n = n (-a+n-1) (-a^2-n a+4 a+n^2+{A}-n-2).
$$
The particular cases that are iterates of corresponding second order operators are given below:
\begin{enumerate}
\item The values $a=-2$, $A=14$ lead to the Legendre operator \cite{littlejohn-krall89} with the weight $p(x)=1$ as
$$
a_4(x)  = (1-x^2)^2,\ \ a_3(x)  = -8x (1-x^2),\ \ a_2(x) = 14 x^2-6 ,\ \ a_1(x)  =4x
$$
with the eigenvalues
$
\lambda_n = n^2 (n+1)^2.
$
\item The special case $a=-1$, $A=7$ is the Chebychev operator of first kind \cite{littlejohn-krall89} with the weight $p(x)=\frac{1}{\sqrt{1-x^2}}$ as
$$
a_4(x)  = (1-x^2)^2,\ \ a_3(x)  = -6x (1-x^2),\ \ a_2(x) = 7 x^2-4 ,\ \ a_1(x)  =x
$$
with the eigenvalues
$
\lambda_n = n^4.
$
\item The special case $a=-3$, $A=23$ is the Chebychev operator of second kind \cite{littlejohn-krall89} with the weight $p(x)=\sqrt{1-x^2}$ as
$$
a_4(x)  = (1-x^2)^2,\ \ a_3(x)  = -10x (1-x^2),\ \ a_2(x) = 23 x^2-8 ,\ \ a_1(x)  =9x
$$
with the eigenvalues
$
\lambda_n = n^2 (n+2)^2.
$
\end{enumerate}
\item If $b\ne 0$, then $a_3(x) = -2(b+(-2+a)x)(-1+x^2) \mb{ with } b-a > 0 \mb{ and } b+a < 0.$ So $a<b<-a$.\\
Solving determining equation  (\ref{det2-n-4}) subject to the constraints (\ref{boundary-n-4}) determines the following fourth order self-adjoint operators (\ref{main-op-n-4}) and the differential equations (\ref{main-eq-n-4}) for this case.
\[
\begin{split}
& p(x) = (1-x)^{\fr{1}{2} (-2-a-b)} (1+x)^{\fr{1}{2}(-2-a+b)} \quad (b-a > 0, \ b+a<0 \mb{ and } b\ne  0) \\
& a_4(x) = (1-x^2)^2 \\
& a_3(x) = -2(b+(-2+a)x)(-1+x^2) \\
& a_2(x) =
\fr{b^3+B+2(-1+a)b^2x-Bx^2+b(-2+a+2x^2-3ax^2+a^2x^2)}{b} \\
& a_1(x)= B + \fr{aB x}{b} \end{split} \]
with the eigenvalues
$$
\lambda_n = \frac{(a-n+1) n \left(-b n^2+a b n+b n-a b+{B}\right)}{b}.
$$
\end{enumerate}

%
%
\subsection{Self-adjoint operators of order $6$}
This section provides examples of the self-adjoint operators
\begin{equation}\label{main-op-n-6}
L=a_6(x) y^{(6)}+a_5(x) y^{(5)}+a_4(x) y^{(4)}+a_3(x) y''' + a_2(x) y''+a_1(x) y'
\end{equation}
with an admissible weight $p(x)=\dfrac{e^{\frac{1}{3}\int \frac{a_5(x)}{a_6(x)}dx}%
}{\left\vert a_6(x)\right\vert }$,
satisfying the differential equation
\begin{equation}\label{main-eq-n-6}
L(y)=\lambda y.
\end{equation}
By Proposition 2.1, the determining equations in this case are
\begin{eqnarray}
& 3(a_6p)' =  (a_5 p) \label{det1-n-6} \\
& 5(a_5p)'' - 6(a_4p)' + 3(a_3p)=0 \label{det2-n-6n}\\
& (a_4p)'''-3(a_3p)'' +5(a_2p)' - 5(a_1p)=0 \label{det3-n-6n}
\end{eqnarray}
on $I$, subject to the vanishing of
\begin{equation}\label{boundary-n-6}
(a_6 p),\quad (a_5 p),   \quad (a_5p)', \quad (a_4p), \quad  (a_4p)' - 3(a_3p), \quad
(a_4p)''-3(a_3p)' +5(a_2p)
\end{equation}
 on the boundary $\partial I$.\\
Equations (\ref{det2-n-6n}) and (\ref{det3-n-6n}) are equivalent to
\begin{equation}
\frac{10}{27}\left(\frac{a_{5}}{a_6}\right)^3+\frac{10}{3}\left(\frac{a_{5}}{a_6}\right)\left(\frac{a_{5}}{a_6}\right)'+
\frac{10}{3}\left(\frac{a_{5}}{a_6}\right)''-4 \left(\frac{a_{4}}{a_6}\right)' -\frac{4}{3} \frac{a_{5}}{a_6}\frac{a_{4}}{a_6}+ 2\frac{a_{3}}{a_6}=0
\end{equation}
and
\begin{eqnarray}
&-&81\left(\frac{{a_3} }{{a_6} }\right)''+\frac{{a_4} }{{a_6} } \left(\frac{{a_5} }{{a_6} }\right)^3-9 \frac{{a_3} }{{a_6} } \left(\frac{{a_5} }{{a_6} }\right)^2+45 \frac{{a_2} }{{a_6} } \left(\frac{{a_5} }{{a_6} }\right)-135 \frac{{a_1} }{{a_6} }+135 \left(\frac{{a_2} }{{a_6} }\right)'+27 \left(\frac{{a_4} }{{a_6} }\right)''' \nonumber \\
&+&9 \frac{{a_4} }{{a_6} } \left(\frac{{a_5} }{{a_6} }\right)''+27 \frac{{a_5} }{{a_6} } \left(\frac{{a_4} }{{a_6} }\right)''
-54 \left(\frac{{a_5} }{{a_6} }\right) \left(\frac{{a_3} }{{a_6} }\right)'+9 \left(\frac{{a_5} }{{a_6} }\right)^2 \left(\frac{{a_4} }{{a_6} }\right)'
+9 \frac{{a_4} }{{a_6} } \frac{{a_5} }{{a_6} } \left(\frac{{a_5} }{{a_6} }\right)' \nonumber \\
&-&27 \left(\frac{{a_3} }{{a_6} }\right) \left(\frac{{a_5} }{{a_6} }\right)'+27 \left(\frac{{a_4} }{{a_6} }\right)' \left(\frac{{a_5} }{{a_6} }\right)'=0.
\end{eqnarray}
These identities give the congruences
\begin{equation}\label{divide-n-6-1}
5 {a_5}  \left({a_5} -6 {a_6}' \right) \left({a_5} -3 {a_6}' \right)\equiv 0 ({\rm mod }a_6)
\end{equation}
and
\begin{equation}\label{divide-n-6-2}
{a_4}  \left( {a_5} -9  {a_6}' \right) \left( {a_5} -6  {a_6}' \right) \left( {a_5} -3  {a_6}' \right)\equiv 0 ({\rm mod }a_6).
\end{equation}
Before giving examples in detail, we note that, in case the leading term has a real root,there are many operators that satisfy the determining equations but for which one of the boundary condition fails at one or both
end points of interval $I$. In case the leading term has no real roots, the boundary conditions will be satisfied - because of the form of the weight - but where one of the determining equations will not be satisfied. Here are some typical examples.\\

\noindent {\bf Example 4.1}\\
The operator in (\ref{main-op-n-6}) with
\[
\begin{split}
 p(x) & = (x-1)^2 (x+1), \\
 a_6 (x)& = (x-1)^2 (x+1)^4, \quad a_5 (x)  = 3 (x-1) (x+1)^3 (9 x-1),\\
 a_4(x) & = 60 x (x+1)^2 (5 x-3), \quad   a_3(x)  = 240 \left(7 x^3+6 x^2-2 x-1\right), \\
 a_2(x) & = 720 x (5 x+3), \quad a_1(x) = 360 (x+1)
\end{split} \]
satisfies the determining equations (\ref{det1-n-6}), (\ref{det2-n-6n}), (\ref{det3-n-6n}) and all the boundary conditons in (\ref{boundary-n-6}) except that the last boundary condition fails at the end point $1$ of $I=[-1,1]$.
\\

 \noindent {\bf Example 4.2}\\
The operator in (\ref{main-op-n-6}) with
\[
\begin{split}
 p(x) & = e^{-m^2 x} x^2 \quad (m\ne 0) \\
 a_6 (x)& = x^2, \quad a_5 (x)  = -3 x (m^2 x-4), \\
 a_4(x) & = -30  m^2 x+{A} x^2+30, \quad a_3(x)  = 5 x^2 m^6-2 ({A} x^2+30 ) m^2+8 {A} x, \\
 a_2(x) & = x ({C}-3 m^8 x)+{A} (m^4 x^2 +12),  \\
 a_1(x) & = 18  m^8 x-60 m^6-8 {A}  m^4 x+24 {A} m^2+{C} (3-m^2 x)
\end{split} \]
satisfies the determining equations (\ref{det1-n-6}), (\ref{det2-n-6n}), (\ref{det3-n-6n}) and all the boundary conditons in (\ref{boundary-n-6}) except that the last boundary condition fails at the end point $0$ of $I = [0,\infty)$.\\

\noindent {\bf Example 4.3}\\
The operator in (\ref{main-op-n-6}) with
\[
\begin{split}
 p(x) & = \frac{e^{-m^2 x^2}}{x^2+1} \quad (m\ne 0) \\
 a_6 (x)& = x^2+1, \quad a_5 (x)  = -6 m^2 x (x^2+1), \\
 a_4(x) & = 10 x^4 m^4-10 m^4+{A} x^2+{A}, \quad a_3(x)  = -4 m^2 (-10 m^4+5 m^2+{A}) x (x^2+1), \\
 a_2(x) & = {C_2} x^2+{C_1} x+{C_0},  \quad a_1(x)  = {D_0}+{D_1} x
\end{split} \]
satisfies the determining equations (\ref{det1-n-6}), (\ref{det2-n-6n}) and all the boundary conditons in (\ref{boundary-n-6}) for $I = (-\infty,\infty)$ but fails to satisfy the remaining determining equation (\ref{det3-n-6n}).\\

In the following sections,  we  give examples of sixth order self-adjoint operators with the weights of the form $p(x)=e^{-x^2}$, $p(x)=|x|^{n} e^{m x}$ and $p(x)=\frac{(1+x)^m}{(1-x)^n}$ by solving  the determining equations and boundary conditions using Mathematica.
\subsubsection{Sixth order operators (\ref{main-op-n-6}) with $a_6=1$ and weight of the form $e^{-x^2}$}
\[
\begin{split}
& p(x)  = e^{-x^2} \\
& a_6(x)   = 1, a_5(x) = -6x, a_4(x) = 12x^2 - \fr{D}{4} - \fr{C}{2} - 6,\\
& a_3(x)   = -8x^3 + (D + 2C + 12)x, a_2(x) =  (-D - 2C)x^2+C, a_1(x)  = D x
\end{split}\]
with eigenvalues
$$
\lambda_n = -n \left(8 n^2+2 {C} n+{D} n-24 n-2 {C}-2 {D}+16\right)  .
$$
The specific case $D=-8$, $C=16$ gives the standard sixth order Hermite operator
\begin{eqnarray*}
&& a_6(x)=1,\quad a_5(x)= -6x,\quad   a_4(x)=12 x^2-12,\\
&& a_3(x) =36 x-8 x^3, \quad a_2(x)=16-24 x^2,  \quad a_1(x)=-8x
\end{eqnarray*}
with the eigenvalues $\lambda_n = -8 n^3$, and coincides with third iteration of the classical second order Hermite operator.
%
\subsubsection{Sixth order operators (\ref{main-op-n-6}) with $a_6=x^3$ and weight of the form $|x|^{n} e^{m x}$}
\[
\begin{split}
& p(x)  = |x|^{\fr{a}{3}-3} e^{\fr{b}{3}x} \mb{ with } a>6, b<0 \\
& a_6(x) = x^3, \quad a_5(x) = x^2(a+bx), \\
& a_4(x)  = \fr{1}{3} x(-3a+a^2+x(3A+b^2 x)) \\
& a_3(x) = \fr{1}{27} (a^3-9a^2(1+bx)+9a(2+2Ax+3bx-b^2x^2)+x(b^3x^2+18A(-3+bx)))\\
& a_2(x) = \fr{1}{27} (18a^2b-2a^3b+27Cx-2ab(18+b^2x^2)+3A(18-9a+a^2+b^2x^2))\\
& a_1(x) = \fr{1}{81} (-6(-3+a)Ab-12ab^2+4a^2b^2+27C)(-6+a+bx)
\end{split}\]
with eigenvalues
$$
\lambda_n = \frac{b}{81}  n \left(4 a^2 b^2+3 n^2 b^2-6 a b^2-6 a n b^2-9 n b^2+6 b^2-6 a {A} b+9 {A} b+9 {A} n b+27 {C}\right)  .
$$
The special case $a=9, b = -3, A=-21, C=13$ is the standard sixth order Laguerre operator with the weight $p(x) = e^{-x}$ as
\[
\begin{split}
 p(x) & = e^{-x} \\
 a_6(x) &  = x^3, \ \  a_5(x)   = -3(-3+x)x^2, \ \ a_4(x)  = 3 x (x^2-7 x+6), \\
 a_3(x)  &  = -x^3+15 x^2-30 x+6, \ \ a_2(x)    = -3 x^2+13 x-6, \ \ a_1(x)    = 1-x
\end{split}\]
with eigenvalues
$
\lambda_n=n^3,
$
which coincides with the third iteration of the classical second order Laguerre operator.
%
\subsubsection{Sixth order operators (\ref{main-op-n-6}) with $a_6=(1-x^2)^3 $ and weight of the form $\fr{(1+x)^{\fr{b-a-2}{2}}}{(1-x)^{\fr{b+a+2}{2}}}$ with $(b-a-2)>4 \mbox{ or } (b-a-2) \in \left\{0,2,4 \right\}$, and $(b+a+2)<-4 \mbox{ or } (b+a+2) \in \left\{-4,-2,0 \right\}$
}

\noindent {\bf 4.3.3.1 The case $J^{6.I}:$\quad $(b-a-2) \in
\left\{0,2,4 \right\},\ (b+a+2) \in \left\{-4,-2,0 \right\} $}
\begin{itemize}
\item $J^{6.I.a}:$ \ $b-a-2=0,\ b+a+2=0$

\[
\begin{split}
 & p(x)  =1  \\
 & a_6 (x) ={\left( 1 - x^2 \right) }^3  ,\\
 & a_5 (x)  =-18\,x\,{\left( -1 + x^2 \right) }^2  ,\\
 & a_4(x)  =\frac{-\left( \left( -1 + x^2 \right) \,\left( 2\,C\,\left( -1 + x^2 \right)  + D\,\left( -1 + x^2 \right)  +
        144\,\left( -1 + 5\,x^2 \right)  \right)  \right) }{8}  ,\\
 & a_3(x)  =x\,\left( 72 + D - 120\,x^2 - D\,x^2 - 2\,C\,\left( -1 + x^2 \right)  \right)  ,\\
 & a_2(x)  =C + \left( -3\,C - D \right) \,x^2  ,\\
 & a_1(x)  =D\,x
\end{split} \]
with eigenvalues
\begin{small}
$$
\lambda_n =  -\frac{1}{8} n (n+1) \left(8 n^4+16 n^3+2 C n^2+D n^2-56 n^2+2 C n+D n-64 n-4 C-6 D+96\right) .
$$
\end{small}
The special case $C=36, D=-8$ gives the standard sixth order Legendre operator with the weight $p(x) = 1$ as
\[
\begin{split}
 a_6(x) &  = {\left( 1 - x^2 \right) }^3, \ \  a_5(x)   =-18\,x \left( 1-{x}^{2} \right) ^{2} , \ \ a_4(x)  = 2\, \left( 1-{x}^{2} \right)  \left( 49\,{x}^{2}-13 \right), \\
 a_3(x)  &  =-8\,x \left( -17+23\,{x}^{2} \right) , \ \ a_2(x)    = 36-100\,{x}^{2}, \ \ a_1(x)    = -8\,x
\end{split}\]
with eigenvalues
$
\lambda_n=-n^3 (n+1)^3,
$
which coincides with the third iteration of the classical second order Legendre operator.
\item $J^{6.I.b}:$ \ $b-a-2=0,\ b+a+2=-2$
\begin{small}
\[
\begin{split}
 & p(x)  = 1-x \\
 & a_6 (x) ={\left( 1 - x^2 \right) }^3  ,\\
 & a_5 (x)  =-3\,\left( 1 + 7\,x \right) \,{\left( -1 + x^2 \right) }^2  ,\\
 & a_4(x)  =\frac{-\left( \left( -1 + x^2 \right) \,\left( C\,\left( -1 + x^2 \right)  + D\,\left( -1 + x^2 \right)  +
        72\,\left( -1 + 2\,x + 7\,x^2 \right)  \right)  \right) }{4}  ,\\
 & a_3(x)  =\frac{12\,\left( 3 + 15x - 15x^2 - 35x^3 \right)  - C\,\left( -1 - 5x + x^2 + 5x^3 \right)  -
    D\,\left( -1 - 5x + x^2 + 5x^3 \right) }{2}  ,\\
 & a_2(x)  =C + \left( -2\,C - 2\,D \right) \,x + \left( \frac{-5\,C}{3} + \frac{5\,\left( -2\,C - 2\,D \right) }{3} - \frac{2\,D}{3} \right)
     \,x^2  ,\\
 & a_1(x)  =D + 3\,D\,x
\end{split} \]
\end{small}
with eigenvalues
\begin{small}
$$
\lambda_n = -\frac{1}{4} n (n+2) \left(4 n^4+16 n^3+C n^2+D n^2-28 n^2+2 C n+2 D n-88 n-3 C-7 D+96\right)  .
$$
\end{small}
\item $J^{6.I.c}:$ \ $b-a-2=0,\ b+a+2=-4$
\begin{small}
\[
\begin{split}
 & p(x)  = (1-x)^2 \\
 & a_6 (x) ={\left( 1 - x^2 \right) }^3  ,\\
 & a_5 (x)  =-6\,\left( 1 + 4\,x \right) \,{\left( -1 + x^2 \right) }^2  ,\\
 & a_4(x)  =\frac{-\left( \left( -1 + x^2 \right) \,\left( 2\,C\,\left( -1 + x^2 \right)  + D\,\left( -1 + x^2 \right)  +
        48\,\left( -1 + 7\,x + 14\,x^2 \right)  \right)  \right) }{4}  ,\\
 & a_3(x)  =12\left( 3 + 6\,x - 21\,x^2 - 28\,x^3 \right)  - 2C\left( -1 - 3\,x + x^2 + 3\,x^3 \right)  -
  D\left( -1 - 3\,x + x^2 + 3\,x^3 \right)  ,\\
 & a_2(x)  =C + \left( -10\,C - 5\,D \right) \,x + \left( -15\,C - 7\,D \right) \,x^2  ,\\
 & a_1(x)  =D + 2\,D\,x
\end{split} \]
\end{small}
with eigenvalues
\begin{small}
$$
\lambda_n = -\frac{1}{4} n (n+3) \left(4 n^4+24 n^3+2 C n^2+D n^2-20 n^2+6 C n+3 D n-168 n-8 C-6 D+160\right)  .
$$
\end{small}
\item $J^{6.I.d}:$ \ $b-a-2=2,\ b+a+2=0$
\begin{small}
\[
\begin{split}
 & p(x)  = 1+x  \\
 & a_6 (x) ={\left( 1 - x^2 \right) }^3  ,\\
 & a_5 (x)  =-3\,\left( -1 + 7\,x \right) \,{\left( -1 + x^2 \right) }^2  ,\\
 & a_4(x)  =\frac{\left( -1 + x^2 \right) \,\left( C - C\,x^2 + 72\,\left( 1 + 2\,x - 7\,x^2 \right)  + D\,\left( -1 + x^2 \right)  \right) }
  {4}  ,\\
 & a_3(x)  =\frac{12\,\left( -3 + 15\,x + 15\,x^2 - 35\,x^3 \right)  + C\,\left( -1 + 5\,x + x^2 - 5\,x^3 \right)  +
    D\,\left( 1 - 5\,x - x^2 + 5\,x^3 \right) }{2}  ,\\
 & a_2(x)  =2\,D\,x\,\left( -1 + 2\,x \right)  + C\,\left( 1 + 2\,x - 5\,x^2 \right)  ,\\
 & a_1(x)  =D - 3\,D\,x
\end{split} \]
\end{small}
with eigenvalues
\begin{small}
$$
\lambda_n = -\frac{1}{4} n (n+2) \left(4 n^4+16 n^3+C n^2-D n^2-28 n^2+2 C n-2 D n-88 n-3 C+7 D+96\right)  .
$$
\end{small}
\item $J^{6.I.e}:$ \ $b-a-2=2,\ b+a+2=-2$
\[
\begin{split}
 & p(x)  =\left( 1 - x \right) \,\left( 1 + x \right)  \\
 & a_6 (x) ={\left( 1 - x^2 \right) }^3  ,\\
 & a_5 (x)  =-24\,x\,{\left( -1 + x^2 \right) }^2  ,\\
 & a_4(x)  =\frac{-\left( \left( -1 + x^2 \right) \,\left( 4\,C\,\left( -1 + x^2 \right)  + D\,\left( -1 + x^2 \right)  +
        576\,\left( -1 + 7\,x^2 \right)  \right)  \right) }{24}  ,\\
 & a_3(x)  =\frac{x\,\left( D - D\,x^2 + 96\,\left( 3 - 7\,x^2 \right)  - 4\,C\,\left( -1 + x^2 \right)  \right) }{2}  ,\\
 & a_2(x)  =C - \left( 5\,C + D \right) \,x^2  ,\\
 & a_1(x)  = D\,x
\end{split} \]
with eigenvalues
\begin{footnotesize}
$$
\lambda_n =  -\frac{1}{24} n (n+3) \left(24 n^4+144 n^3+4 C n^2+D n^2-120 n^2+12 C n+3 D n-1008 n-16 C-10 D+960\right).
$$
\end{footnotesize}
\item $J^{6.I.f}:$ \ $b-a-2=2,\ b+a+2=-4$
\begin{small}
\[
\begin{split}
 & p(x)  ={\left( 1 - x \right) }^2\,\left( 1 + x \right)  \\
 & a_6 (x) ={\left( 1 - x^2 \right) }^3  ,\\
 & a_5 (x)  =-3\,\left( 1 + 9\,x \right) \,{\left( -1 + x^2 \right) }^2  ,\\
 & a_4(x)  =\frac{-\left( \left( -1 + x^2 \right) \,\left( C\,\left( -1 + x^2 \right)  + D\,\left( -1 + x^2 \right)  +
        144\,\left( -1 + 2\,x + 9\,x^2 \right)  \right)  \right) }{6}  ,\\
 & a_3(x)  =\frac{-\left( C\left( -1 - 7x + x^2 + 7x^3 \right)  \right)  - D\left( -1 - 7x + x^2 + 7x^3 \right)  -
    72\left( -1 - 7x + 7x^2 + 21x^3 \right) }{3}  ,\\
 & a_2(x)  =-2\,D\,x\,\left( 1 + 3\,x \right)  + C\,\left( 1 - 2\,x - 7\,x^2 \right)  ,\\
 & a_1(x)  =D + 5 D x
\end{split} \]
\end{small}
with eigenvalues
\begin{small}
$$
\lambda_n = -\frac{1}{6} n (n+4) \left(6 n^4+48 n^3+C n^2+D n^2-6 n^2+4 C n+4 D n-408 n-5 C-11 D+360\right)  .
$$
\end{small}
\item $J^{6.I.g}:$ \ $b-a-2=4,\ b+a+2=0$
\begin{small}
\[
\begin{split}
 & p(x)  ={\left( 1 + x \right) }^2  \\
 & a_6 (x) ={\left( 1 - x^2 \right) }^3  ,\\
 & a_5 (x)  =-6\,\left( -1 + 4\,x \right) \,{\left( -1 + x^2 \right) }^2  ,\\
 & a_4(x)  =\frac{\left( -1 + x^2 \right) \,\left( 48\,\left( 1 + 7\,x - 14\,x^2 \right)  - 2\,C\,\left( -1 + x^2 \right)  +
      D\,\left( -1 + x^2 \right)  \right) }{4}  ,\\
 & a_3(x)  =12\left( -3 + 6\,x + 21\,x^2 - 28\,x^3 \right)  + C\left( -2 + 6\,x + 2\,x^2 - 6\,x^3 \right)  +
  D\left( 1 - 3\,x - x^2 + 3\,x^3 \right),\\
 & a_2(x)  =C + \left( 10\,C - 5\,D \right) \,x + \left( -15\,C + 7\,D \right) \,x^2  ,\\
 & a_1(x)  =D - 2\,D\,x
\end{split} \]
\end{small}
with eigenvalues
\begin{small}
$$
\lambda_n = -\frac{1}{4} n (n+3) \left(4 n^4+24 n^3+2 C n^2-D n^2-20 n^2+6 C n-3 D n-168 n-8 C+6 D+160\right)  .
$$
\end{small}
\item $J^{6.I.h}:$ \ $b-a-2=4,\ b+a+2=-2$
\begin{small}
\[
\begin{split}
 & p(x)  =\left( 1 - x \right) \,{\left( 1 + x \right) }^2  \\
 & a_6 (x) ={\left( 1 - x^2 \right) }^3 ,\\
 & a_5 (x)  =-3\,\left( -1 + 9\,x \right) \,{\left( -1 + x^2 \right) }^2  ,\\
 & a_4(x)  =\frac{\left( -1 + x^2 \right) \,\left( C - C\,x^2 + D\,\left( -1 + x^2 \right)  - 144\,\left( -1 - 2\,x + 9\,x^2 \right)  \right)
      }{6}  ,\\
 & a_3(x)  =\frac{C\,\left( -1 + 7\,x + x^2 - 7\,x^3 \right)  + D\,\left( 1 - 7\,x - x^2 + 7\,x^3 \right)  -
    72\,\left( 1 - 7\,x - 7\,x^2 + 21\,x^3 \right) }{3}  ,\\
 & a_2(x)  =C + \left( 2\,C - 2\,D \right) \,x + \left( -7\,C + 6\,D \right) \,x^2  ,\\
 & a_1(x)  = D - 5\,D\,x
\end{split} \]
\end{small}
with eigenvalues
\begin{small}
$$
\lambda_n = -\frac{1}{6} n (n+4) \left(6 n^4+48 n^3+C n^2-D n^2-6 n^2+4 C n-4 D n-408 n-5 C+11 D+360\right)  .
$$
\end{small}
\item $J^{6.I.i}:$ \ $b-a-2=4,\ b+a+2=-4$
\[
\begin{split}
 & p(x)  ={\left( 1 - x \right) }^2\,{\left( 1 + x \right) }^2  \\
 & a_6 (x) ={\left( 1 - x^2 \right) }^3  ,\\
 & a_5 (x)  =-30\,x\,{\left( -1 + x^2 \right) }^2  ,\\
 & a_4(x)  =\frac{-\left( \left( -1 + x^2 \right) \,\left( 6\,C\,\left( -1 + x^2 \right)  + D\,\left( -1 + x^2 \right)  +
        1440\,\left( -1 + 9\,x^2 \right)  \right)  \right) }{48}  ,\\
 & a_3(x)  =\frac{x\,\left( 720 + 6\,C + D - \left( 2160 + 6\,C + D \right) \,x^2 \right) }{3}  ,\\
 & a_2(x)  =C - \left( 7\,C + D \right) \,x^2  ,\\
 & a_1(x)  =D\,x
\end{split} \]
with eigenvalues
\begin{footnotesize}
$$
\lambda_n = -\frac{1}{48} n (n+5) \left(48 n^4+480 n^3+6 C n^2+D n^2+240 n^2+30 C n+5 D n-4800 n-36 C-14 D+4032\right).
$$
\end{footnotesize}
\end{itemize}
%

\noindent {\bf 4.3.3.2 The case $J^{6.II}:$\quad $(b-a-2) \in
\left\{0,2,4 \right\},\ (b+a+2)<-4 $}
\begin{itemize}
\item $J^{6.II.a}:$ \ $b-a-2=0,\ (b+a+2)<-4$
\[
\begin{split}
 & p(x)  ={\left( 1 - x \right) }^{-2 - a}  \\
 & a_6 (x) ={\left( 1 - x^2 \right) }^3  ,\\
 & a_5 (x)  =3\,\left( 2 + a - 4\,x + a\,x \right) \,{\left( -1 + x^2 \right) }^2  ,\\
 & a_4(x)  = \frac{-\left( -1 + x^2 \right) }{2\,\left( -2 + a + a^2 \right) } X_4  \mbox{ with }\\
 & \quad X_4 = 6\,a^4\,{\left( 1 + x \right) }^2 - \left( -1 + x \right) \,\left( C + 144\,x + C\,x \right)  -
  36\,a^3\,\left( -1 + x^2 \right) \\
  & \quad + 18\,a^2\,\left( 1 - 6\,x + x^2 \right)  + 12\,a\,\left( -5 - 4\,x + 13\,x^2
  \right) ,\\
 & a_3(x)  = \frac{X_3}{-2 + a + a^2}  \mbox{ with }\\
 & \quad X_3 = a^5\,{\left( 1 + x \right) }^3 - 2\,a^4\,{\left( 1 + x \right) }^2\,\left( -5 + 4\,x \right)  +
  2\,{\left( -1 + x \right) }^2\,\left( 12 + C + 24\,x + C\,x \right) \\
  & \quad   + 3\,a^3\,\left( 5 - 9\,x - 9\,x^2 + 5\,x^3 \right)  +
  a^2\,\left( -22 - 48\,x + 42\,x^2 + 20\,x^3 \right) \\
  & \quad  +
  a\,\left( -28 + 60\,x + 60\,x^2 - 76\,x^3 - C\,\left( -1 + x \right) \,{\left( 1 + x \right) }^2
  \right), \\
 & a_2(x)  = \frac{X_2}{2\,a\,\left( -2 + a + a^2 \right) } \mbox{ with } \\
 & \quad X_2 = a^3\,C\,{\left( 1 + x \right) }^2 + 4\,D\,\left( -1 + x^2 \right)  +
  2\,a\,\left( D + C\,{\left( -1 + x \right) }^2 - D\,x^2 \right)\\
  & \quad   -
  a^2\,\left( 1 + x \right) \,\left( 2\,D\,\left( -1 + x \right)  + C\,\left( -5 + 3\,x \right)
  \right), \\
 & a_1(x)  =\frac{D\,\left( 2 + a + a\,x \right) }{a}
\end{split} \]
with eigenvalues
$$
\lambda_n =  -\frac{n (-a+n-1) }{2 (a-1) a (a+2)} \Lambda_n
$$
where
\begin{small}
\begin{eqnarray*}
\Lambda_n &=& 2 n^2 a^5-6 n a^5+4 a^5-4 n^3 a^4+12 n^2 a^4-8 n a^4+2 n^4 a^3-8 n^3 a^3+4 n^2 a^3\\
&& +14 n a^3 -12 a^3+2 n^4 a^2+4 n^3 a^2-22 n^2 a^2-C a^2+2 D a^2+C n
   a^2+8 n a^2  \\
&& +8 a^2-4 n^4 a +8 n^3 a-C n^2 a+4 n^2 a+2 D a+C n a-8 n a-4 D
\end{eqnarray*}
\end{small}
%
The zeros of denominator terms in above expressions consist of
$a=-2, a=0$ and $a=1$. Since $(b-a-2)=0$ and $(b+a+2)<-4 $ imply
$a<-4$, the denominator terms in above expressions never vanish.
\item $J^{6.II.b}:$ \ $b-a-2=2,\ (b+a+2)<-4$
\[
\begin{split}
 & p(x)  ={\left( 1 - x \right) }^{-3 - a}\,\left( 1 + x \right)  \\
 & a_6 (x) ={\left( 1 - x^2 \right) }^3  ,\\
 & a_5 (x)  =3\,\left( 4 + a - 4\,x + a\,x \right) \,{\left( -1 + x^2 \right) }^2  ,\\
 & a_4(x)  =-\left( \frac{-1 + x^2}{a\,\left( 2 - 3\,a + a^2 \right) } \right) X_4  \mbox{ with } \\
 & \quad X_4 = D - D\,x^2 + 3\,a^5\,{\left( 1 + x \right) }^2 - a\,\left( -1 + x \right) \,\left( 72 + C - 72\,x + C\,x \right)\\
 & \quad   -
  6\,a^4\,\left( -3 + 2\,x + 5\,x^2 \right)  - 6\,a^2\,\left( 9 - 38\,x + 25\,x^2 \right)  +
  3\,a^3\,\left( -13 - 26\,x + 35\,x^2 \right), \\
 & a_3(x)  = \frac{X_3}{a\,\left( 2 - 3\,a + a^2 \right) }  \mbox{ with }\\
 & \quad X_3 = -12\,a^5\,\left( -1 + x \right) \,{\left( 1 + x \right) }^2 + a^6\,{\left( 1 + x \right) }^3 +
  4\,D\,\left( 2 - x - 2\,x^2 + x^3 \right)\\
  & \quad   - 24\,a^3\,\left( 4 - x - 8\,x^2 + 5\,x^3 \right)  +
  a^4\,\left( 7 - 75\,x - 27\,x^2 + 55\,x^3 \right)\\
  & \quad   +
  2\,a\,\left( -1 + x \right) \,\left( -24\,{\left( -1 + x \right) }^2 - D\,{\left( 1 + x \right) }^2 +
     2\,C\,\left( -2 - x + x^2 \right)  \right)\\
     & \quad   - 2\,a^2\,
   \left( C\,\left( -1 + x \right) \,{\left( 1 + x \right) }^2 - 2\,\left( 7 + 45\,x - 75\,x^2 + 31\,x^3 \right)
   \right), \\
 & a_2(x)  =\frac{X_2}{a\,\left( 2 - 3\,a + a^2 \right) }  \mbox{ with }\\
 & \quad X_2 = -8\,D\,\left( -2 + x \right)  + a^2\,\left( 2\,D - 3\,C\,\left( -3 + x \right)  \right) \,\left( 1 + x \right)  +
  a^3\,C\,{\left( 1 + x \right) }^2\\
  & \quad  + 2\,a\,\left( 3\,D\,\left( 1 + x \right)  + C\,\left( 7 - 4\,x + x^2 \right)  \right), \\
 & a_1(x)  = \frac{D\,\left( 4 + a + a\,x \right) }{a}
\end{split} \]
with eigenvalues
$$
\lambda_n =  -\frac{n (-a+n-1) }{(a-2) (a-1) a} \Lambda_n
$$
where
\begin{small}
\begin{eqnarray*}
\Lambda_n &=&
 n^2 a^5-3 n a^5+2 a^5-2 n^3 a^4+2 n^2 a^4+8 n a^4-8 a^4+n^4 a^3+4 n^3 a^3-14 n^2 a^3-n a^3\\
 && +10 a^3-3 n^4 a^2+2 n^3 a^2+13 n^2 a^2-C a^2+D a^2+C n a^2-8 n
   a^2-4 a^2+2 n^4 a\\
   && -4 n^3 a-C n^2 a-2 n^2 a-4 D a+C n a+D n a+4 n a-D n^2+2 D+D n
\end{eqnarray*}
\end{small}
The zeros of denominator terms in above expressions consist of $a=0,
a=1$ and $a=2$. Since $(b-a-2)=2$ and $(b+a+2)<-4 $ imply $a<-5$,
the denominator terms in above expressions never vanish.
\item $J^{6.II.c}:$ \ $b-a-2=4,\ (b+a+2)<-4$
\[
\begin{split}
 & p(x)  ={\left( 1 - x \right) }^{-4 - a}\,{\left( 1 + x \right) }^2  \\
 & a_6 (x) ={\left( 1 - x^2 \right) }^3  ,\\
 & a_5 (x)  =3\,\left( 6 + a - 4\,x + a\,x \right) \,{\left( -1 + x^2 \right) }^2  ,\\
 & a_4(x)  =-\left( \frac{-1 + x^2}{a\,\left( 2 - 3\,a + a^2 \right) } \right) X_4  \mbox{ with }\\
 & \quad X_4 = D - D\,x^2 + 3\,a^5\,{\left( 1 + x \right) }^2 - 30\,a^4\,\left( -1 + x^2 \right)  -
  30\,a^2\,\left( 7 - 12\,x + 5\,x^2 \right)\\
  & \quad  + 15\,a^3\,\left( -1 - 10\,x + 7\,x^2 \right)  +
  a\,\left( C - C\,x^2 + 24\,\left( 8 - 9\,x + 3\,x^2 \right)
  \right) ,\\
 & a_3(x)  = \frac{X_3}{a\,\left( 2 - 3\,a + a^2 \right) } \mbox{ with }\\
 & \quad X_3 = a^6\,{\left( 1 + x \right) }^3 - 6\,a^5\,{\left( 1 + x \right) }^2\,\left( -3 + 2\,x \right)  +
  4\,D\,\left( 3 - x - 3\,x^2 + x^3 \right)\\
  & \quad  - 30\,a^3\,\left( 5 + 6\,x - 11\,x^2 + 4\,x^3 \right)  +
  a^4\,\left( 55 - 75\,x - 75\,x^2 + 55\,x^3 \right)\\
  & \quad  -
  2\,a^2\,\left( 118 - 306\,x + 234\,x^2 - 62\,x^3 + C\,\left( -1 + x \right) \,{\left( 1 + x \right) }^2 \right) \\
  & \quad  +
  2a\left( -\left( D\left( -1 + x \right) {\left( 1 + x \right) }^2 \right)  +
     2C\left( 3 - x - 3\,x^2 + x^3 \right)  - 12\left( -13 + 16\,x - 9\,x^2 + 2\,x^3 \right)
     \right),\\
 & a_2(x)  =\frac{X_2}{a\,\left( 2 - 3\,a + a^2 \right) }  \mbox{ with }\\
 & \quad X_2 = -12\,D\,\left( -3 + x \right)  + a^3\,C\,{\left( 1 + x \right) }^2 -
  a^2\,\left( 1 + x \right) \,\left( -2\,D + C\,\left( -13 + 3\,x \right)  \right)\\
  & \quad  +
  2\,a\,\left( 5\,D\,\left( 1 + x \right)  + C\,\left( 17 - 6\,x + x^2 \right)
  \right),\\
 & a_1(x)  = \frac{D\,\left( 6 + a + a\,x \right) }{a}
\end{split} \]
with eigenvalues
$$
\lambda_n = -\frac{n (-a+n-1) }{(a-2) (a-1) a} \Lambda_n
$$
where
\begin{small}
\begin{eqnarray*}
\Lambda_n &=& n^2 a^5-3 n a^5+2 a^5-2 n^3 a^4+2 n^2 a^4+8 n a^4-8 a^4+n^4 a^3+4 n^3 a^3-14 n^2 a^3-n a^3\\
&& +10 a^3-3 n^4 a^2+2 n^3 a^2+13 n^2 a^2-C a^2+D a^2+C n a^2-8 n
   a^2-4 a^2+2 n^4 a-4 n^3 a\\
   && -C n^2 a-2 n^2 a-4 D a+C n a+D n a+4 n a-D n^2+2 D+D n
\end{eqnarray*}
\end{small}
%
The zeros of denominator terms in above expressions consist of $a=0,
a=1$ and $a=2$. Since $(b-a-2)=4$ and $(b+a+2)<-4 $ imply $a<-6$,
the denominator terms in above expressions never vanish.
\end{itemize}
\newpage
%

\noindent {\bf 4.3.3.3 The case $J^{6.III}:$\quad $(b-a-2)>4, \
(b+a+2) \in \left\{-4,-2,0 \right\} $}
\begin{itemize}
\item $J^{6.III.a}:$ \ $(b-a-2)>4,\ b+a+2=0$
\[
\begin{split}
 & p(x)  ={\left( 1 + x \right) }^{-2 - a}  \\
 & a_6 (x) ={\left( 1 - x^2 \right) }^3  ,\\
 & a_5 (x)  = 3\,\left( -2 + a\,\left( -1 + x \right)  - 4\,x \right) \,{\left( -1 + x^2 \right) }^2 ,\\
 & a_4(x)  =\frac{-\left( -1 + x^2 \right) }{2\,\left( -2 + a + a^2 \right) } X_4  \mbox{ with }\\
 & \quad X_4 = 6\,a^4\,{\left( -1 + x \right) }^2 + \left( C\,\left( -1 + x \right)  - 144\,x \right) \,\left( 1 + x \right)  -
  36\,a^3\,\left( -1 + x^2 \right) \\
  & \quad + 18\,a^2\,\left( 1 + 6\,x + x^2 \right)  + 12\,a\,\left( -5 + 4\,x + 13\,x^2 \right) ,\\
 & a_3(x)  = -\left( \frac{X_3}{-2 + a + a^2} \right)  \mbox{ with }\\
& \quad X_3 = -\left( C\,\left( -1 + x^2 \right) \,\left( a\,\left(
-1 + x \right)  - 2\,\left( 1 + x \right)  \right)  \right) \\
& \quad +
  \left( -2 - a \right) \,\left( 1 - a \right) \,\left( -\left( a^3\,{\left( -1 + x \right) }^3 \right)  +
     9\,a^2\,{\left( -1 + x \right) }^2\,\left( 1 + x \right) \right)\\
     & \quad  +  \left( -2 - a \right) \,\left( 1 - a \right)
     \left(  12\,{\left( 1 + x \right) }^2\,\left( -1 + 2\,x \right)  +
     a\,\left( 8 + 30\,x - 12\,x^2 - 26\,x^3 \right)  \right)  ,\\
 & a_2(x)  =\frac{X_2}{2\,a\,\left( -2 + a + a^2 \right) }  \mbox{ with }\\
 & \quad X_2 =-\left( a^3\,C\,{\left( -1 + x \right) }^2 \right)  + 4\,D\,\left( -1 + x^2 \right)  -
  2\,a\,\left( 1 + x \right) \,\left( D\,\left( -1 + x \right)  + C\,\left( 1 + x \right)  \right) \\
  & \quad +
  a^2\,\left( -1 + x \right) \,\left( -2\,D\,\left( 1 + x \right)  + C\,\left( 5 + 3\,x \right)  \right)    ,\\
 & a_1(x)  = \frac{D\,\left( -2 + a\,\left( -1 + x \right)  \right) }{a}
\end{split} \]
with eigenvalues
$$
\lambda_n =  -\frac{n (-a+n-1) }{2 (a-1) a (a+2)} \Lambda_n
$$
where
\begin{small}
\begin{eqnarray*}
\Lambda_n &=& 2 n^2 a^5-6 n a^5+4 a^5-4 n^3 a^4+12 n^2 a^4-8 n a^4+2 n^4 a^3-8 n^3 a^3+4 n^2 a^3\\
&&+14 n a^3-12 a^3+2 n^4 a^2+4 n^3 a^2-22 n^2 a^2+C a^2+2 D a^2-C n
   a^2+8 n a^2\\
   && +8 a^2-4 n^4 a+8 n^3 a+C n^2 a+4 n^2 a+2 D a-C n a-8 n a-4 D
\end{eqnarray*}
\end{small}
%
The zeros of denominator terms in above expressions consist of
$a=-2, a=0$ and $a=1$. Since $(b-a-2)>4$ and $(b+a+2)=0 $ imply
$a<-4$, the denominator terms in above expressions never vanish.
\item $J^{6.III.b}:$ \ $(b-a-2)>4,\ b+a+2=-2$
\[
\begin{split}
 & p(x)  =\left( 1 - x \right) \,{\left( 1 + x \right) }^{-3 - a}  \\
 & a_6 (x) ={\left( 1 - x^2 \right) }^3  ,\\
 & a_5 (x)  =3\,{\left( -1 + x^2 \right) }^2\,\left( a\,\left( -1 + x \right)  - 4\,\left( 1 + x \right)  \right)  ,\\
 & a_4(x)  =-\left( \frac{-1 + x^2}{a\,\left( 2 - 3\,a + a^2 \right) } \right) X_4  \mbox{ with }\\
 & \quad X_4 = D + 3\,a^5\,{\left( -1 + x \right) }^2 - D\,x^2 - 6\,a^4\,\left( -3 - 2\,x + 5\,x^2 \right)  -
  6\,a^2\,\left( 9 + 38\,x + 25\,x^2 \right)\\
  & \quad  + 3\,a^3\,\left( -13 + 26\,x + 35\,x^2 \right)  -
  a\,\left( 1 + x \right) \,\left( C\,\left( -1 + x \right)  - 72\,\left( 1 + x \right)  \right) ,\\
 & a_3(x)  = \frac{X_3}{a\,\left( 2 - 3\,a + a^2 \right) }  \mbox{ with }\\
 & \quad X_3 = a^6\,{\left( -1 + x \right) }^3 - 12\,a^5\,{\left( -1 + x \right) }^2\,\left( 1 + x \right)  +
  4\,D\,\left( -2 - x + 2\,x^2 + x^3 \right) \\
  & \quad - 24\,a^3\,\left( -4 - x + 8\,x^2 + 5\,x^3 \right)  +
  a^4\,\left( -7 - 75\,x + 27\,x^2 + 55\,x^3 \right)\\
  & \quad  -
  2\,a\,\left( 1 + x \right) \,\left( D\,{\left( -1 + x \right) }^2 + 24\,{\left( 1 + x \right) }^2 -
     2\,C\,\left( -2 + x + x^2 \right)  \right) \\
     & \quad - 2\,a^2\,
   \left( C\,{\left( -1 + x \right) }^2\,\left( 1 + x \right)  - 2\,\left( -7 + 45\,x + 75\,x^2 + 31\,x^3 \right)  \right)  ,\\
 & a_2(x)  = \frac{X_2}{a\,\left( 2 - 3\,a + a^2 \right) }  \mbox{ with }\\
 & \quad X_2 = a^3\,C\,{\left( -1 + x \right) }^2 + 8\,D\,\left( 2 + x \right)  -
  a^2\,\left( -1 + x \right) \,\left( 2\,D + 3\,C\,\left( 3 + x \right)  \right) \\
  & \quad +
  2\,a\,\left( -3\,D\,\left( -1 + x \right)  + C\,\left( 7 + 4\,x + x^2 \right)  \right)  ,\\
 & a_1(x)  = \frac{D\,\left( -4 + a\,\left( -1 + x \right)  \right) }{a}
\end{split} \]
with eigenvalues
$$
\lambda_n = -\frac{n (-a+n-1) }{(a-2) (a-1) a}  \Lambda_n
$$
where
\begin{small}
\begin{eqnarray*}
\Lambda_n &=&  n^2 a^5-3 n a^5+2 a^5-2 n^3 a^4+2 n^2 a^4+8 n a^4-8 a^4+n^4 a^3+4 n^3 a^3-14 n^2 a^3-n a^3\\
&& +10 a^3-3 n^4 a^2+2 n^3 a^2+13 n^2 a^2-C a^2+D a^2+C n a^2-8 n
   a^2-4 a^2+2 n^4 a\\
   &&-4 n^3 a-C n^2 a-2 n^2 a-4 D a+C n a+D n a+4 n a-D n^2+2 D+D n
\end{eqnarray*}
\end{small}
%
The zeros of denominator terms in above expressions consist of $a=0,
a=1$ and $a=2$. Since $(b-a-2)>4$ and $(b+a+2)=-2 $ imply $a<-5$,
the denominator terms in above expressions never vanish.
\item $J^{6.III.c}:$ \ $(b-a-2)>4,\ b+a+2=-4$
\[
\begin{split}
 & p(x)  ={\left( 1 - x \right) }^2\,{\left( 1 + x \right) }^{-4 - a}  \\
 & a_6 (x) ={\left( 1 - x^2 \right) }^3  ,\\
 & a_5 (x)  =3\,\left( -6 + a\,\left( -1 + x \right)  - 4\,x \right) \,{\left( -1 + x^2 \right) }^2  ,\\
 & a_4(x)  = -\left( \frac{-1 + x^2}{a\,\left( 2 - 3\,a + a^2 \right) } \right) X_4  \mbox{ with }\\
 & \quad X_4 =D + 3\,a^5\,{\left( -1 + x \right) }^2 - D\,x^2 - 30\,a^4\,\left( -1 + x^2 \right)  -
  30\,a^2\,\left( 7 + 12\,x + 5\,x^2 \right) \\
  & \quad + 15\,a^3\,\left( -1 + 10\,x + 7\,x^2 \right)  +
  a\,\left( C - C\,x^2 + 24\,\left( 8 + 9\,x + 3\,x^2 \right)  \right)    ,\\
 & a_3(x)  = \frac{X_3}{a\,\left( 2 - 3\,a + a^2 \right) }  \mbox{ with }\\
 & \quad X_3 =a^6\,{\left( -1 + x \right) }^3 - 6\,a^5\,{\left( -1 + x \right) }^2\,\left( 3 + 2\,x \right)  +
  4\,D\,\left( -3 - x + 3\,x^2 + x^3 \right) \\
  & \quad - 30\,a^3\,\left( -5 + 6\,x + 11\,x^2 + 4\,x^3 \right)  +
  a^4\,\left( -55 - 75\,x + 75\,x^2 + 55\,x^3 \right) \\
  & \quad +
  2\,a^2\,\left( 118 + 306\,x + 234\,x^2 + 62\,x^3 - C\,{\left( -1 + x \right) }^2\,\left( 1 + x \right)  \right)\\
  & \quad  -
  2\,a\,\left( D\,{\left( -1 + x \right) }^2\,\left( 1 + x \right)  - 2\,C\,\left( -3 - x + 3\,x^2 + x^3 \right)  +
     12\,\left( 13 + 16\,x + 9\,x^2 + 2\,x^3 \right)  \right)   ,\\
 & a_2(x)  = \frac{X_2}{a\,\left( 2 - 3\,a + a^2 \right) }  \mbox{ with }\\
 & \quad X_2 = a^3\,C\,{\left( -1 + x \right) }^2 + 12\,D\,\left( 3 + x \right)  -
  a^2\,\left( -1 + x \right) \,\left( 2\,D + C\,\left( 13 + 3\,x \right)  \right) \\
  & \quad +
  2\,a\,\left( -5\,D\,\left( -1 + x \right)  + C\,\left( 17 + 6\,x + x^2 \right)  \right)   ,\\
 & a_1(x)  = \frac{D\,\left( -6 + a\,\left( -1 + x \right)  \right) }{a}
\end{split} \]
with eigenvalues
$$
\lambda_n = -\frac{n (-a+n-1) }{(a-2) (a-1) a}  \Lambda_n
$$
where
\begin{small}
\begin{eqnarray*}
\Lambda_n &=&n^2 a^5-3 n a^5+2 a^5-2 n^3 a^4+2 n^2 a^4+8 n a^4-8 a^4+n^4 a^3+4 n^3 a^3-14 n^2 a^3-n a^3\\
&& +10 a^3-3 n^4 a^2+2 n^3 a^2+13 n^2 a^2-C a^2+D a^2+C n a^2-8 n
   a^2-4 a^2+2 n^4 a\\
   && -4 n^3 a-C n^2 a-2 n^2 a-4 D a+C n a+D n a+4 n a-D n^2+2 D+D n
\end{eqnarray*}
\end{small}
%
The zeros of denominator terms in above expressions consist of $a=0,
a=1$ and $a=2$. Since $(b-a-2)>4$ and $(b+a+2)=-4 $ imply $a<-6$,
the denominator terms in above expressions never vanish.
\end{itemize}
\newpage
\noindent {\bf 4.3.3.4 The case $J^{6.IV}:$\quad $(b-a-2)>4,\
(b+a+2)<-4 $}

\[
\begin{split}
 & p(x)  =\fr{(1+x)^{\fr{b-a-2}{2}}}{(1-x)^{\fr{b+a+2}{2}}}  \\
 & a_6 (x) =(1-x^2)^3  ,\\
 & a_5 (x)  =3\,\left( b + \left( -4 + a \right) \,x \right) \,{\left( -1 + x^2 \right) }^2  ,\\
 & a_4(x)  =-\left( \frac{-1 + x^2}{a\,\left( 2 - 3\,a + a^2 \right) } \right) X_4 \mbox{ with }\\
 & \quad X_4 = D + 3\,a^5\,x^2 - D\,x^2 + a^2\,\left( 42 - 9\,b^2 + 66\,b\,x - 150\,x^2 \right)  + a^4\,\left( 3 + 6\,b\,x - 30\,x^2 \right)
 \\ & \quad +
  3\,a^3\,\left( -7 + b^2 - 12\,b\,x + 35\,x^2 \right)  + a\,\left( -24 + 6\,b^2 + C - 36\,b\,x + 72\,x^2 - C\,x^2
  \right) \\
 & a_3(x)  = \frac{X_3}{a\,\left( 2 - 3\,a + a^2 \right) }  \mbox{ with }\\
 & \quad X_3 = a^6\,x^3 - 2\,D\,\left( b - 2\,x \right) \,\left( -1 + x^2 \right)  + 3\,a^5\,\left( x + b\,x^2 - 4\,x^3 \right) \\
 &  \quad +
  a^4\left( 3b - 27x + 3b^2x - 24bx^2 + 55x^3 \right)
  +
  a^3\left( b^3 + 84x - 15b^2x - 120x^3 + b\left( -19 + 69x^2 \right)  \right) \\
  & \quad +
  2a\left( b^3 - 6b^2x + \left( -24 + 2C - D \right) x\left( -1 + x^2 \right)  +
     b\,\left( -10 + C + 18\,x^2 - C\,x^2 \right)  \right)  \\
     & \quad+
  a^2\,\left( -3\,b^3 + 24\,b^2\,x + b\,\left( 36 - 84\,x^2 \right)  + 2\,x\,\left( -54 + C + 62\,x^2 - C\,x^2 \right)
  \right), \\
 & a_2(x)  = \frac{X_2}{a\,\left( 2 - 3\,a + a^2 \right) }  \mbox{ with }\\
& \quad X_2 = b\,D\,\left( b - 2\,x \right)  + a^3\,C\,x^2
a^2\,\left( C + D + 2\,b\,C\,x - 3\,C\,x^2 \right)  + \\
 & \quad + a\,\left( 2\,D\,\left( -1 + b\,x \right)  + C\,\left( -2 + b^2 - 2\,b\,x + 2\,x^2 \right)
  \right) \\
 & a_1(x)  =\frac{D\,\left( b + a\,x \right) }{a}
\end{split} \]
with eigenvalues
$$
\lambda_n =  -\frac{n (-a+n-1) }{(a-2) (a-1) a} \Lambda_n
$$
where
\begin{small}
\begin{eqnarray*}
\Lambda_n &=&n^2 a^5-3 n a^5+2 a^5-2 n^3 a^4+2 n^2 a^4+8 n a^4-8 a^4+n^4 a^3+4 n^3 a^3-14 n^2 a^3-n a^3\\
&&+10 a^3-3 n^4 a^2+2 n^3 a^2+13 n^2 a^2-C a^2+D a^2+C n a^2-8 n
   a^2-4 a^2+2 n^4 a\\
   && -4 n^3 a-C n^2 a-2 n^2 a-4 D a+C n a+D n a+4 n a-D n^2+2 D+D n
\end{eqnarray*}
\end{small}
%
The zeros of denominator terms in above expressions consist of $a=0,
a=1$ and $a=2$. Since $(b-a-2)>4$ and $(b+a+2)<-4 $ imply $a<-6$,
the denominator terms in above expressions never vanish.

%
%
\subsection{Self-adjoint operators of order $8$}
This section provides examples of the self-adjoint operators
\begin{equation}\label{main-op-n-8}
L=a_8(x) y^{(8)}+a_7(x) y^{(7)}+a_6(x) y^{(6)}+a_5(x) y^{(5)}+a_4(x) y^{(4)}+a_3(x) y''' + a_2(x) y''+a_1(x) y'
\end{equation}
with an admissible weight $p(x)=\dfrac{e^{\frac{1}{4}\int \frac{a_7(x)}{a_8(x)}dx}}{\left\vert a_8(x)\right\vert }$,
satisfying the differential equation
\begin{equation}\label{main-eq-n-8}
L(y)=\lambda y.
\end{equation}
By Proposition 2.1, the determining equations in this case are
\begin{eqnarray}
& 4(a_8p)' =  (a_7 p) \label{det1-n-8} \\
& 7(a_7p)'' - 6(a_6p)' + 2(a_5p)=0 \label{det2-n-8n}\\
& (a_6p)'''-2(a_5p)'' +2(a_4p)' - (a_3p)=0 \label{det3-n-8n}\\
& (a_5p)^{(4)} - 4 (a_4p)'''+ 9(a_3p)'' -14(a_2p)' + 14 (a_1p)=0       \label{det4-n-8n}
\end{eqnarray}
on $I$, subject to the vanishing of
\begin{eqnarray}\label{boundary-n-8}
&& (a_8 p),\ (a_7 p),   \ (a_7p)', \ (a_6p), \  (a_6p)', \ (a_5 p), \
5(a_6p)''-11(a_5p)' +14(a_4p),  \nonumber \\
&& 9(a_6p)''-17(a_5p)' +14(a_4p), \ (a_5p)''-4(a_4p)' +9(a_3p), \nonumber \\
&& (a_5p)''' - 4(a_4p)''+9(a_3p)' -14(a_2p)
\end{eqnarray}
 on the boundary $\partial I$.\\

 In the following sections, we  give examples of eighth order self-adjoint operators with the weights of the form $p(x)=e^{-x^2}$, $p(x)=|x|^{n} e^{m x}$ and $p(x)=\frac{(1+x)^m}{(1-x)^n}$ by solving  the determining equations and boundary conditions using Mathematica.
\subsubsection{Eighth order operators (\ref{main-op-n-8}) with $a_8=1$ and weight of the form $e^{-x^2}$}
\[
\begin{split}
 & p(x)  =e^{-x^2}  \\
 & a_8 (x) =1 , \qquad  a_7 (x)  =-8\,x  , \qquad a_6 (x) =\frac{-144 + D - 2\,F - G + 288\,x^2}{12} ,\\
 & a_5 (x)  =\frac{x\,\left( 96 - D + 2\,F + G - 64\,x^2 \right) }{2}  ,\\
 & a_4(x)  = \frac{48 + G - 192\,x^2 - 4\,G\,x^2 + 64\,x^4 + F\,\left( 2 - 8\,x^2 \right)  + D\,\left( -2 + 4\,x^2 \right) }{4}  ,\\
 & a_3(x)  = D\,x - \frac{2\,\left( D - 2\,F - G \right) \,x^3}{3} ,\qquad a_2(x)  = F - \left( 2\,F + G \right) \,x^2  ,\qquad a_1(x)  = G\,x
\end{split}\]
with eigenvalues
$$
\lambda_n = \frac{1}{3} n (48 n^3-2 D n^2+4 F n^2+2 G n^2-288 n^2+6 D n-18 F n-9 G n+528 n-4 D+14 F+10 G-288)  .
$$
The specific case $D=-256$, $F=-64$, $G=16$ gives the standard eighth order Hermite operator
\begin{eqnarray*}
&& a_8(x)=1,\ a_7(x)= -8x,\   a_6(x)=24\,{x}^{2}-24, \ a_5(x)=120x-32{x}^{3}\\
&& a_4(x)= 112-192{x}^{2}+16{x}^{4},  \ a_3(x) =-256x+96{x}^{3}, \ a_2(x)=-64+112{x}^{2},  \ a_1(x)=16x
\end{eqnarray*}
with the eigenvalues $\lambda_n = 16 n^4$, and coincides with fourth iteration of the classical second order Hermite operator.

%
\subsubsection{Eighth order operators (\ref{main-op-n-8}) with $a_8=x^4$ and weight of the form $|x|^{n} e^{m x}$}
\[
\begin{split}
& p(x)  = |x|^{\fr{a}{4}-4} e^{\fr{b}{4}x} \mb{ with } a>12, b<0 \\
  & a_8 (x) =x^4 ,\\
 & a_7 (x)  =x^3\,\left( a + b\,x \right)  ,\\
 & a_6 (x) = \frac{x^2\,\left( -12\,a + 3\,a^2 + x\,\left( 8\,A + 3\,b^2\,x \right)  \right) }{8} ,\\
 & a_5 (x)  = \frac{x\,\left( a^3 - 6\,a^2\,\left( 2 + b\,x \right)  + 2\,a\,\left( 16 + 6\,A\,x + 12\,b\,x - 3\,b^2\,x^2 \right)  +
      x\,\left( b^3\,x^2 + 12\,A\,\left( -4 + b\,x \right)  \right)  \right) }{16} ,\\
 & a_4(x)  = \frac{X_4}{256} \mbox{ with }\\
 & \quad X_4 =  a^4 - 8\,a^3\,\left( 3 + 4\,b\,x \right)  + 16\,a^2\,\left( 11 + 3\,A\,x + 24\,b\,x \right)  -
  32\,a\,\left( 12 + 18\,A\,x + 32\,b\,x + b^3\,x^3 \right) \\
  & \quad  +
  x\,\left( 256\,C\,x + b^4\,x^3 + 48\,A\,\left( 32 + b^2\,x^2 \right)  \right) ,\\
 & a_3(x)  = \frac{X_3}{256}  \mbox{ with } \\
 & \quad X_3 = -3\,a^4\,b + 128\,C\,x\,\left( -8 + b\,x \right)  + 24\,a^3\,b\,\left( 3 + b\,x \right)  +
  24\,a^2\,b\,\left( -22 - 12\,b\,x + b^2\,x^2 \right) \\
  & \quad +
  a\,\left( 1152\,b + 768\,b^2\,x + 128\,C\,x - 96\,b^3\,x^2 - 3\,b^4\,x^3 \right) \\
  & \quad +
  4\,A\,\left( -384 + a^3 - 288\,b\,x + 36\,b^2\,x^2 + b^3\,x^3 - 3\,a^2\,\left( 8 + 3\,b\,x \right)  +
     a\,\left( 176 + 108\,b\,x - 9\,b^2\,x^2 \right)  \right) ,\\
 & a_2(x)  = \frac{X_2}{2048} \mbox{ with } \\
 & \quad X_2 =  -48\,\left( -384 + 176\,a - 24\,a^2 + a^3 \right) \,A\,b + 33\,a\,\left( -384 + 176\,a - 24\,a^2 + a^3 \right) \,b^2 \\
 & \quad +
  128\,\left( 96 - 20\,a + a^2 \right) \,C + 2048\,F\,x +
  b^2\,\left( -48\,\left( -4 + a \right) \,A\,b - 132\,a\,b^2 + 33\,a^2\,b^2 + 128\,C \right) \,x^2 ,\\
 & a_1(x)  = \left( \frac{-12 + a + b\,x}{4096} \right) X_1 \mbox{ with } \\
 & \quad X_1 = 48\,\left( 32 - 12\,a + a^2 \right) \,A\,b^2 + 396\,a^2\,b^3 - 33\,a^3\,b^3 - 32\,a\,\left( 33\,b^3 + 4\,b\,C \right)  +
  1024\,\left( b\,C + F \right)
\end{split}\]
with eigenvalues
$$
\lambda_n=  \frac{b n }{4096} \Lambda_n
$$
where
\begin{small}
\begin{eqnarray*}
\Lambda_n &=&   -33 a^3 b^3+16 n^3 b^3+330 a^2 b^3-48 a n^2 b^3-96 n^2 b^3-888 a b^3+66 a^2 n b^3-120 a n b^3\\
&& +176 n b^3 -96 b^3+64 A n^2 b^2+48 a^2 A b^2-480 a A b^2+1280 A
   b^2-96 a A n b^2+192 A n b^2\\
   && -128 a C b +768 C b+256 C n b+1024 F
\end{eqnarray*}
\end{small}

The special case $a=16, b = -4, A=-54, C=187, F=-29$ is the standard eighth order Laguerre operator with the weight $p(x) = e^{-x}$ as
\[
\begin{split}
  & a_8 (x) = x^4 ,\\
 & a_7 (x)  = -4\,\left( -4 + x \right) \,x^3 ,\\
 & a_6 (x) =6 x^2 \left(x^2-9 x+12\right)  ,\\
 & a_5 (x)  = -2 x \left(2 x^3-33 x^2+99 x-48\right),\\
 & a_4(x)  = x^4-34 x^3+187 x^2-204 x+24 ,\\
 & a_3(x)  = 6 x^3-68 x^2+136 x-36 ,\\
 & a_2(x)  = 7 x^2-29 x+14 ,\\
 & a_1(x)  = x-1
\end{split}\]
with eigenvalues
$
\lambda_n=n^4,
$
which coincides with the fourth iteration of the classical second order Laguerre operator.
\newpage
%
%
\subsubsection{Eighth order operators (\ref{main-op-n-8}) with $a_8=(1-x^2)^4 $ and weight of the form $\fr{(1+x)^{\fr{b-a-2}{2}}}{(1-x)^{\fr{b+a+2}{2}}}$ with $(b-a-2)>6 \mbox{ or } (b-a-2) \in \left\{0,2,4,6 \right\}$, and $(b+a+2)<-6 \mbox{ or } (b+a+2) \in \left\{-6,-4,-2,0 \right\}$
}

\noindent {\bf 4.4.3.1 The case $J^{8.I}:$\quad $(b-a-2) \in
\left\{0,2,4,6 \right\},\ (b+a+2) \in \left\{-6,-4,-2,0 \right\} $}
\begin{itemize}
\item $J^{8.I.a}:$ \ $b-a-2=0,\ b+a+2=0$
\begin{small}
\[
\begin{split}
 & p(x)  = 1  \\
 & a_8 (x) = {\left( 1 - x^2 \right) }^4 ,\\
 & a_7 (x)  = -32\,x\,{\left( 1 - x^2 \right) }^3  ,\\
 & a_6 (x) = \frac{{\left( -1 + x^2 \right) }^2\,\left( -3456 + D - G + 24192\,x^2 - D\,x^2 + G\,x^2 + 2\,F\,\left( -1 + x^2 \right)  \right)
      }{72} ,\\
 & a_5 (x)  = \frac{x\,\left( -1 + x^2 \right) \,\left( -2304 + D - G + 5376\,x^2 - D\,x^2 + G\,x^2 + 2\,F\,\left( -1 + x^2 \right)  \right) }
  {4} ,\\
 & a_4(x)  = \frac{X_4}{8} \mbox{ with } \\
 & \quad X_4 = 1152 + G - 11520\,x^2 - 10\,G\,x^2 + 13440\,x^4 + 9\,G\,x^4 - 2\,D\,\left( 1 - 6\,x^2 + 5\,x^4 \right) \\
 & \quad  +
  2\,F\,\left( 1 - 10\,x^2 + 9\,x^4 \right) ,\\
 & a_3(x)  = \frac{2\,\left( 2\,F + G \right) \,x^3 + D\,\left( 3\,x - 5\,x^3 \right) }{3}  ,\\
 & a_2(x)  = F - \left( 3\,F + G \right) \,x^2 ,\\
 & a_1(x)  = G\,x
\end{split} \]
\end{small}
with eigenvalues
$$
\lambda_n = \frac{1}{72} n (n+1) \Lambda_n
$$
where
\begin{small}
\begin{eqnarray*}
\Lambda_n &=&   72 n^6+216 n^5-D n^4+2 F n^4+G n^4-1224 n^4-2 D n^3+4 F n^3+2 G n^3\\
&& -2808 n^3+7 D n^2-32 F n^2-16 G n^2+6336 n^2+8 D n-34 F n-17 G n\\
&& +7776 n-12 D+60
   F+66 G-10368
\end{eqnarray*}
\end{small}
The special case $D=-1856,\ F=-216,\ G=16$ gives the standard eighth order Legendre operator with the weight $p(x) = 1$ as
\[
\begin{split}
 a_8 (x) & = {\left( 1 - x^2 \right) }^4,\ \ a_7 (x)  = -32\,x\,{\left( 1 - x^2 \right) }^3, \ \  a_6(x)   = 4 \left(x^2-1\right)^2 \left(89 x^2-17\right), \\ a_5(x)  &  =24 x \left(x^2-1\right) \left(71 x^2-39\right) , \ \ a_4(x)  = 4 \left(883 x^4-926 x^2+139\right), \\
 a_3(x)  &  =64 x \left(44 x^2-29\right) , \ \ a_2(x)    = 8 \left(79 x^2-27\right), \ \ a_1(x)    = 16 x
\end{split}\]
with eigenvalues
$
\lambda_n=n^4 (n+1)^4,
$
which coincides with the fourth iteration of the classical second order Legendre operator.
\item $J^{8.I.b}:$ \ $b-a-2=0,\ b+a+2=-2$
\begin{small}
\[
\begin{split}
 & p(x)  = 1 - x \\
 & a_8 (x) ={\left( 1 - x^2 \right) }^4 ,\\
 & a_7 (x)  = 4\,\left( 1 + 9\,x \right) \,{\left( -1 + x^2 \right) }^3 ,\\
 & a_6 (x) = \frac{{\left( -1 + x^2 \right) }^2\left( -1728 - G + 3456x + 15552x^2 + Gx^2 - 2D\left( -1 + x^2 \right)  +
      F\left( -1 + x^2 \right)  \right) }{36},\\
 & a_5 (x)  = \frac{(1 - x^2)}{12}  X_5 \mbox{ with } \\
 & \quad X_5 = 1152 + G + 8064\,x + 7\,G\,x - 8064\,x^2 - G\,x^2 - 24192\,x^3 - 7\,G\,x^3 \\
 & \quad + 2\,D\,\left( -1 - 7\,x + x^2 + 7\,x^3 \right)  -
  F\,\left( -1 - 7\,x + x^2 + 7\,x^3 \right) ,\\
 & a_4(x)  = \frac{X_4}{4} \mbox{ with } \\
 & \quad X_4 = 576 + G - 2304\,x - 4\,G\,x - 8064\,x^2 - 14\,G\,x^2 + 5376\,x^3 + 4\,G\,x^3 + 12096\,x^4 \\
 & \quad + 13\,G\,x^4 -
  4\,D\,\left( 1 - 2\,x - 8\,x^2 + 2\,x^3 + 7\,x^4 \right)  + F\,\left( 1 - 4\,x - 14\,x^2 + 4\,x^3 + 13\,x^4 \right) ,\\
 & a_3(x)  = \frac{2\,\left( F + G \right) \,x^2\,\left( 3 + 5\,x \right) }{3} + D\,\left( 1 + 5\,x - 5\,x^2 - \frac{35\,x^3}{3} \right) ,\\
 & a_2(x)  = -2\,G\,x\,\left( 1 + 2\,x \right)  + F\,\left( 1 - 2\,x - 5\,x^2 \right) ,\\
 & a_1(x)  = G + 3\,G\,x
\end{split} \]
\end{small}
with eigenvalues
$$
\lambda_n = \frac{1}{36} n (n+2) \Lambda_n
$$
where
\begin{small}
\begin{eqnarray*}
\Lambda_n &=&   36 n^6+216 n^5-2 D n^4+F n^4+G n^4-504 n^4-8 D n^3+4 F n^3+4 G n^3-3456 n^3\\
&& +14 D n^2-16 F n^2-16 G n^2+3060 n^2+44 D n-40 F n-40 G n+13608 n\\
&& -48 D+51
   F+87 G-12960
\end{eqnarray*}
\end{small}
\item $J^{8.I.c}:$ \ $b-a-2=0,\ b+a+2=-4$
\begin{small}
\[
\begin{split}
 & p(x)  = {\left( -1 + x \right) }^2  \\
 & a_8 (x) ={\left( 1 - x^2 \right) }^4 ,\\
 & a_7 (x)  = 8\,\left( 1 + 5\,x \right) \,{\left( -1 + x^2 \right) }^3 ,\\
 & a_6 (x) = \frac{{\left( -1 + x^2 \right) }^2\,\left( -1296 + D - G + 7776\,x + 19440\,x^2 - D\,x^2 + G\,x^2 +
      2\,F\,\left( -1 + x^2 \right)  \right) }{36} ,\\
 & a_5 (x)  = \frac{(1 - x^2)}{6}   X_5 \mbox{ with } \\
 & \quad X_5 =  1152 + G + 3456\,x + 4\,G\,x - 10368\,x^2 - G\,x^2 - 17280\,x^3 - 4\,G\,x^3 \\
 & \quad  + D\,\left( -1 - 4\,x + x^2 + 4\,x^3 \right)  -
  2\,F\,\left( -1 - 4\,x + x^2 + 4\,x^3 \right)   ,\\
 & a_4(x)  =  \frac{X_4}{12}   \mbox{ with }  \\
 & \quad X_4 =   576 + G - 16128\,x - 28\,G\,x - 24192\,x^2 - 54\,G\,x^2 + 48384\,x^3 + 28\,G\,x^3 \\
 & \quad  + 60480\,x^4 + 53\,G\,x^4 -
  4\,D\,\left( 1 - 7\,x - 15\,x^2 + 7\,x^3 + 14\,x^4 \right)  \\
  & \quad  + 2\,F\,\left( 1 - 28\,x - 54\,x^2 + 28\,x^3 + 53\,x^4 \right)  ,\\
 & a_3(x)  =  \frac{\left( 2\,F + G \right) \,x\,\left( 3 + 18\,x + 19\,x^2 \right) }{3} +
  D\,\left( 1 + 2\,x - 7\,x^2 - \frac{28\,x^3}{3} \right) ,\\
 & a_2(x)  =  -\left( G\,x\,\left( 5 + 7\,x \right)  \right)  + F\,\left( 1 - 10\,x - 15\,x^2 \right) ,\\
 & a_1(x)  =  G + 2\,G\,x
\end{split} \]
\end{small}
with eigenvalues
$$
\lambda_n = \frac{1}{36} n (n+3) \Lambda_n
$$
where
\begin{small}
\begin{eqnarray*}
\Lambda_n &=& 36 n^6+324 n^5-D n^4+2 F n^4+G n^4-180 n^4-6 D n^3+12 F n^3+6 G n^3-5940 n^3\\
&& +5 D n^2-28 F n^2-14 G n^2+144 n^2+42 D n-138 F n-69 G n+31536 n-40 D\\
&& +152
   F+94 G-25920
\end{eqnarray*}
\end{small}
\item $J^{8.I.d}:$ \ $b-a-2=0,\ b+a+2=-6$
\begin{small}
\[
\begin{split}
 & p(x)  = {\left( 1 - x \right) }^3 \\
 & a_8 (x) ={\left( 1 - x^2 \right) }^4 ,\\
 & a_7 (x)  = 4\,\left( 3 + 11\,x \right) \,{\left( -1 + x^2 \right) }^3 ,\\
 & a_6 (x) = \frac{-\left( {\left( -1 + x^2 \right) }^2\left( 576 - G - 17280x - 31680x^2 + Gx^2 + D\left( -1 + x^2 \right)  +
        3F\left( -1 + x^2 \right)  \right)  \right) }{48}  ,\\
 & a_5 (x)  =  \frac{-3\,\left( -1 + x^2 \right) }{16}   X_5 \mbox{ with } \\
 & \quad X_5 =  1408 - G + 1152\,x - 3\,G\,x - 17280\,x^2 + G\,x^2 - 21120\,x^3 + 3\,G\,x^3 \\
 & \quad  + D\,\left( -1 - 3\,x + x^2 + 3\,x^3 \right)  +
  3\,F\,\left( -1 - 3\,x + x^2 + 3\,x^3 \right)   ,\\
 & a_4(x)  =  \frac{X_4}{6} \mbox{ with }  \\
 & \quad X_4 =  F\,\left( 3 + 54\,x + 75\,x^2 - 54\,x^3 - 78\,x^4 \right)  + G\,\left( 1 + 18\,x + 25\,x^2 - 18\,x^3 - 26\,x^4 \right) \\
 & \quad   -
  9\,\left( 96 - 2\,\left( -704 + D \right) \,x - 3\,\left( -192 + D \right) \,x^2 + 2\,\left( -2880 + D \right) \,x^3 +
     3\,\left( -1760 + D \right) \,x^4 \right)  ,\\
 & a_3(x)  = \frac{\left( -3\,F - G \right) \,x\,\left( 14 + 57\,x + 49\,x^2 \right)  + D\,\left( 6 - 63\,x^2 - 63\,x^3 \right) }{6} ,\\
 & a_2(x)  =  \frac{2\,G\,x\,\left( 9 + 11\,x \right) }{3} + F\,\left( 1 + 18\,x + 21\,x^2 \right) ,\\
 & a_1(x)  =  G + \frac{5\,G\,x}{3}
\end{split} \]
\end{small}
with eigenvalues
$$
\lambda_n = \frac{1}{48} n (n+4) \Lambda_n
$$
where
\begin{small}
\begin{eqnarray*}
\Lambda_n &=& 48 n^6+576 n^5-D n^4-3 F n^4-G n^4+480 n^4-8 D n^3-24 F n^3-8 G n^3-11520 n^3\\
&& +D n^2+27 F n^2+9 G n^2-9168 n^2+68 D n+300 F n+100 G n+80064 n-60 D\\
&& -300
   F-84 G-60480
\end{eqnarray*}
\end{small}
%
\item $J^{8.I.e}:$ \ $b-a-2=2,\ b+a+2=0$
\begin{small}
\[
\begin{split}
 & p(x)  = 1 + x \\
 & a_8 (x) ={\left( 1 - x^2 \right) }^4 ,\\
 & a_7 (x)  =  4\,\left( -1 + 9\,x \right) \,{\left( -1 + x^2 \right) }^3 ,\\
 & a_6 (x) =  \frac{{\left( -1 + x^2 \right) }^2\left( -1728 + G - 3456x + 15552x^2 - Gx^2 + 2D\left( -1 + x^2 \right)  +
      F\left( -1 + x^2 \right)  \right) }{36} ,\\
 & a_5 (x)  =  \frac{-1 + x^2}{12}  X_5 \mbox{ with }  \\
 & \quad X_5 =   1152 - G - 8064\,x + 7\,G\,x - 8064\,x^2 + G\,x^2 + 24192\,x^3 - 7\,G\,x^3 \\
 & \quad + 2\,D\,\left( 1 - 7\,x - x^2 + 7\,x^3 \right)  +
  F\,\left( 1 - 7\,x - x^2 + 7\,x^3 \right)  ,\\
 & a_4(x)  = \frac{X_4}{4} \mbox{ with } \\
 & \quad X_4 =   576 - G + 2304\,x - 4\,G\,x - 8064\,x^2 + 14\,G\,x^2 - 5376\,x^3 + 4\,G\,x^3 + 12096\,x^4 \\
 & \quad - 13\,G\,x^4 +
  4\,D\,\left( 1 + 2\,x - 8\,x^2 - 2\,x^3 + 7\,x^4 \right)  + F\,\left( 1 + 4\,x - 14\,x^2 - 4\,x^3 + 13\,x^4 \right)  ,\\
 & a_3(x)  = \frac{2\,\left( F - G \right) \,x^2\,\left( -3 + 5\,x \right) }{3} + D\,\left( 1 - 5\,x - 5\,x^2 + \frac{35\,x^3}{3} \right) ,\\
 & a_2(x)  = 2\,G\,x\,\left( -1 + 2\,x \right)  + F\,\left( 1 + 2\,x - 5\,x^2 \right) ,\\
 & a_1(x)  = G - 3\,G\,x
\end{split} \]
\end{small}
with eigenvalues
$$
\lambda_n = \frac{1}{36} n (n+2) \Lambda_n
$$
where
\begin{small}
\begin{eqnarray*}
\Lambda_n &=& 36 n^6+216 n^5+2 D n^4+F n^4-G n^4-504 n^4+8 D n^3+4 F n^3-4 G n^3-3456 n^3\\
&& -14 D n^2-16 F n^2+16 G n^2+3060 n^2-44 D n-40 F n+40 G n+13608 n+48 D\\
&&+51
   F-87 G-12960
\end{eqnarray*}
\end{small}
\item $J^{8.I.f}:$ \ $b-a-2=2,\ b+a+2=-2$
\begin{small}
\[
\begin{split}
 & p(x)  = 1 - x^2 \\
 & a_8 (x) ={\left( 1 - x^2 \right) }^4 ,\\
 & a_7 (x)  = 40\,x\,{\left( -1 + x^2 \right) }^3 ,\\
 & a_6 (x) = \frac{{\left( -1 + x^2 \right) }^2\,\left( -17280 - G + 155520\,x^2 + G\,x^2 - 2\,D\,\left( -1 + x^2 \right)  +
      4\,F\,\left( -1 + x^2 \right)  \right) }{288} ,\\
 & a_5 (x)  = \frac{x\,\left( -1 + x^2 \right) \,\left( -11520 - G + 34560\,x^2 + G\,x^2 - 2\,D\,\left( -1 + x^2 \right)  +
      4\,F\,\left( -1 + x^2 \right)  \right) }{12} ,\\
 & a_4(x)  = \frac{X_4}{24} \mbox{ with }  \\
& \quad X_4 =  5760 + G - 80640\,x^2 - 14\,G\,x^2 + 120960\,x^4 +
13\,G\,x^4 - 4\,D\,\left( 1 - 8\,x^2 + 7\,x^4 \right) \\
& \quad  +
  F\,\left( 4 - 56\,x^2 + 52\,x^4 \right)  ,\\
 & a_3(x)  = \frac{2\,\left( 4\,F + G \right) \,x^3 + D\,\left( 3\,x - 7\,x^3 \right) }{3} ,\\
 & a_2(x)  = F - \left( 5\,F + G \right) \,x^2 ,\\
 & a_1(x)  = G\,x
\end{split} \]
\end{small}
with eigenvalues
$$
\lambda_n = \frac{1}{288} n (n+3) \Lambda_n
$$
where
\begin{small}
\begin{eqnarray*}
\Lambda_n &=& 288 n^6+2592 n^5-2 D n^4+4 F n^4+G n^4-1440 n^4-12 D n^3+24 F n^3+6 G n^3\\
&& -47520 n^3+10 D n^2-68 F n^2-17 G n^2+1152 n^2+84 D n-312 F n-78 G n\\
&& +252288
   n-80 D+352 F +160 G-207360
\end{eqnarray*}
\end{small}
\item $J^{8.I.g}:$ \ $b-a-2=2,\ b+a+2=-4$
\begin{small}
\[
\begin{split}
 & p(x)  = {\left( 1 - x \right) }^2\,\left( 1 + x \right) \\
 & a_8 (x) ={\left( 1 - x^2 \right) }^4 ,\\
 & a_7 (x)  = 4\,\left( 1 + 11\,x \right) \,{\left( -1 + x^2 \right) }^3 ,\\
 & a_6 (x) = \frac{{\left( -1 + x^2 \right) }^2\left( -4320 - G + 8640x + 47520x^2 + Gx^2 - 3D\left( -1 + x^2 \right)  +
      F\left( -1 + x^2 \right)  \right) }{72}  ,\\
 & a_5 (x)  = \frac{(1 - x^2)}{24}  X_5 \mbox{ with } \\
 & \quad X_5 =   2880 + G + 25920\,x + 9\,G\,x - 25920\,x^2 - G\,x^2 - 95040\,x^3 - 9\,G\,x^3 \\
 & \quad + 3\,D\,\left( -1 - 9\,x + x^2 + 9\,x^3 \right)  -
  F\,\left( -1 - 9\,x + x^2 + 9\,x^3 \right)   ,\\
 & a_4(x)  =  \frac{X_4}{6}  \mbox{ with }   \\
 & \quad X_4 = 1440 + G - 5760\,x - 4\,G\,x - 25920\,x^2 - 18\,G\,x^2 + 17280\,x^3 + 4\,G\,x^3 + 47520\,x^4 \\
 & \quad + 17\,G\,x^4 -
  6\,D\,\left( 1 - 2\,x - 10\,x^2 + 2\,x^3 + 9\,x^4 \right)  + F\,\left( 1 - 4\,x - 18\,x^2 + 4\,x^3 + 17\,x^4 \right)  ,\\
 & a_3(x)  = \frac{2\,\left( F + G \right) \,x^2\,\left( 3 + 7\,x \right) }{3} + D\,\left( 1 + 7\,x - 7\,x^2 - 21\,x^3 \right) ,\\
 & a_2(x)  =  2\,G\,x\,\left( 1 + 3\,x \right)  + F\,\left( 1 - 2\,x - 7\,x^2 \right) ,\\
 & a_1(x)  = G + 5\,G\,x
\end{split} \]
\end{small}
with eigenvalues
$$
\lambda_n = \frac{1}{72} n (n+4) \Lambda_n
$$
where
\begin{small}
\begin{eqnarray*}
\Lambda_n &=& 72 n^6+864 n^5-3 D n^4+F n^4+G n^4+720 n^4-24 D n^3+8 F n^3+8 G n^3-17280 n^3\\
&& +3 D n^2-13 F n^2-13 G n^2-13752 n^2+204 D n-116 F n-116 G n+120096
   n\\
   && -180 D+120 F+192 G-90720
\end{eqnarray*}
\end{small}
\item $J^{8.I.h}:$ \ $b-a-2=2,\ b+a+2=-6$
\begin{small}
\[
\begin{split}
 & p(x)  = {\left( 1 - x \right) }^3\,\left( 1 + x \right) \\
 & a_8 (x) ={\left( 1 - x^2 \right) }^4 ,\\
 & a_7 (x)  = 8\,\left( 1 + 6\,x \right) \,{\left( -1 + x^2 \right) }^3 ,\\
 & a_6 (x) = \frac{{\left( -1 + x^2 \right) }^2\left( -4608 - G + 25344x + 76032x^2 + Gx^2 - 2D\left( -1 + x^2 \right)  +
      2F\left( -1 + x^2 \right)  \right) }{96} ,\\
 & a_5 (x)  = \frac{(1 - x^2)}{16} X_5 \mbox{ with } \\
 & \quad X_5 =   3840 + G + 15360\,x + 5\,G\,x - 42240\,x^2 - G\,x^2 - 84480\,x^3 - 5\,G\,x^3 \\
 & \quad  + 2\,D\,\left( -1 - 5\,x + x^2 + 5\,x^3 \right)  -
  2\,F\,\left( -1 - 5\,x + x^2 + 5\,x^3 \right)   ,\\
 & a_4(x)  = \frac{X_4}{16} \mbox{ with }   \\
 & \quad X_4 = 1920 + G - 34560\,x - 18\,G\,x - 69120\,x^2 - 44\,G\,x^2 + 126720\,x^3 + 18\,G\,x^3 \\
 & \quad + 190080\,x^4 + 43\,G\,x^4 -
  6\,D\,\left( 1 - 6\,x - 16\,x^2 + 6\,x^3 + 15\,x^4 \right)  \\
  & \quad  + F\,\left( 2 - 36\,x - 88\,x^2 + 36\,x^3 + 86\,x^4 \right)  ,\\
 & a_3(x)  = \frac{\left( 2\,F + G \right) \,x\,\left( 1 + 8\,x + 11\,x^2 \right) }{2} + D\,\left( 1 + 3\,x - 9\,x^2 - 15\,x^3 \right) ,\\
 & a_2(x)  = \frac{-\left( G\,x\,\left( 7 + 13\,x \right)  \right) }{2} + F\,\left( 1 - 7\,x - 14\,x^2 \right) ,\\
 & a_1(x)  = G + 3\,G\,x
\end{split} \]
\end{small}
with eigenvalues
$$
\lambda_n = \frac{1}{96} n (n+5) \Lambda_n
$$
where
\begin{small}
\begin{eqnarray*}
\Lambda_n &=& 96 n^6+1440 n^5-2 D n^4+2 F n^4+G n^4+2976 n^4-20 D n^3+20 F n^3+10 G n^3\\
&& -30240 n^3-10 D n^2-14 F n^2-7 G n^2-51456 n^2+200 D n-320 F n-160 G
   n\\
   && +270720 n-168 D+312 F+204 G-193536
\end{eqnarray*}
\end{small}
%
\item $J^{8.I.i}:$ \ $b-a-2=4,\ b+a+2=0$
\begin{small}
\[
\begin{split}
 & p(x)  =  {\left( 1 + x \right) }^2 \\
 & a_8 (x) ={\left( 1 - x^2 \right) }^4 ,\\
 & a_7 (x)  = 8\,\left( -1 + 5\,x \right) \,{\left( -1 + x^2 \right) }^3 ,\\
 & a_6 (x) =  \frac{{\left( -1 + x^2 \right) }^2\left( -1296 + G - 7776x + 19440x^2 - Gx^2 + D\left( -1 + x^2 \right)  +
      2F\left( -1 + x^2 \right)  \right) }{36} ,\\
 & a_5 (x)  = \frac{(-1 + x^2)}{6} X_5 \mbox{ with } \\
 & \quad X_5 =   1152 - G - 3456\,x + 4\,G\,x - 10368\,x^2 + G\,x^2 + 17280\,x^3 - 4\,G\,x^3 \\
 & \quad + D\,\left( 1 - 4\,x - x^2 + 4\,x^3 \right)  +
  F\,\left( 2 - 8\,x - 2\,x^2 + 8\,x^3 \right)   ,\\
 & a_4(x)  = \frac{X_4}{12}  \mbox{ with }   \\
 & \quad X_4 =  576 - G + 16128\,x - 28\,G\,x - 24192\,x^2 + 54\,G\,x^2 - 48384\,x^3 + 28\,G\,x^3 \\
 & \quad  + 60480\,x^4 - 53\,G\,x^4 +
  4\,D\,\left( 1 + 7\,x - 15\,x^2 - 7\,x^3 + 14\,x^4 \right)  \\
  & \quad  + 2\,F\,\left( 1 + 28\,x - 54\,x^2 - 28\,x^3 + 53\,x^4 \right) ,\\
 & a_3(x)  =  \frac{\left( 2\,F - G \right) \,x\,\left( 3 - 18\,x + 19\,x^2 \right) }{3} +
  D\,\left( 1 - 2\,x - 7\,x^2 + \frac{28\,x^3}{3} \right) ,\\
 & a_2(x)  = G\,x\,\left( -5 + 7\,x \right)  + F\,\left( 1 + 10\,x - 15\,x^2 \right) ,\\
 & a_1(x)  = G - 2\,G\,x
\end{split} \]
\end{small}
with eigenvalues
$$
\lambda_n = \frac{1}{36} n (n+3) \Lambda_n
$$
where
\begin{small}
\begin{eqnarray*}
\Lambda_n &=& 36 n^6+324 n^5+D n^4+2 F n^4-G n^4-180 n^4+6 D n^3+12 F n^3-6 G n^3-5940 n^3\\
&& -5 D n^2-28 F n^2+14 G n^2+144 n^2-42 D n-138 F n+69 G n+31536 n+40 D\\
&& +152
   F-94 G-25920
\end{eqnarray*}
\end{small}
\item $J^{8.I.j}:$ \ $b-a-2=4,\ b+a+2=-2$
\begin{small}
\[
\begin{split}
 & p(x)  = \left( 1 - x \right) \,{\left( 1 + x \right) }^2 \\
 & a_8 (x) ={\left( 1 - x^2 \right) }^4 ,\\
 & a_7 (x)  = 4\,\left( -1 + 11\,x \right) \,{\left( -1 + x^2 \right) }^3 ,\\
 & a_6 (x) = \frac{{\left( -1 + x^2 \right) }^2\left( -4320 + G - 8640x + 47520x^2 - Gx^2 + 3D\left( -1 + x^2 \right)  +
      F\left( -1 + x^2 \right)  \right) }{72} ,\\
 & a_5 (x)  = \frac{(-1 + x^2)}{24} X_5 \mbox{ with }  \\
& \quad X_5 =  2880 - G - 25920\,x + 9\,G\,x - 25920\,x^2 + G\,x^2 +
95040\,x^3 - 9\,G\,x^3 \\
& \quad + 3\,D\,\left( 1 - 9\,x - x^2 + 9\,x^3 \right)  +
  F\,\left( 1 - 9\,x - x^2 + 9\,x^3 \right)    ,\\
 & a_4(x)  = \frac{X_4}{6}  \mbox{ with } \\
 & \quad X_4 =  1440 - G + 5760\,x - 4\,G\,x - 25920\,x^2 + 18\,G\,x^2 - 17280\,x^3 + 4\,G\,x^3 \\
 & \quad + 47520\,x^4 - 17\,G\,x^4 +
  6\,D\,\left( 1 + 2\,x - 10\,x^2 - 2\,x^3 + 9\,x^4 \right)  \\
  & \quad + F\,\left( 1 + 4\,x - 18\,x^2 - 4\,x^3 + 17\,x^4 \right) ,\\
 & a_3(x)  =  \frac{2\,\left( F - G \right) \,x^2\,\left( -3 + 7\,x \right) }{3} + D\,\left( 1 - 7\,x - 7\,x^2 + 21\,x^3 \right) ,\\
 & a_2(x)  =  2\,G\,x\,\left( -1 + 3\,x \right)  + F\,\left( 1 + 2\,x - 7\,x^2 \right) ,\\
 & a_1(x)  =  G - 5\,G\,x
\end{split} \]
\end{small}
with eigenvalues
$$
\lambda_n = \frac{1}{72} n (n+4) \Lambda_n
$$
where
\begin{small}
\begin{eqnarray*}
\Lambda_n &=& 72 n^6+864 n^5+3 D n^4+F n^4-G n^4+720 n^4+24 D n^3+8 F n^3-8 G n^3\\
&& -17280 n^3-3 D n^2-13 F n^2+13 G n^2-13752 n^2-204 D n-116 F n\\
&& +116 G n+120096
   n+180 D+120 F-192 G-90720
\end{eqnarray*}
\end{small}
\item $J^{8.I.k}:$ \ $b-a-2=4,\ b+a+2=-4$
\begin{small}
\[
\begin{split}
 & p(x)  = {\left( 1 - x \right) }^2\,{\left( 1 + x \right) }^2 \\
 & a_8 (x) ={\left( 1 - x^2 \right) }^4 ,\\
 & a_7 (x)  = 48\,x\,{\left( -1 + x^2 \right) }^3 ,\\
 & a_6 (x) = \frac{{\left( -1 + x^2 \right) }^2\,\left( -51840 - G + 570240\,x^2 + G\,x^2 - 3\,D\,\left( -1 + x^2 \right)  +
      6\,F\,\left( -1 + x^2 \right)  \right) }{720} ,\\
 & a_5 (x)  =  \frac{x\,\left( -1 + x^2 \right) \,\left( -34560 - G + 126720\,x^2 + G\,x^2 - 3\,D\,\left( -1 + x^2 \right)  +
      6\,F\,\left( -1 + x^2 \right)  \right) }{24}   ,\\
 & a_4(x)  = \frac{X_4}{48}  \mbox{ with }   \\
 & \quad X_4 =  17280 + G - 311040\,x^2 - 18\,G\,x^2 + 570240\,x^4 + 17\,G\,x^4 \\
 & \quad - 6\,D\,\left( 1 - 10\,x^2 + 9\,x^4 \right)  +
  6\,F\,\left( 1 - 18\,x^2 + 17\,x^4 \right) ,\\
 & a_3(x)  = D\,x + \left( -3\,D + 4\,F + \frac{2\,G}{3} \right) \,x^3 ,\\
 & a_2(x)  =  F - \left( 7\,F + G \right) \,x^2,\\
 & a_1(x)  = G\,x
\end{split} \]
\end{small}
with eigenvalues
$$
\lambda_n = \frac{1}{720} n (n+5) \Lambda_n
$$
where
\begin{small}
\begin{eqnarray*}
\Lambda_n &=& 720 n^6+10800 n^5-3 D n^4+6 F n^4+G n^4+22320 n^4-30 D n^3+60 F n^3+10 G n^3\\
&& -226800 n^3-15 D n^2-60 F n^2-10 G n^2-385920 n^2+300 D n-1050 F n-175 G
   n\\
   && +2030400 n-252 D+1044 F+294 G-1451520
\end{eqnarray*}
\end{small}
\item $J^{8.I.l}:$ \ $b-a-2=4,\ b+a+2=-6$
\begin{small}
\[
\begin{split}
 & p(x)  = {\left( 1 - x \right) }^3\,{\left( 1 + x \right) }^2 \\
 & a_8 (x) ={\left( 1 - x^2 \right) }^4 ,\\
 & a_7 (x)  = 4\,\left( 1 + 13\,x \right) \,{\left( -1 + x^2 \right) }^3 ,\\
 & a_6 (x) = \frac{{\left( -1 + x^2 \right) }^2\left( -8640 - G + 17280x + 112320x^2 + Gx^2 - 4D\left( -1 + x^2 \right)  +
      F\left( -1 + x^2 \right)  \right) }{120}  ,\\
 & a_5 (x)  = \frac{(1 - x^2)}{40} X_5 \mbox{ with } \\
 & \quad X_5 =   5760 + G + 63360\,x + 11\,G\,x - 63360\,x^2 - G\,x^2 - 274560\,x^3 - 11\,G\,x^3 \\
 & \quad +
  4\,D\,\left( -1 - 11\,x + x^2 + 11\,x^3 \right)  - F\,\left( -1 - 11\,x + x^2 + 11\,x^3 \right)   ,\\
 & a_4(x)  = \frac{X_4}{8}  \mbox{ with }   \\
 & \quad X_4 =  2880 + G - 11520\,x - 4\,G\,x - 63360\,x^2 - 22\,G\,x^2 + 42240\,x^3 + 4\,G\,x^3 \\
 & \quad + 137280\,x^4 + 21\,G\,x^4 -
  8\,D\,\left( 1 - 2\,x - 12\,x^2 + 2\,x^3 + 11\,x^4 \right)  \\
  & \quad + F\,\left( 1 - 4\,x - 22\,x^2 + 4\,x^3 + 21\,x^4 \right) ,\\
 & a_3(x)  = 2\,\left( F + G \right) \,x^2\,\left( 1 + 3\,x \right)  + D\,\left( 1 + 9\,x - 9\,x^2 - 33\,x^3 \right) ,\\
 & a_2(x)  = -2\,G\,x\,\left( 1 + 4\,x \right)  + F\,\left( 1 - 2\,x - 9\,x^2 \right) ,\\
 & a_1(x)  = G + 7\,G\,x
\end{split} \]
\end{small}
with eigenvalues
$$
\lambda_n = \frac{1}{120} n (n+6) \Lambda_n
$$
where
\begin{small}
\begin{eqnarray*}
\Lambda_n &=& 120 n^6+2160 n^5-4 D n^4+F n^4+G n^4+6960 n^4-48 D n^3+12 F n^3+12 G n^3\\
&& -46080 n^3-52 D n^2-2 F n^2-2 G n^2-128040 n^2+552 D n-228 F n-228 G
   n\\
   && +527760 n-448 D+217 F+337 G-362880
\end{eqnarray*}
\end{small}
%
\item $J^{8.I.m}:$ \ $b-a-2=6,\ b+a+2=0$
\begin{small}
\[
\begin{split}
 & p(x)  =  {\left( 1 + x \right) }^3 \\
 & a_8 (x) ={\left( 1 - x^2 \right) }^4 ,\\
 & a_7 (x)  = 4\,\left( -3 + 11\,x \right) \,{\left( -1 + x^2 \right) }^3 ,\\
 & a_6 (x) = \frac{{\left( -1 + x^2 \right) }^2\left( -576 - G - 17280x + 31680x^2 + Gx^2 + D\,\left( -1 + x^2 \right)  -
      3F\left( -1 + x^2 \right)  \right) }{48} ,\\
 & a_5 (x)  = \frac{3\,\left( -1 + x^2 \right) }{16}  X_5 \mbox{ with } \\
 & \quad X_5 =   1408 + G - 1152\,x - 3\,G\,x - 17280\,x^2 - G\,x^2 + 21120\,x^3 + 3\,G\,x^3 \\
 & \quad + F\,\left( -3 + 9\,x + 3\,x^2 - 9\,x^3 \right)  +
  D\,\left( 1 - 3\,x - x^2 + 3\,x^3 \right)   ,\\
 & a_4(x)  = \frac{X_4}{6} \mbox{ with }   \\
 & \quad X_4 =  F\,\left( 3 - 54\,x + 75\,x^2 + 54\,x^3 - 78\,x^4 \right)  + G\,\left( -1 + 18\,x - 25\,x^2 - 18\,x^3 + 26\,x^4 \right) \\
 & \quad +
  9\,\left( -96 + 2\,\left( 704 + D \right) \,x - 3\,\left( 192 + D \right) \,x^2 - 2\,\left( 2880 + D \right) \,x^3 +
     3\,\left( 1760 + D \right) \,x^4 \right) ,\\
 & a_3(x)  = \frac{\left( -3\,F + G \right) \,x\,\left( 14 - 57\,x + 49\,x^2 \right)  + D\,\left( 6 - 63\,x^2 + 63\,x^3 \right) }{6} ,\\
 & a_2(x)  =  \frac{2\,G\,\left( 9 - 11\,x \right) \,x}{3} + F\,\left( 1 - 18\,x + 21\,x^2 \right) ,\\
 & a_1(x)  = G - \frac{5\,G\,x}{3}
\end{split} \]
\end{small}
with eigenvalues
$$
\lambda_n = \frac{1}{48} n (n+4) \Lambda_n
$$
where
\begin{small}
\begin{eqnarray*}
\Lambda_n &=& 48 n^6+576 n^5+D n^4-3 F n^4+G n^4+480 n^4+8 D n^3-24 F n^3+8 G n^3-11520 n^3\\
&& -D n^2+27 F n^2-9 G n^2-9168 n^2-68 D n+300 F n-100 G n+80064 n+60 D\\
&& -300
   F+84 G-60480
\end{eqnarray*}
\end{small}
\item $J^{8.I.n}:$ \ $b-a-2=6,\ b+a+2=-2$
\begin{small}
\[
\begin{split}
 & p(x)  = \left( 1 - x \right) \,{\left( 1 + x \right) }^3 \\
 & a_8 (x) ={\left( 1 - x^2 \right) }^4 ,\\
 & a_7 (x)  = 8\,\left( -1 + 6\,x \right) \,{\left( -1 + x^2 \right) }^3 ,\\
 & a_6 (x) = \frac{{\left( -1 + x^2 \right) }^2\left( -4608 + G - 25344x + 76032x^2 - Gx^2 + 2D\left( -1 + x^2 \right)  +
      2F\left( -1 + x^2 \right)  \right) }{96} ,\\
 & a_5 (x)  = \frac{(-1 + x^2)}{16} X_5 \mbox{ with } \\
 & \quad X_5 =   3840 - G - 15360\,x + 5\,G\,x - 42240\,x^2 + G\,x^2 + 84480\,x^3 - 5\,G\,x^3 \\
 & \quad + 2\,D\,\left( 1 - 5\,x - x^2 + 5\,x^3 \right)  +
  2\,F\,\left( 1 - 5\,x - x^2 + 5\,x^3 \right)   ,\\
 & a_4(x)  = \frac{X_4}{16} \mbox{ with }   \\
 & \quad X_4 =  1920 - G + 34560\,x - 18\,G\,x - 69120\,x^2 + 44\,G\,x^2 - 126720\,x^3 + 18\,G\,x^3 \\
 & \quad + 190080\,x^4 - 43\,G\,x^4 +
  6\,D\,\left( 1 + 6\,x - 16\,x^2 - 6\,x^3 + 15\,x^4 \right)  \\
  & \quad + F\,\left( 2 + 36\,x - 88\,x^2 - 36\,x^3 + 86\,x^4 \right) ,\\
 & a_3(x)  = \frac{\left( 2\,F - G \right) \,x\,\left( 1 - 8\,x + 11\,x^2 \right) }{2} + D\,\left( 1 - 3\,x - 9\,x^2 + 15\,x^3 \right) ,\\
 & a_2(x)  = \frac{G\,x\,\left( -7 + 13\,x \right) }{2} + F\,\left( 1 + 7\,x - 14\,x^2 \right) ,\\
 & a_1(x)  = G - 3\,G\,x
\end{split} \]
\end{small}
with eigenvalues
$$
\lambda_n = \frac{1}{96} n (n+5) \Lambda_n
$$
where
\begin{small}
\begin{eqnarray*}
\Lambda_n &=& 96 n^6+1440 n^5+2 D n^4+2 F n^4-G n^4+2976 n^4+20 D n^3+20 F n^3-10 G n^3\\
&& -30240 n^3+10 D n^2-14 F n^2+7 G n^2-51456 n^2-200 D n-320 F n+160 G
   n\\
   && +270720 n+168 D+312 F-204 G-193536
\end{eqnarray*}
\end{small}
\item $J^{8.I.o}:$ \ $b-a-2=6,\ b+a+2=-4$
\begin{small}
\[
\begin{split}
 & p(x)  = {\left( 1 - x \right) }^2\,{\left( 1 + x \right) }^3 \\
 & a_8 (x) ={\left( 1 - x^2 \right) }^4 ,\\
 & a_7 (x)  = 4\,\left( -1 + 13\,x \right) \,{\left( -1 + x^2 \right) }^3 ,\\
 & a_6 (x) = \frac{{\left( -1 + x^2 \right) }^2\left( -8640 + G - 17280x + 112320x^2 - Gx^2 + 4D\left( -1 + x^2 \right)  +
      F\left( -1 + x^2 \right)  \right) }{120} ,\\
 & a_5 (x)  = \frac{(-1 + x^2)}{40}  X_5 \mbox{ with } \\
 & \quad X_5 =   5760 - G - 63360\,x + 11\,G\,x - 63360\,x^2 + G\,x^2 + 274560\,x^3 - 11\,G\,x^3 \\
 & \quad + F\,\left( 1 - 11\,x - x^2 + 11\,x^3 \right)  +
  D\,\left( 4 - 44\,x - 4\,x^2 + 44\,x^3 \right)   ,\\
 & a_4(x)  = \frac{X_4}{8} \mbox{ with }   \\
 & \quad X_4 =  2880 - G + 11520\,x - 4\,G\,x - 63360\,x^2 + 22\,G\,x^2 - 42240\,x^3 + 4\,G\,x^3 \\
 & \quad + 137280\,x^4 - 21\,G\,x^4 +
  8\,D\,\left( 1 + 2\,x - 12\,x^2 - 2\,x^3 + 11\,x^4 \right)  \\
  & \quad + F\,\left( 1 + 4\,x - 22\,x^2 - 4\,x^3 + 21\,x^4 \right) ,\\
 & a_3(x)  = 2\,\left( F - G \right) \,x^2\,\left( -1 + 3\,x \right)  + D\,\left( 1 - 9\,x - 9\,x^2 + 33\,x^3 \right) ,\\
 & a_2(x)  = 2\,G\,x\,\left( -1 + 4\,x \right)  + F\,\left( 1 + 2\,x - 9\,x^2 \right)  ,\\
 & a_1(x)  = G - 7\,G\,x
\end{split} \]
\end{small}
with eigenvalues
$$
\lambda_n = \frac{1}{120} n (n+6) \Lambda_n
$$
where
\begin{small}
\begin{eqnarray*}
\Lambda_n &=& 120 n^6+2160 n^5+4 D n^4+F n^4-G n^4+6960 n^4+48 D n^3+12 F n^3-12 G n^3\\
&& -46080 n^3+52 D n^2-2 F n^2+2 G n^2-128040 n^2-552 D n-228 F n+228 G
   n\\
   && +527760 n+448 D+217 F-337 G-362880
\end{eqnarray*}
\end{small}
\item $J^{8.I.p}:$ \ $b-a-2=6,\ b+a+2=-6$
\begin{small}
\[
\begin{split}
 & p(x)  = {\left( 1 - x \right) }^3\,{\left( 1 + x \right) }^3 \\
 & a_8 (x) ={\left( 1 - x^2 \right) }^4 ,\\
 & a_7 (x)  = 56\,x\,{\left( -1 + x^2 \right) }^3 ,\\
 & a_6 (x) = \frac{{\left( -1 + x^2 \right) }^2\,\left( -120960 - G + 1572480\,x^2 + G\,x^2 - 4\,D\,\left( -1 + x^2 \right)  +
      8\,F\,\left( -1 + x^2 \right)  \right) }{1440} ,\\
 & a_5 (x)  =  \frac{x\,\left( -1 + x^2 \right) \,\left( -80640 - G + 349440\,x^2 + G\,x^2 - 4\,D\,\left( -1 + x^2 \right)  +
      8\,F\,\left( -1 + x^2 \right)  \right) }{40} ,\\
 & a_4(x)  = \frac{X_4}{80} \mbox{ with }   \\
 & \quad X_4 =  40320 + G - 887040\,x^2 - 22\,G\,x^2 + 1921920\,x^4 + 21\,G\,x^4 \\
 & \quad + D\,\left( -8 + 96\,x^2 - 88\,x^4 \right)  +
  8\,F\,\left( 1 - 22\,x^2 + 21\,x^4 \right) ,\\
 & a_3(x)  = D\,x - \frac{\left( 11\,D - 2\,\left( 8\,F + G \right)  \right) \,x^3}{3} ,\\
 & a_2(x)  = F - \left( 9\,F + G \right) \,x^2 ,\\
 & a_1(x)  = G\,x
\end{split} \]
\end{small}
with eigenvalues
$$
\lambda_n = \frac{1}{1440} n (n+7) \Lambda_n
$$
where
\begin{small}
\begin{eqnarray*}
\Lambda_n &=& 1440 n^6+30240 n^5-4 D n^4+8 F n^4+G n^4+131040 n^4-56 D n^3+112 F n^3+14 G n^3\\
&& -635040 n^3-92 D n^2+40 F n^2+5 G n^2-2620800 n^2+728 D n-2464 F n-308 G n\\
&& +9313920 n-576 D+2304
   F+468 G-6220800
\end{eqnarray*}
\end{small}
\end{itemize}
\newpage
\noindent {\bf 4.4.3.2 The case $J^{8.II}:$\quad $(b-a-2) \in
\left\{0,2,4,6 \right\},\ (b+a+2)<-6 $}
\begin{itemize}
\item $J^{8.II.a}:$ \ $b-a-2=0,\ (b+a+2)<-6$
\begin{footnotesize}
\[
\begin{split}
 & p(x)  = {\left( 1 - x \right) }^{-2 - a} \\
 & a_8 (x) ={\left( 1 - x^2 \right) }^4 ,\quad  a_7 (x)  = -4\,\left( 2 + a - 6\,x + a\,x \right) \,{\left( -1 + x^2 \right) }^3 ,\\
 & a_6 (x) = \frac{{\left( -1 + x^2 \right) }^2}{\left( 2 + a \right) \,\left( 24 - 50\,a + 35\,a^2 - 10\,a^3 + a^4 \right) }  X_6 \mbox{ with } \\
& \quad X_6 =   -576 + 2352\,a - 2232\,a^2 - 36\,a^3 + 666\,a^4 -
162\,a^5 - 18\,a^6 + 6\,a^7 - 2\,D + a\,D \\
& \quad + a^2\,D - 2\,F + a\,F - 5760\,x +
  7392\,a\,x + 912\,a^2\,x - 3432\,a^3\,x + 660\,a^4\,x + 348\,a^5\,x \\
  & \quad - 132\,a^6\,x + 12\,a^7\,x + 8640\,x^2 - 16848\,a\,x^2 +
  8904\,a^2\,x^2 + 924\,a^3\,x^2 - 2310\,a^4\,x^2 \\
  & \quad + 798\,a^5\,x^2 - 114\,a^6\,x^2 + 6\,a^7\,x^2 + 2\,D\,x^2 - a\,D\,x^2 -
  a^2\,D\,x^2 + 2\,F\,x^2 - a\,F\,x^2   ,\\
 & a_5 (x)  = \frac{-1 + x^2}{\left( 2 + a \right) \,\left( 24 - 50\,a + 35\,a^2 - 10\,a^3 + a^4 \right) } X_5 \mbox{ with } \\
 & \quad X_5 =   4608 - 7680\,a + 800\,a^2 + 3824\,a^3 - 1304\,a^4 - 460\,a^5 + 220\,a^6 - 4\,a^7 - 4\,a^8 \\
 & \quad + 12\,D - 9\,a^2\,D - 3\,a^3\,D +
  12\,F - 3\,a^2\,F - 4608\,x + 19968\,a\,x - 22560\,a^2\,x \\
  & \quad + 4176\,a^3\,x + 5400\,a^4\,x - 2628\,a^5\,x + 180\,a^6\,x +
  84\,a^7\,x - 12\,a^8\,x - 24\,D\,x + 18\,a\,D\,x \\
  & \quad + 9\,a^2\,D\,x - 3\,a^3\,D\,x - 24\,F\,x + 18\,a\,F\,x - 3\,a^2\,F\,x -
  23040\,x^2 + 35328\,a\,x^2 \\
  & \quad - 3744\,a^2\,x^2 - 14640\,a^3\,x^2 + 6072\,a^4\,x^2 + 732\,a^5\,x^2 - 876\,a^6\,x^2 +
  180\,a^7\,x^2 - 12\,a^8\,x^2 \\
  & \quad - 12\,D\,x^2 + 9\,a^2\,D\,x^2 + 3\,a^3\,D\,x^2 - 12\,F\,x^2 + 3\,a^2\,F\,x^2 + 23040\,x^3 -
  50688\,a\,x^3 \\
  & \quad + 34976\,a^2\,x^3 - 3472\,a^3\,x^3 - 6776\,a^4\,x^3 + 3668\,a^5\,x^3 - 836\,a^6\,x^3 + 92\,a^7\,x^3 -
  4\,a^8\,x^3 \\
  & \quad + 24\,D\,x^3 - 18\,a\,D\,x^3 - 9\,a^2\,D\,x^3 + 3\,a^3\,D\,x^3 + 24\,F\,x^3 - 18\,a\,F\,x^3 + 3\,a^2\,F\,x^3  ,\\
 & a_4(x)  = \frac{X_4}{2\,\left( 2 + a \right) \,\left( 24 - 50\,a + 35\,a^2 - 10\,a^3 + a^4 \right) } \mbox{ with }   \\
 & \quad X_4 =  -2304 - 4608\,a + 14320\,a^2 - 6312\,a^3 - 3308\,a^4 + 2426\,a^5 - 80\,a^6 - 148\,a^7 \\
 & \quad + 12\,a^8 + 2\,a^9 - 60\,a\,D +
  18\,a^2\,D + 36\,a^3\,D + 6\,a^4\,D - 24\,F - 34\,a\,F + 9\,a^2\,F + 7\,a^3\,F \\
  & \quad + 27648\,x - 55296\,a\,x + 20160\,a^2\,x +
  21344\,a^3\,x - 15472\,a^4\,x - 152\,a^5\,x + 2240\,a^6\,x \\
  & \quad - 464\,a^7\,x - 16\,a^8\,x + 8\,a^9\,x + 144\,D\,x - 48\,a\,D\,x -
  108\,a^2\,D\,x + 12\,a^4\,D\,x + 144\,F\,x \\
  & \quad - 48\,a\,F\,x - 36\,a^2\,F\,x + 12\,a^3\,F\,x - 13824\,x^2 + 64512\,a\,x^2 -
  87648\,a^2\,x^2 + 35088\,a^3\,x^2 \\
  & \quad + 12024\,a^4\,x^2 - 13284\,a^5\,x^2 + 3168\,a^6\,x^2 + 72\,a^7\,x^2 - 120\,a^8\,x^2 +
  12\,a^9\,x^2 - 144\,D\,x^2 \\
  & \quad + 216\,a\,D\,x^2 - 72\,a^3\,D\,x^2 - 96\,F\,x^2 + 164\,a\,F\,x^2 - 54\,a^2\,F\,x^2 -
  2\,a^3\,F\,x^2 - 46080\,x^3 \\
  & \quad + 86016\,a\,x^3 - 31040\,a^2\,x^3 - 26784\,a^3\,x^3 + 21904\,a^4\,x^3 - 2584\,a^5\,x^3 -
  2240\,a^6\,x^3 \\
  & \quad + 944\,a^7\,x^3 - 144\,a^8\,x^3 + 8\,a^9\,x^3 - 144\,D\,x^3 + 48\,a\,D\,x^3 + 108\,a^2\,D\,x^3 -
  12\,a^4\,D\,x^3 \\
  & \quad - 144\,F\,x^3 + 48\,a\,F\,x^3 + 36\,a^2\,F\,x^3 - 12\,a^3\,F\,x^3 + 34560\,x^4 - 87552\,a\,x^4 +
  77808\,a^2\,x^4 \\
  & \quad - 22696\,a^3\,x^4 - 8428\,a^4\,x^4 + 8890\,a^5\,x^4 - 3088\,a^6\,x^4 + 556\,a^7\,x^4 - 52\,a^8\,x^4 +
  2\,a^9\,x^4 \\
  & \quad + 144\,D\,x^4 - 156\,a\,D\,x^4 - 18\,a^2\,D\,x^4 + 36\,a^3\,D\,x^4 - 6\,a^4\,D\,x^4 + 120\,F\,x^4 -
  130\,a\,F\,x^4 \\
  & \quad + 45\,a^2\,F\,x^4 - 5\,a^3\,F\,x^4 ,\\
  \end{split} \]
 \[
\begin{split}
 & a_3(x)  = \frac{X_3}{48 - 76\,a + 20\,a^2 + 15\,a^3 - 8\,a^4 + a^5} \mbox{ with }   \\
 & \quad X_3 =  24\,D - 28\,a\,D - 22\,a^2\,D + 15\,a^3\,D + 10\,a^4\,D + a^5\,D - 24\,F - 2\,a^2\,F + 2\,a^4\,F \\
 & \quad + 60\,a\,D\,x - 48\,a^2\,D\,x -
  27\,a^3\,D\,x + 12\,a^4\,D\,x + 3\,a^5\,D\,x + 48\,F\,x - 16\,a\,F\,x \\
  & \quad - 4\,a^2\,F\,x - 8\,a^3\,F\,x + 4\,a^4\,F\,x -
  72\,D\,x^2 + 60\,a\,D\,x^2 + 42\,a^2\,D\,x^2 - 27\,a^3\,D\,x^2 \\
  & \quad - 6\,a^4\,D\,x^2 + 3\,a^5\,D\,x^2 - 24\,F\,x^2 +
  32\,a\,F\,x^2 - 2\,a^2\,F\,x^2 - 8\,a^3\,F\,x^2 + 2\,a^4\,F\,x^2 \\
  & \quad + 48\,D\,x^3 - 76\,a\,D\,x^3 + 20\,a^2\,D\,x^3 +
  15\,a^3\,D\,x^3 - 8\,a^4\,D\,x^3 + a^5\,D\,x^3 ,\\
 & a_2(x)  = \frac{X_2}{2\,a\,\left( -2 + a + a^2 \right) }  \mbox{ with } \\
 & \quad X_2 =  2\,a\,F + 5\,a^2\,F + a^3\,F - 4\,G + 2\,a\,G + 2\,a^2\,G - 4\,a\,F\,x + 2\,a^2\,F\,x + 2\,a^3\,F\,x \\
 & \quad + 2\,a\,F\,x^2 -
  3\,a^2\,F\,x^2 + a^3\,F\,x^2 + 4\,G\,x^2 - 2\,a\,G\,x^2 - 2\,a^2\,G\,x^2 ,\\
 & a_1(x)  = \frac{G\,\left( 2 + a + a\,x \right) }{a}
\end{split} \]
\end{footnotesize}
with eigenvalues
$$
\lambda_n = \frac{n (-a+n-1)}{2 (a-4) (a-3) (a-2) (a-1) a (a+2)} \Lambda_n
$$
where
\begin{small}
\begin{eqnarray*}
\Lambda_n &=& -2 n^3 a^9+12 n^2 a^9-22 n a^9+12 a^9+6 n^4 a^8-14 n^3 a^8-66 n^2 a^8+206 n a^8-132 a^8\\
&& -6 n^5 a^7-24 n^4 a^7+212 n^3 a^7-132 n^2 a^7-542 n a^7+492 a^7+2 n^6 a^6+42 n^5 a^6\\
&& -112
   n^4 a^6-476 n^3 a^6+1274 n^2 a^6-238 n a^6-492 a^6-16 n^6 a^5-42 n^5 a^5+560 n^4 a^5\\
   && -658 n^3 a^5-2 D n^2 a^5-1456 n^2 a^5-4 D a^5+F a^5-2 G a^5+6 D n a^5-F n a^5+2884 n
   a^5\\
   && -1272 a^5+30 n^6 a^4-210 n^5 a^4-126 n^4 a^4+4 D n^3 a^4+2674 n^3 a^4-12 D n^2 a^4-F n^2 a^4\\
   && -3024 n^2 a^4-13 F a^4+16 G a^4+8 D n a^4+14 F n a^4-3136 n a^4+3792 a^4+40
   n^6 a^3\\
   && +336 n^5 a^3-2 D n^4 a^3-1736 n^4 a^3+8 D n^3 a^3+4 F n^3 a^3-992 n^3 a^3-4 D n^2 a^3\\
   && -15 F n^2 a^3+7072 n^2 a^3+12 D a^3+38 F a^3-30 G a^3-14 D n a^3-27 F n a^3-1168
   n a^3\\
   && -3552 a^3-152 n^6 a^2+168 n^5 a^2-2 D n^4 a^2-2 F n^4 a^2+1912 n^4 a^2-4 D n^3 a^2-4 F n^3 a^2\\
   &&-2184 n^3 a^2+22 D n^2 a^2+48 F n^2 a^2-4064 n^2 a^2-8 D a^2-32 F a^2-40 G
   a^2-8 D n a^2\\
   && -10 F n a^2+3168 n a^2+1152 a^2+96 n^6 a-288 n^5 a+4 D n^4 a+4 F n^4 a-480 n^4 a\\
   && -8 D n^3 a-8 F n^3 a+1440 n^3 a-4 D n^2 a-28 F n^2 a+384 n^2 a+152 G a+8 D n
   a\\
   && +32 F n a-1152 n a-96 G
\end{eqnarray*}
\end{small}
\item $J^{8.II.b}:$ \ $b-a-2=2,\ (b+a+2)<-6$
\begin{footnotesize}
\[
\begin{split}
 & p(x)  = {\left( 1 - x \right) }^{-3 - a}\,\left( 1 + x \right) \\
 & a_8 (x) ={\left( 1 - x^2 \right) }^4 ,\quad a_7 (x)  = -4\,\left( 4 + a - 6\,x + a\,x \right) \,{\left( -1 + x^2 \right) }^3 ,\\
 & a_6 (x) =  \frac{{\left( -1 + x^2 \right) }^2}{-96 + 176\,a - 90\,a^2 + 5\,a^3 + 6\,a^4 - a^5} X_6 \mbox{ with } \\
& \quad X_6 =   -5760 + 5376\,a + 3528\,a^2 - 3504\,a^3 + 90\,a^4 +
294\,a^5 - 18\,a^6 - 6\,a^7 + 4\,D \\
& \quad - 3\,a\,D - a^2\,D + 2\,F - a\,F +
  23040\,x - 41088\,a\,x + 18336\,a^2\,x + 1992\,a^3\,x \\
  & \quad - 2580\,a^4\,x + 228\,a^5\,x + 84\,a^6\,x - 12\,a^7\,x - 17280\,x^2 +
  38016\,a\,x^2 - 28392\,a^2\,x^2 \\
  & \quad + 7896\,a^3\,x^2 + 210\,a^4\,x^2 - 546\,a^5\,x^2 + 102\,a^6\,x^2 - 6\,a^7\,x^2 - 4\,D\,x^2 +
  3\,a\,D\,x^2 \\
  & \quad + a^2\,D\,x^2 - 2\,F\,x^2 + a\,F\,x^2   ,\\
 & a_5 (x)  = \frac{-1 + x^2}{-96 + 176\,a - 90\,a^2 + 5\,a^3 + 6\,a^4 - a^5} X_5 \mbox{ with } \\
 & \quad X_5 =  16896\,a - 25216\,a^2 + 5664\,a^3 + 3816\,a^4 - 996\,a^5 - 204\,a^6 + 36\,a^7 + 4\,a^8 - 48\,D \\
 & \quad + 24\,a\,D + 21\,a^2\,D +
  3\,a^3\,D - 24\,F + 6\,a\,F + 3\,a^2\,F - 46080\,x + 54528\,a\,x + 17472\,a^2\,x \\
  & \quad - 35088\,a^3\,x + 7728\,a^4\,x +
  2172\,a^5\,x - 732\,a^6\,x - 12\,a^7\,x + 12\,a^8\,x + 48\,D\,x - 48\,a\,D\,x \\
  & \quad - 3\,a^2\,D\,x + 3\,a^3\,D\,x + 24\,F\,x -
  18\,a\,F\,x + 3\,a^2\,F\,x + 92160\,x^2 - 187392\,a\,x^2 \\
  & \quad + 114432\,a^2\,x^2 - 10368\,a^3\,x^2 - 12312\,a^4\,x^2 +
  3492\,a^5\,x^2 + 108\,a^6\,x^2 - 132\,a^7\,x^2 \\
  & \quad + 12\,a^8\,x^2 + 48\,D\,x^2 - 24\,a\,D\,x^2 - 21\,a^2\,D\,x^2 -
  3\,a^3\,D\,x^2 + 24\,F\,x^2 - 6\,a\,F\,x^2 - 3\,a^2\,F\,x^2 \\
  & \quad - 46080\,x^3 + 112896\,a\,x^3 - 101056\,a^2\,x^3 +
  39984\,a^3\,x^3 - 4704\,a^4\,x^3 - 1596\,a^5\,x^3 \\
  & \quad + 636\,a^6\,x^3 - 84\,a^7\,x^3 + 4\,a^8\,x^3 - 48\,D\,x^3 + 48\,a\,D\,x^3 +
  3\,a^2\,D\,x^3 - 3\,a^3\,D\,x^3 - 24\,F\,x^3 \\
  & \quad + 18\,a\,F\,x^3 - 3\,a^2\,F\,x^3    ,\\
 & a_4(x)  = \frac{X_4}{2\,\left( -96 + 176\,a - 90\,a^2 + 5\,a^3 + 6\,a^4 - a^5 \right) } \mbox{ with }   \\
 & \quad X_4 = 23040 - 61056\,a + 33248\,a^2 + 18824\,a^3 - 14328\,a^4 - 1002\,a^5 + 1272\,a^6 + 36\,a^7 \\
 & \quad - 32\,a^8 - 2\,a^9 + 288\,D -
  210\,a^2\,D - 72\,a^3\,D - 6\,a^4\,D + 168\,F + 10\,a\,F - 33\,a^2\,F \\
  & \quad - 7\,a^3\,F + 101376\,a\,x - 185088\,a^2\,x +
  84416\,a^3\,x + 11568\,a^4\,x - 13608\,a^5\,x + 768\,a^6\,x \\
  & \quad + 624\,a^7\,x - 48\,a^8\,x - 8\,a^9\,x - 576\,D\,x +
  480\,a\,D\,x + 156\,a^2\,D\,x - 48\,a^3\,D\,x - 12\,a^4\,D\,x \\
  & \quad - 288\,F\,x + 168\,a\,F\,x + 12\,a^2\,F\,x - 12\,a^3\,F\,x -
  138240\,x^2 + 209664\,a\,x^2 - 2112\,a^2\,x^2 \\
  & \quad - 122736\,a^3\,x^2 + 58272\,a^4\,x^2 - 1212\,a^5\,x^2 - 4368\,a^6\,x^2 +
  696\,a^7\,x^2 + 48\,a^8\,x^2 - 12\,a^9\,x^2 \\
  & \quad - 384\,a\,D\,x^2 + 288\,a^2\,D\,x^2 + 96\,a^3\,D\,x^2 - 48\,F\,x^2 -
  140\,a\,F\,x^2 + 78\,a^2\,F\,x^2 + 2\,a^3\,F\,x^2 \\
  & \quad + 184320\,x^3 - 436224\,a\,x^3 + 353792\,a^2\,x^3 - 97024\,a^3\,x^3 -
  17712\,a^4\,x^3 + 15192\,a^5\,x^3 \\
  & \quad - 2112\,a^6\,x^3 - 336\,a^7\,x^3 + 112\,a^8\,x^3 - 8\,a^9\,x^3 + 576\,D\,x^3 -
  480\,a\,D\,x^3 - 156\,a^2\,D\,x^3 \\
  & \quad + 48\,a^3\,D\,x^3 + 12\,a^4\,D\,x^3 + 288\,F\,x^3 - 168\,a\,F\,x^3 - 12\,a^2\,F\,x^3 +
  12\,a^3\,F\,x^3 - 69120\,x^4 \\
  & \quad + 192384\,a\,x^4 - 208032\,a^2\,x^4 + 110504\,a^3\,x^4 - 27048\,a^4\,x^4 - 42\,a^5\,x^4 +
  1752\,a^6\,x^4 \\
  & \quad - 444\,a^7\,x^4 + 48\,a^8\,x^4 - 2\,a^9\,x^4 - 288\,D\,x^4 + 384\,a\,D\,x^4 - 78\,a^2\,D\,x^4 -
  24\,a^3\,D\,x^4 \\
  & \quad + 6\,a^4\,D\,x^4 - 120\,F\,x^4 + 130\,a\,F\,x^4 - 45\,a^2\,F\,x^4 + 5\,a^3\,F\,x^4  ,\\
   \end{split} \]
 \[
\begin{split}
 & a_3(x)  = \frac{X_3}{96 - 176\,a + 90\,a^2 - 5\,a^3 - 6\,a^4 + a^5} \mbox{ with }   \\
 & \quad X_3 =  -96\,D - 128\,a\,D + 114\,a^2\,D + 91\,a^3\,D + 18\,a^4\,D + a^5\,D - 144\,F + 4\,a\,F \\
 & \quad + 10\,a^2\,F + 8\,a^3\,F + 2\,a^4\,F +
  288\,D\,x - 144\,a\,D\,x - 210\,a^2\,D\,x + 33\,a^3\,D\,x  + 30\,a^4\,D\,x \\
  & \quad + 3\,a^5\,D\,x + 192\,F\,x - 112\,a\,F\,x -
  16\,a^2\,F\,x + 4\,a^3\,F\,x + 4\,a^4\,F\,x  - 288\,D\,x^2 \\
  & \quad + 384\,a\,D\,x^2 - 42\,a^2\,D\,x^2 - 63\,a^3\,D\,x^2 +
  6\,a^4\,D\,x^2 + 3\,a^5\,D\,x^2 - 48\,F\,x^2 + 76\,a\,F\,x^2 \\
  & \quad - 26\,a^2\,F\,x^2 - 4\,a^3\,F\,x^2 + 2\,a^4\,F\,x^2 +
  96\,D\,x^3 - 176\,a\,D\,x^3 + 90\,a^2\,D\,x^3 - 5\,a^3\,D\,x^3 \\
  & \quad - 6\,a^4\,D\,x^3 + a^5\,D\,x^3 ,\\
 & a_2(x)  = \frac{X_2}{2\,a\,\left( -4 + 3\,a + a^2 \right) } \mbox{ with } \\
& \quad X_2 =   14\,a\,F + 9\,a^2\,F + a^3\,F - 8\,G + 6\,a\,G +
2\,a^2\,G - 8\,a\,F\,x + 6\,a^2\,F\,x + 2\,a^3\,F\,x \\
& \quad + 2\,a\,F\,x^2 -
  3\,a^2\,F\,x^2 + a^3\,F\,x^2 + 8\,G\,x^2 - 6\,a\,G\,x^2 - 2\,a^2\,G\,x^2  ,\\
 & a_1(x)  = \frac{G\,\left( 4 + a + a\,x \right) }{a}
\end{split} \]
\end{footnotesize}
with eigenvalues
$$
\lambda_n = \frac{n (-a+n-1)}{2 (a-4) (a-3) (a-2) (a-1) a (a+4)} \Lambda_n
$$
where
\begin{small}
\begin{eqnarray*}
\Lambda_n &=& -2 n^3 a^9+12 n^2 a^9-22 n a^9+12 a^9+6 n^4 a^8-18 n^3 a^8-42 n^2 a^8+162 n a^8-108 a^8\\
&& -6 n^5 a^7-12 n^4 a^7+192 n^3 a^7-312 n^2 a^7-42 n a^7+180 a^7+2 n^6 a^6+30 n^5 a^6\\
&& -184
   n^4 a^6-12 n^3 a^6+1370 n^2 a^6-2322 n a^6+1116 a^6-12 n^6 a^5+66 n^5 a^5+480 n^4 a^5\\
   && -2538 n^3 a^5-2 D n^2 a^5+900 n^2 a^5-4 D a^5+F a^5-2 G a^5+6 D n a^5-F n a^5+6576 n
   a^5\\
   && -5472 a^5-10 n^6 a^4-510 n^5 a^4+1154 n^4 a^4+4 D n^3 a^4+5118 n^3 a^4-16 D n^2 a^4-F n^2 a^4\\
   && -10648 n^2 a^4-8 D a^4-13 F a^4+12 G a^4+20 D n a^4+14 F n a^4-4752 n
   a^4+9648 a^4\\
   && +180 n^6 a^3+516 n^5 a^3-2 D n^4 a^3-4548 n^4 a^3+16 D n^3 a^3+4 F n^3 a^3-532 n^3 a^3\\
   && -20 D n^2 a^3-15 F n^2 a^3+16272 n^2 a^3+28 D a^3+38 F a^3+10 G a^3-22 D n
   a^3-27 F n a^3\\
   && -4208 n a^3-7680 a^3-352 n^6 a^2+480 n^5 a^2-6 D n^4 a^2-2 F n^4 a^2+4064 n^4 a^2-4 D n^3 a^2\\
   && -4 F n^3 a^2-5088 n^3 a^2+46 D n^2 a^2+48 F n^2 a^2-8320 n^2
   a^2-16 D a^2-32 F a^2-180 G a^2\\
   && -20 D n a^2-10 F n a^2+6912 n a^2+2304 a^2+192 n^6 a-576 n^5 a+8 D n^4 a+4 F n^4 a\\
   && -960 n^4 a-16 D n^3 a-8 F n^3 a+2880 n^3 a-8 D n^2 a-28 F
   n^2 a+768 n^2 a+352 G a\\
   && +16 D n a+32 F n a-2304 n a-192 G
\end{eqnarray*}
\end{small}
\item $J^{8.II.c}:$ \ $b-a-2=4,\ (b+a+2)<-6$
\begin{footnotesize}
\[
\begin{split}
 & p(x)  = {\left( 1 - x \right) }^{-4 - a}\,{\left( 1 + x \right) }^2 \\
 & a_8 (x) ={\left( 1 - x^2 \right) }^4 ,\\
 & a_7 (x)  = -4\,\left( 6 + a - 6\,x + a\,x \right) \,{\left( -1 + x^2 \right) }^3 ,\\
 & a_6 (x) = \frac{{\left( -1 + x^2 \right) }^2}{3\,\left( -36 + 60\,a - 25\,a^2 + a^4 \right) }  X_6 \mbox{ with } \\
& \quad X_6 =   -19440 + 23976\,a - 108\,a^2 - 4770\,a^3 + 90\,a^4 +
234\,a^5 + 18\,a^6 - D + a\,D + F \\
& \quad + 38880\,x - 66096\,a\,x +
  27864\,a^2\,x + 1260\,a^3\,x - 1980\,a^4\,x + 36\,a^5\,x + 36\,a^6\,x \\
  & \quad - 19440\,x^2 + 39528\,a\,x^2 - 26028\,a^2\,x^2 +
  6030\,a^3\,x^2 + 90\,a^4\,x^2 - 198\,a^5\,x^2 + 18\,a^6\,x^2 \\
  & \quad + D\,x^2 - a\,D\,x^2 - F\,x^2   ,\\
 & a_5 (x)  = \frac{-1 + x^2}{-36 + 60\,a - 25\,a^2 + a^4}  X_5 \mbox{ with } \\
 & \quad X_5 =   17280 - 12960\,a - 11376\,a^2 + 6104\,a^3 + 1380\,a^4 - 340\,a^5 - 84\,a^6 - 4\,a^7 + 6\,D \\
 & \quad - 5\,a\,D - a^2\,D - 6\,F - a\,F -
  51840\,x + 76896\,a\,x - 16272\,a^2\,x - 12648\,a^3\,x + 3420\,a^4\,x \\
  & \quad + 564\,a^5\,x - 108\,a^6\,x - 12\,a^7\,x - 4\,D\,x +
  5\,a\,D\,x - a^2\,D\,x + 4\,F\,x - a\,F\,x + 51840\,x^2 \\
  & \quad - 101088\,a\,x^2 + 59184\,a^2\,x^2 - 7608\,a^3\,x^2 - 3060\,a^4\,x^2 +
  708\,a^5\,x^2 + 36\,a^6\,x^2 - 12\,a^7\,x^2 \\
  & \quad - 6\,D\,x^2 + 5\,a\,D\,x^2 + a^2\,D\,x^2 + 6\,F\,x^2 + a\,F\,x^2 - 17280\,x^3 +
  39456\,a\,x^3 - 31920\,a^2\,x^3 \\
  & \quad + 11144\,a^3\,x^3 - 1260\,a^4\,x^3 - 196\,a^5\,x^3 + 60\,a^6\,x^3 - 4\,a^7\,x^3 + 4\,D\,x^3 -
  5\,a\,D\,x^3 + a^2\,D\,x^3 \\
  & \quad - 4\,F\,x^3 + a\,F\,x^3   ,\\
 & a_4(x)  = \frac{X_4}{2\,\left( -36 + 60\,a - 25\,a^2 + a^4 \right) } \mbox{ with }   \\
 & \quad X_4 = -25920 - 8208\,a + 48744\,a^2 - 6300\,a^3 - 8902\,a^4 + 48\,a^5 + 476\,a^6 + 60\,a^7 + 2\,a^8 \\
 & \quad - 64\,D + 38\,a\,D + 24\,a^2\,D +
  2\,a^3\,D + 70\,F + 21\,a\,F + 3\,a^2\,F + 103680\,x - 112320\,a\,x \\
  & \quad - 42336\,a^2\,x + 59376\,a^3\,x - 3928\,a^4\,x -
  4800\,a^5\,x + 176\,a^6\,x + 144\,a^7\,x + 8\,a^8\,x \\
  & \quad + 72\,D\,x - 84\,a\,D\,x + 8\,a^2\,D\,x + 4\,a^3\,D\,x - 72\,F\,x +
  12\,a\,F\,x + 4\,a^2\,F\,x - 155520\,x^2 \\
  & \quad + 282528\,a\,x^2 - 125712\,a^2\,x^2 - 21672\,a^3\,x^2 + 22908\,a^4\,x^2 -
  1728\,a^5\,x^2 - 888\,a^6\,x^2 \\
  & \quad + 72\,a^7\,x^2 + 12\,a^8\,x^2 + 40\,D\,x^2 - 40\,a^2\,D\,x^2 - 52\,F\,x^2 - 30\,a\,F\,x^2 -
  2\,a^2\,F\,x^2 \\
  & \quad + 103680\,x^3 - 236736\,a\,x^3 + 185760\,a^2\,x^3 - 54672\,a^3\,x^3 - 1048\,a^4\,x^3 + 3456\,a^5\,x^3 \\
  & \quad -
  400\,a^6\,x^3 - 48\,a^7\,x^3 + 8\,a^8\,x^3 - 72\,D\,x^3 + 84\,a\,D\,x^3 - 8\,a^2\,D\,x^3 - 4\,a^3\,D\,x^3 \\
  & \quad + 72\,F\,x^3 -
  12\,a\,F\,x^3 - 4\,a^2\,F\,x^3 - 25920\,x^4 + 67824\,a\,x^4 - 67608\,a^2\,x^4 + 32676\,a^3\,x^4 \\
  & \quad - 7462\,a^4\,x^4 +
  336\,a^5\,x^4 + 188\,a^6\,x^4 - 36\,a^7\,x^4 + 2\,a^8\,x^4 + 24\,D\,x^4 - 38\,a\,D\,x^4 \\
  & \quad + 16\,a^2\,D\,x^4 - 2\,a^3\,D\,x^4 -
  18\,F\,x^4 + 9\,a\,F\,x^4 - a^2\,F\,x^4  ,\\
 & a_3(x)  = \frac{X_3}{3\,\left( -36 + 60\,a - 25\,a^2 + a^4 \right) } \mbox{ with }   \\
 & \quad X_3 =  -156\,D + 40\,a\,D + 95\,a^2\,D + 20\,a^3\,D + a^4\,D + 264\,F + 44\,a\,F + 24\,a^2\,F + 4\,a^3\,F \\
 & \quad + 192\,D\,x - 210\,a\,D\,x -
  15\,a^2\,D\,x + 30\,a^3\,D\,x + 3\,a^4\,D\,x - 228\,F\,x + 66\,a\,F\,x \\
  & \quad + 12\,a^2\,F\,x + 6\,a^3\,F\,x - 108\,D\,x^2 +
  180\,a\,D\,x^2 - 75\,a^2\,D\,x^2 + 3\,a^4\,D\,x^2 + 24\,D\,x^3 \\
  & \quad - 50\,a\,D\,x^3 + 35\,a^2\,D\,x^3 - 10\,a^3\,D\,x^3 +
  a^4\,D\,x^3 + 12\,F\,x^3 - 22\,a\,F\,x^3 + 12\,a^2\,F\,x^3 \\
  & \quad - 2\,a^3\,F\,x^3 ,\\
   \end{split} \]
 \[
\begin{split}
 & a_2(x)  = \frac{X_2}{2\,a\,\left( -6 + 5\,a + a^2 \right) } \mbox{ with } \\
& \quad X_2 =   34\,a\,F + 13\,a^2\,F + a^3\,F - 12\,G + 10\,a\,G +
2\,a^2\,G - 12\,a\,F\,x + 10\,a^2\,F\,x \\
& \quad + 2\,a^3\,F\,x + 2\,a\,F\,x^2 -
  3\,a^2\,F\,x^2 + a^3\,F\,x^2 + 12\,G\,x^2 - 10\,a\,G\,x^2 - 2\,a^2\,G\,x^2  ,\\
 & a_1(x)  = \frac{G\,\left( 6 + a + a\,x \right) }{a}
\end{split} \]
\end{footnotesize}
with eigenvalues
$$
\lambda_n = \frac{n (-a+n-1)}{6 (a-3) (a-2) (a-1) a (a+6)} \Lambda_n
$$
where
\begin{small}
\begin{eqnarray*}
\Lambda_n &=& -6 n^3 a^8+36 n^2 a^8-66 n a^8+36 a^8+18 n^4 a^7-90 n^3 a^7+90 n^2 a^7+90 n a^7-108 a^7\\
&& -18 n^5 a^6+72 n^4 a^6+156 n^3 a^6-1116 n^2 a^6+1734 n a^6-828 a^6+6 n^6 a^5-18 n^5
   a^5\\
   &&-480 n^4 a^5+1980 n^3 a^5-66 n^2 a^5-6282 n a^5+4860 a^5+450 n^5 a^4-720 n^4 a^4\\
   && -5334 n^3 a^4-2 D n^2 a^4+9504 n^2 a^4-4 D a^4+3 F a^4-6 G a^4+6 D n a^4-3 F n a^4\\
   && +5676 n
   a^4-9576 a^4-150 n^6 a^3-630 n^5 a^3+4422 n^4 a^3+4 D n^3 a^3+1350 n^3 a^3\\
   && -8 D n^2 a^3+F n^2 a^3-16800 n^2 a^3+8 D a^3-19 F a^3-4 D n a^3+18 F n a^3+3600 n a^3\\
   && +8208 a^3+360
   n^6 a^2-432 n^5 a^2-2 D n^4 a^2-4392 n^4 a^2+4 F n^3 a^2+5184 n^3 a^2\\
   && +12 D n^2 a^2-25 F n^2 a^2+9216 n^2 a^2-4 D a^2+22 F a^2+150 G a^2-6 D n a^2-F n a^2\\
   && -7344 n a^2-2592
   a^2-216 n^6 a+648 n^5 a+2 D n^4 a-2 F n^4 a+1080 n^4 a-4 D n^3 a\\
   && +4 F n^3 a-3240 n^3 a-2 D n^2 a+20 F n^2 a-864 n^2 a-360 G a+4 D n a-22 F n a\\
   && +2592 n a+216 G
\end{eqnarray*}
\end{small}
\item $J^{8.II.d}:$ \ $b-a-2=6,\ (b+a+2)<-6$
\begin{footnotesize}
\[
\begin{split}
 & p(x)  = {\left( 1 - x \right) }^{-5 - a}\,{\left( 1 + x \right) }^3 \\
 & a_8 (x) ={\left( 1 - x^2 \right) }^4 ,\\
 & a_7 (x)  = -4\,\left( 8 + a - 6\,x + a\,x \right) \,{\left( -1 + x^2 \right) }^3 ,\\
 & a_6 (x) =  \frac{{\left( -1 + x^2 \right) }^2}{3\,\left( -48 + 82\,a - 37\,a^2 + 2\,a^3 + a^4 \right) } X_6 \mbox{ with } \\
& \quad X_6 =  -50112 + 70920\,a - 14400\,a^2 - 7758\,a^3 + 990\,a^4
+ 342\,a^5 + 18\,a^6 - D + a\,D \\
& \quad + F + 69120\,x - 123264\,a\,x +
  60408\,a^2\,x - 3924\,a^3\,x - 2556\,a^4\,x + 180\,a^5\,x + 36\,a^6\,x \\
  & \quad - 25920\,x^2 + 53784\,a\,x^2 - 37080\,a^2\,x^2 +
  9882\,a^3\,x^2 - 522\,a^4\,x^2 - 162\,a^5\,x^2 + 18\,a^6\,x^2 \\
  & \quad + D\,x^2 - a\,D\,x^2 - F\,x^2    ,\\
 & a_5 (x)  = \frac{-1 + x^2}{-48 + 82\,a - 37\,a^2 + 2\,a^3 + a^4}  X_5 \mbox{ with } \\
 & \quad X_5 =  73728 - 87552\,a - 3584\,a^2 + 17864\,a^3 + 532\,a^4 - 868\,a^5 - 116\,a^6 - 4\,a^7 + 8\,D \\
 & \quad - 7\,a\,D - a^2\,D - 8\,F - a\,F -
  133632\,x + 222528\,a\,x - 85680\,a^2\,x - 11088\,a^3\,x \\
  & \quad + 7812\,a^4\,x + 252\,a^5\,x - 180\,a^6\,x - 12\,a^7\,x - 4\,D\,x +
  5\,a\,D\,x - a^2\,D\,x + 4\,F\,x - a\,F\,x \\
  & \quad + 92160\,x^2 - 187392\,a\,x^2 + 121632\,a^2\,x^2 - 25368\,a^3\,x^2 -
  2100\,a^4\,x^2 + 1092\,a^5\,x^2 \\
  & \quad - 12\,a^6\,x^2 - 12\,a^7\,x^2 - 8\,D\,x^2 + 7\,a\,D\,x^2 + a^2\,D\,x^2 + 8\,F\,x^2 +
  a\,F\,x^2 - 23040\,x^3 \\
  & \quad + 53568\,a\,x^3 - 44912\,a^2\,x^3 + 17024\,a^3\,x^3 - 2660\,a^4\,x^3 - 28\,a^5\,x^3 + 52\,a^6\,x^3 -
  4\,a^7\,x^3 \\
  & \quad + 4\,D\,x^3 - 5\,a\,D\,x^3 + a^2\,D\,x^3 - 4\,F\,x^3 + a\,F\,x^3    ,\\
 & a_4(x)  = \frac{X_4}{2\,\left( -48 + 82\,a - 37\,a^2 + 2\,a^3 + a^4 \right) } \mbox{ with }   \\
 & \quad X_4 =  -228096 + 202080\,a + 100952\,a^2 - 64120\,a^3 - 14966\,a^4 + 3080\,a^5 + 988\,a^6 + 80\,a^7 \\
 & \quad + 2\,a^8 - 120\,D + 86\,a\,D +
  32\,a^2\,D + 2\,a^3\,D + 126\,F + 29\,a\,F + 3\,a^2\,F + 442368\,x \\
  & \quad - 672768\,a\,x + 153600\,a^2\,x + 114352\,a^3\,x -
  32536\,a^4\,x - 6272\,a^5\,x + 1040\,a^6\,x \\
  & \quad + 208\,a^7\,x + 8\,a^8\,x + 96\,D\,x - 116\,a\,D\,x + 16\,a^2\,D\,x +
  4\,a^3\,D\,x - 96\,F\,x + 20\,a\,F\,x \\
  & \quad + 4\,a^2\,F\,x - 400896\,x^2 + 801216\,a\,x^2 - 479568\,a^2\,x^2 + 52416\,a^3\,x^2 +
  34524\,a^4\,x^2 \\
  & \quad - 7056\,a^5\,x^2 - 792\,a^6\,x^2 + 144\,a^7\,x^2 + 12\,a^8\,x^2 + 96\,D\,x^2 - 48\,a\,D\,x^2 -
  48\,a^2\,D\,x^2 \\
  & \quad - 108\,F\,x^2 - 38\,a\,F\,x^2 - 2\,a^2\,F\,x^2 + 184320\,x^3 - 436224\,a\,x^3 + 368192\,a^2\,x^3 \\
  & \quad -
  131824\,a^3\,x^3 + 12712\,a^4\,x^3 + 3584\,a^5\,x^3 - 752\,a^6\,x^3 - 16\,a^7\,x^3 + 8\,a^8\,x^3 - 96\,D\,x^3 \\
  & \quad +
  116\,a\,D\,x^3 - 16\,a^2\,D\,x^3 - 4\,a^3\,D\,x^3 + 96\,F\,x^3 - 20\,a\,F\,x^3 - 4\,a^2\,F\,x^3 - 34560\,x^4 \\
  & \quad + 91872\,a\,x^4 -
  94152\,a^2\,x^4 + 47992\,a^3\,x^4 - 12502\,a^4\,x^4 + 1288\,a^5\,x^4 + 92\,a^6\,x^4 \\
  & \quad - 32\,a^7\,x^4 + 2\,a^8\,x^4 +
  24\,D\,x^4 - 38\,a\,D\,x^4 + 16\,a^2\,D\,x^4 - 2\,a^3\,D\,x^4 - 18\,F\,x^4 \\
  & \quad + 9\,a\,F\,x^4 - a^2\,F\,x^4 ,\\
 & a_3(x)  = \frac{X_3}{3\,\left( -48 + 82\,a - 37\,a^2 + 2\,a^3 + a^4 \right) } \mbox{ with }   \\
 & \quad X_3 =  -432\,D + 226\,a\,D + 179\,a^2\,D + 26\,a^3\,D + a^4\,D + 576\,F + 104\,a\,F + 36\,a^2\,F \\
 & \quad + 4\,a^3\,F + 360\,D\,x -
  438\,a\,D\,x + 33\,a^2\,D\,x + 42\,a^3\,D\,x + 3\,a^4\,D\,x - 396\,F\,x + 126\,a\,F\,x \\
  & \quad + 24\,a^2\,F\,x + 6\,a^3\,F\,x -
  144\,D\,x^2 + 246\,a\,D\,x^2 - 111\,a^2\,D\,x^2 + 6\,a^3\,D\,x^2 + 3\,a^4\,D\,x^2 \\
  & \quad + 24\,D\,x^3 - 50\,a\,D\,x^3 +
  35\,a^2\,D\,x^3 - 10\,a^3\,D\,x^3 + a^4\,D\,x^3 + 12\,F\,x^3 - 22\,a\,F\,x^3 \\
  & \quad + 12\,a^2\,F\,x^3 - 2\,a^3\,F\,x^3  ,\\
   \end{split} \]
 \[
\begin{split}
 & a_2(x)  = \frac{X_2}{2\,a\,\left( -8 + 7\,a + a^2 \right) } \mbox{ with } \\
& \quad X_2 =  62\,a\,F + 17\,a^2\,F + a^3\,F - 16\,G + 14\,a\,G +
2\,a^2\,G - 16\,a\,F\,x + 14\,a^2\,F\,x \\
& \quad + 2\,a^3\,F\,x + 2\,a\,F\,x^2 -
  3\,a^2\,F\,x^2 + a^3\,F\,x^2 + 16\,G\,x^2 - 14\,a\,G\,x^2 - 2\,a^2\,G\,x^2   ,\\
 & a_1(x)  =  \frac{G\,\left( 8 + a + a\,x \right) }{a}
\end{split} \]
\end{footnotesize}
with eigenvalues
$$
\lambda_n = \frac{n (-a+n-1)}{6 (a-3) (a-2) (a-1) a (a+8)} \Lambda_n
$$
where
\begin{small}
\begin{eqnarray*}
\Lambda_n &=& -6 n^3 a^8+36 n^2 a^8-66 n a^8+36 a^8+18 n^4 a^7-102 n^3 a^7+162 n^2 a^7-42 n a^7-36 a^7\\
&& -18 n^5 a^6+108 n^4 a^6+48 n^3 a^6-1368 n^2 a^6+2706 n a^6-1476 a^6+6 n^6 a^5-54 n^5
   a^5\\
   && -552 n^4 a^5+2940 n^3 a^5-786 n^2 a^5-8646 n a^5+7092 a^5+12 n^6 a^4+630 n^5 a^4\\
   && -1248 n^4 a^4-7134 n^3 a^4-2 D n^2 a^4+13692 n^2 a^4-4 D a^4+3 F a^4-6 G a^4+6 D n a^4\\
   && -3 F
   n a^4+7296 n a^4-13248 a^4-222 n^6 a^3-810 n^5 a^3+6150 n^4 a^3+4 D n^3 a^3\\
   && +1482 n^3 a^3-8 D n^2 a^3+F n^2 a^3-22920 n^2 a^3+8 D a^3-19 F a^3-12 G a^3-4 D n a^3\\
   && +18 F n
   a^3+5232 n a^3+11088 a^3+492 n^6 a^2-612 n^5 a^2-2 D n^4 a^2-5916 n^4 a^2\\
   && +4 F n^3 a^2+7092 n^3 a^2+12 D n^2 a^2-25 F n^2 a^2+12336 n^2 a^2-4 D a^2+22 F a^2\\
   && +222 G a^2-6 D n
   a^2-F n a^2-9936 n a^2-3456 a^2-288 n^6 a+864 n^5 a+2 D n^4 a\\
   && -2 F n^4 a+1440 n^4 a-4 D n^3 a+4 F n^3 a-4320 n^3 a-2 D n^2 a+20 F n^2 a-1152 n^2 a\\
   && -492 G a+4 D n a-22 F n
   a+3456 n a+288 G
\end{eqnarray*}
\end{small}
\end{itemize}
\newpage
\noindent {\bf 4.4.3.3 The case $J^{8.III}:$\quad $(b-a-2)>6, \
(b+a+2) \in \left\{-6,-4,-2,0 \right\} $}
\begin{itemize}
\item $J^{8.III.a}:$ \ $(b-a-2)>6,\ b+a+2=0$
\begin{footnotesize}
\[
\begin{split}
 & p(x)  = {\left( 1 + x \right) }^{-2 - a} \\
 & a_8 (x) ={\left( 1 - x^2 \right) }^4 ,\quad  a_7 (x)  = -4\,\left( -2 + a\,\left( -1 + x \right)  - 6\,x \right) \,{\left( -1 + x^2 \right) }^3 ,\\
 & a_6 (x) = -\left( \frac{{\left( -1 + x^2 \right) }^2}{\left( 2 + a \right) \,\left( 24 - 50\,a + 35\,a^2 - 10\,a^3 + a^4 \right) } \right)  X_6 \mbox{ with } \\
& \quad X_6 =  576 - 2352\,a + 2232\,a^2 + 36\,a^3 - 666\,a^4 +
162\,a^5 + 18\,a^6 - 6\,a^7 + 2\,D - a\,D \\
& \quad - a^2\,D - 2\,F + a\,F - 5760\,x +
  7392\,a\,x + 912\,a^2\,x - 3432\,a^3\,x + 660\,a^4\,x + 348\,a^5\,x \\
  & \quad - 132\,a^6\,x + 12\,a^7\,x - 8640\,x^2 + 16848\,a\,x^2 -
  8904\,a^2\,x^2 - 924\,a^3\,x^2 + 2310\,a^4\,x^2 \\
  & \quad - 798\,a^5\,x^2 + 114\,a^6\,x^2 - 6\,a^7\,x^2 - 2\,D\,x^2 + a\,D\,x^2 +
  a^2\,D\,x^2 + 2\,F\,x^2 - a\,F\,x^2    ,\\
 & a_5 (x)  = -\left( \frac{-1 + x^2}{\left( 2 + a \right) \,\left( 24 - 50\,a + 35\,a^2 - 10\,a^3 + a^4 \right) } \right) X_5 \mbox{ with } \\
 & \quad X_5 =  4608 - 7680\,a + 800\,a^2 + 3824\,a^3 - 1304\,a^4 - 460\,a^5 + 220\,a^6 - 4\,a^7 - 4\,a^8 + 12\,D \\
 & \quad - 9\,a^2\,D - 3\,a^3\,D -
  12\,F + 3\,a^2\,F + 4608\,x - 19968\,a\,x + 22560\,a^2\,x - 4176\,a^3\,x \\
  & \quad - 5400\,a^4\,x + 2628\,a^5\,x - 180\,a^6\,x -
  84\,a^7\,x + 12\,a^8\,x + 24\,D\,x - 18\,a\,D\,x - 9\,a^2\,D\,x \\
  & \quad + 3\,a^3\,D\,x - 24\,F\,x + 18\,a\,F\,x - 3\,a^2\,F\,x -
  23040\,x^2 + 35328\,a\,x^2 - 3744\,a^2\,x^2 \\
  & \quad - 14640\,a^3\,x^2 + 6072\,a^4\,x^2 + 732\,a^5\,x^2 - 876\,a^6\,x^2 +
  180\,a^7\,x^2 - 12\,a^8\,x^2 - 12\,D\,x^2 \\
  & \quad + 9\,a^2\,D\,x^2 + 3\,a^3\,D\,x^2 + 12\,F\,x^2 - 3\,a^2\,F\,x^2 - 23040\,x^3 +
  50688\,a\,x^3 - 34976\,a^2\,x^3 \\
  & \quad + 3472\,a^3\,x^3 + 6776\,a^4\,x^3 - 3668\,a^5\,x^3 + 836\,a^6\,x^3 - 92\,a^7\,x^3 +
  4\,a^8\,x^3 - 24\,D\,x^3 \\
  & \quad + 18\,a\,D\,x^3 + 9\,a^2\,D\,x^3 - 3\,a^3\,D\,x^3 + 24\,F\,x^3 - 18\,a\,F\,x^3 + 3\,a^2\,F\,x^3    ,\\
 & a_4(x)  = \frac{X_4}{2\,\left( 2 + a \right) \,\left( 24 - 50\,a + 35\,a^2 - 10\,a^3 + a^4 \right) } \mbox{ with }   \\
 & \quad X_4 = -2304 - 4608\,a + 14320\,a^2 - 6312\,a^3 - 3308\,a^4 + 2426\,a^5 - 80\,a^6 - 148\,a^7 \\
 & \quad + 12\,a^8 + 2\,a^9 - 60\,a\,D +
  18\,a^2\,D + 36\,a^3\,D + 6\,a^4\,D + 24\,F + 34\,a\,F - 9\,a^2\,F - 7\,a^3\,F \\
  & \quad - 27648\,x + 55296\,a\,x - 20160\,a^2\,x -
  21344\,a^3\,x + 15472\,a^4\,x + 152\,a^5\,x - 2240\,a^6\,x \\
  & \quad + 464\,a^7\,x + 16\,a^8\,x - 8\,a^9\,x - 144\,D\,x + 48\,a\,D\,x +
  108\,a^2\,D\,x - 12\,a^4\,D\,x + 144\,F\,x \\
  & \quad - 48\,a\,F\,x - 36\,a^2\,F\,x + 12\,a^3\,F\,x - 13824\,x^2 + 64512\,a\,x^2 -
  87648\,a^2\,x^2 + 35088\,a^3\,x^2 \\
  & \quad + 12024\,a^4\,x^2 - 13284\,a^5\,x^2 + 3168\,a^6\,x^2 + 72\,a^7\,x^2 - 120\,a^8\,x^2 +
  12\,a^9\,x^2 - 144\,D\,x^2 \\
  & \quad + 216\,a\,D\,x^2 - 72\,a^3\,D\,x^2 + 96\,F\,x^2 - 164\,a\,F\,x^2 + 54\,a^2\,F\,x^2 +
  2\,a^3\,F\,x^2 + 46080\,x^3 \\
  & \quad - 86016\,a\,x^3 + 31040\,a^2\,x^3 + 26784\,a^3\,x^3 - 21904\,a^4\,x^3 + 2584\,a^5\,x^3 +
  2240\,a^6\,x^3 \\
  & \quad - 944\,a^7\,x^3 + 144\,a^8\,x^3 - 8\,a^9\,x^3 + 144\,D\,x^3 - 48\,a\,D\,x^3 - 108\,a^2\,D\,x^3 +
  12\,a^4\,D\,x^3 \\
  & \quad - 144\,F\,x^3 + 48\,a\,F\,x^3 + 36\,a^2\,F\,x^3 - 12\,a^3\,F\,x^3 + 34560\,x^4 - 87552\,a\,x^4 +
  77808\,a^2\,x^4 \\
  & \quad - 22696\,a^3\,x^4 - 8428\,a^4\,x^4 + 8890\,a^5\,x^4 - 3088\,a^6\,x^4 + 556\,a^7\,x^4 - 52\,a^8\,x^4 +
  2\,a^9\,x^4 \\
  & \quad + 144\,D\,x^4 - 156\,a\,D\,x^4 - 18\,a^2\,D\,x^4 + 36\,a^3\,D\,x^4 - 6\,a^4\,D\,x^4 - 120\,F\,x^4 +
  130\,a\,F\,x^4 \\
  & \quad - 45\,a^2\,F\,x^4 + 5\,a^3\,F\,x^4  ,\\
   \end{split} \]
 \[
\begin{split}
 & a_3(x)  = \frac{X_3}{48 - 76\,a + 20\,a^2 + 15\,a^3 - 8\,a^4 + a^5} \mbox{ with }   \\
 & \quad X_3 = -24\,D + 28\,a\,D + 22\,a^2\,D - 15\,a^3\,D - 10\,a^4\,D - a^5\,D - 24\,F - 2\,a^2\,F + 2\,a^4\,F \\
 & \quad + 60\,a\,D\,x -
  48\,a^2\,D\,x - 27\,a^3\,D\,x + 12\,a^4\,D\,x + 3\,a^5\,D\,x - 48\,F\,x + 16\,a\,F\,x \\
  & \quad + 4\,a^2\,F\,x + 8\,a^3\,F\,x -
  4\,a^4\,F\,x + 72\,D\,x^2 - 60\,a\,D\,x^2 - 42\,a^2\,D\,x^2 + 27\,a^3\,D\,x^2 \\
  & \quad + 6\,a^4\,D\,x^2 - 3\,a^5\,D\,x^2 - 24\,F\,x^2 +
  32\,a\,F\,x^2 - 2\,a^2\,F\,x^2 - 8\,a^3\,F\,x^2 + 2\,a^4\,F\,x^2 \\
  & \quad + 48\,D\,x^3 - 76\,a\,D\,x^3 + 20\,a^2\,D\,x^3 +
  15\,a^3\,D\,x^3 - 8\,a^4\,D\,x^3 + a^5\,D\,x^3  ,\\
 & a_2(x)  = \frac{X_2}{2\,a\,\left( -2 + a + a^2 \right) } \mbox{ with } \\
& \quad X_2 =  -2\,a\,F - 5\,a^2\,F - a^3\,F - 4\,G + 2\,a\,G +
2\,a^2\,G - 4\,a\,F\,x + 2\,a^2\,F\,x + 2\,a^3\,F\,x \\
& \quad - 2\,a\,F\,x^2 +
  3\,a^2\,F\,x^2 - a^3\,F\,x^2 + 4\,G\,x^2 - 2\,a\,G\,x^2 - 2\,a^2\,G\,x^2   ,\\
 & a_1(x)  = \frac{G\,\left( -2 + a\,\left( -1 + x \right)  \right) }{a}
\end{split} \]
\end{footnotesize}
with eigenvalues
$$
\lambda_n = \frac{n (-a+n-1)}{2 (a-4) (a-3) (a-2) (a-1) a (a+2)} \Lambda_n
$$
where
\begin{small}
\begin{eqnarray*}
\Lambda_n &=& -2 n^3 a^9+12 n^2 a^9-22 n a^9+12 a^9+6 n^4 a^8-14 n^3 a^8-66 n^2 a^8+206 n a^8-132 a^8\\
&& -6 n^5 a^7-24 n^4 a^7+212 n^3 a^7-132 n^2 a^7-542 n a^7+492 a^7+2 n^6 a^6+42 n^5 a^6\\
&& -112
   n^4 a^6-476 n^3 a^6+1274 n^2 a^6-238 n a^6-492 a^6-16 n^6 a^5-42 n^5 a^5+560 n^4 a^5\\
   && -658 n^3 a^5-2 D n^2 a^5-1456 n^2 a^5-4 D a^5-F a^5-2 G a^5+6 D n a^5+F n a^5\\
   && +2884 n
   a^5-1272 a^5+30 n^6 a^4-210 n^5 a^4-126 n^4 a^4+4 D n^3 a^4+2674 n^3 a^4\\
   && -12 D n^2 a^4+F n^2 a^4-3024 n^2 a^4+13 F a^4+16 G a^4+8 D n a^4-14 F n a^4-3136 n a^4\\
   && +3792 a^4+40
   n^6 a^3+336 n^5 a^3-2 D n^4 a^3-1736 n^4 a^3+8 D n^3 a^3-4 F n^3 a^3\\
   && -992 n^3 a^3-4 D n^2 a^3+15 F n^2 a^3+7072 n^2 a^3+12 D a^3-38 F a^3-30 G a^3-14 D n a^3\\
   && +27 F n a^3-1168
   n a^3-3552 a^3-152 n^6 a^2+168 n^5 a^2-2 D n^4 a^2+2 F n^4 a^2\\
   && +1912 n^4 a^2-4 D n^3 a^2+4 F n^3 a^2-2184 n^3 a^2+22 D n^2 a^2-48 F n^2 a^2-4064 n^2 a^2\\
   && -8 D a^2+32 F a^2-40 G
   a^2-8 D n a^2+10 F n a^2+3168 n a^2+1152 a^2+96 n^6 a\\
   && -288 n^5 a+4 D n^4 a-4 F n^4 a-480 n^4 a-8 D n^3 a+8 F n^3 a+1440 n^3 a-4 D n^2 a\\
   && +28 F n^2 a+384 n^2 a+152 G a+8 D n
   a-32 F n a-1152 n a-96 G
\end{eqnarray*}
\end{small}
\item $J^{8.III.b}:$ \ $(b-a-2)>6,\ b+a+2=-2$
\begin{footnotesize}
\[
\begin{split}
 & p(x)  = \left( 1 - x \right) \,{\left( 1 + x \right) }^{-3 - a} \\
 & a_8 (x) ={\left( 1 - x^2 \right) }^4 ,\\
 & a_7 (x)  = -4\,\left( -4 + a\,\left( -1 + x \right)  - 6\,x \right) \,{\left( -1 + x^2 \right) }^3 ,\\
 & a_6 (x) = \frac{{\left( -1 + x^2 \right) }^2}{-96 + 176\,a - 90\,a^2 + 5\,a^3 + 6\,a^4 - a^5}  X_6 \mbox{ with } \\
& \quad X_6 =  -5760 + 5376\,a + 3528\,a^2 - 3504\,a^3 + 90\,a^4 +
294\,a^5 - 18\,a^6 - 6\,a^7 + 4\,D - 3\,a\,D - a^2\,D - 2\,F \\
& \quad + a\,F
-
  23040\,x + 41088\,a\,x - 18336\,a^2\,x - 1992\,a^3\,x + 2580\,a^4\,x - 228\,a^5\,x - 84\,a^6\,x + 12\,a^7\,x \\
  & \quad - 17280\,x^2 +
  38016\,a\,x^2 - 28392\,a^2\,x^2 + 7896\,a^3\,x^2 + 210\,a^4\,x^2 - 546\,a^5\,x^2 + 102\,a^6\,x^2 - 6\,a^7\,x^2\\
  & \quad - 4\,D\,x^2 +
  3\,a\,D\,x^2 + a^2\,D\,x^2 + 2\,F\,x^2 - a\,F\,x^2    ,\\
 & a_5 (x)  = \frac{-1 + x^2}{-96 + 176\,a - 90\,a^2 + 5\,a^3 + 6\,a^4 - a^5} X_5 \mbox{ with } \\
 & \quad X_5 =  -16896\,a + 25216\,a^2 - 5664\,a^3 - 3816\,a^4 + 996\,a^5 + 204\,a^6 - 36\,a^7 - 4\,a^8 + 48\,D - 24\,a\,D\\
 & \quad - 21\,a^2\,D -
  3\,a^3\,D - 24\,F + 6\,a\,F + 3\,a^2\,F - 46080\,x + 54528\,a\,x + 17472\,a^2\,x - 35088\,a^3\,x \\
  & \quad + 7728\,a^4\,x +
  2172\,a^5\,x - 732\,a^6\,x - 12\,a^7\,x + 12\,a^8\,x + 48\,D\,x - 48\,a\,D\,x - 3\,a^2\,D\,x + 3\,a^3\,D\,x\\
  & \quad - 24\,F\,x +
  18\,a\,F\,x - 3\,a^2\,F\,x - 92160\,x^2 + 187392\,a\,x^2 - 114432\,a^2\,x^2 + 10368\,a^3\,x^2 + 12312\,a^4\,x^2\\
  & \quad -
  3492\,a^5\,x^2 - 108\,a^6\,x^2 + 132\,a^7\,x^2 - 12\,a^8\,x^2 - 48\,D\,x^2 + 24\,a\,D\,x^2 + 21\,a^2\,D\,x^2 +
  3\,a^3\,D\,x^2\\
  & \quad + 24\,F\,x^2 - 6\,a\,F\,x^2 - 3\,a^2\,F\,x^2 - 46080\,x^3 + 112896\,a\,x^3 - 101056\,a^2\,x^3 +
  39984\,a^3\,x^3 - 4704\,a^4\,x^3\\
  & \quad - 1596\,a^5\,x^3 + 636\,a^6\,x^3 - 84\,a^7\,x^3 + 4\,a^8\,x^3 - 48\,D\,x^3 + 48\,a\,D\,x^3 +
  3\,a^2\,D\,x^3 - 3\,a^3\,D\,x^3 + 24\,F\,x^3\\
  & \quad - 18\,a\,F\,x^3 + 3\,a^2\,F\,x^3    ,\\
 & a_4(x)  = \frac{X_4}{2\,\left( -96 + 176\,a - 90\,a^2 + 5\,a^3 + 6\,a^4 - a^5 \right) } \mbox{ with }   \\
 & \quad X_4 = 23040 - 61056\,a + 33248\,a^2 + 18824\,a^3 - 14328\,a^4 - 1002\,a^5 + 1272\,a^6 + 36\,a^7 - 32\,a^8 - 2\,a^9\\
 & \quad + 288\,D -
  210\,a^2\,D - 72\,a^3\,D - 6\,a^4\,D - 168\,F - 10\,a\,F + 33\,a^2\,F + 7\,a^3\,F - 101376\,a\,x + 185088\,a^2\,x\\
  & \quad -
  84416\,a^3\,x - 11568\,a^4\,x + 13608\,a^5\,x - 768\,a^6\,x - 624\,a^7\,x + 48\,a^8\,x + 8\,a^9\,x + 576\,D\,x -
  480\,a\,D\,x\\
  & \quad - 156\,a^2\,D\,x + 48\,a^3\,D\,x + 12\,a^4\,D\,x - 288\,F\,x + 168\,a\,F\,x + 12\,a^2\,F\,x - 12\,a^3\,F\,x -
  138240\,x^2\\
  & \quad + 209664\,a\,x^2 - 2112\,a^2\,x^2 - 122736\,a^3\,x^2 + 58272\,a^4\,x^2 - 1212\,a^5\,x^2 - 4368\,a^6\,x^2 +
  696\,a^7\,x^2\\
  & \quad + 48\,a^8\,x^2 - 12\,a^9\,x^2 - 384\,a\,D\,x^2 + 288\,a^2\,D\,x^2 + 96\,a^3\,D\,x^2 + 48\,F\,x^2 +
  140\,a\,F\,x^2 - 78\,a^2\,F\,x^2\\
  & \quad - 2\,a^3\,F\,x^2 - 184320\,x^3 + 436224\,a\,x^3 - 353792\,a^2\,x^3 + 97024\,a^3\,x^3 +
  17712\,a^4\,x^3 - 15192\,a^5\,x^3\\
  & \quad + 2112\,a^6\,x^3 + 336\,a^7\,x^3 - 112\,a^8\,x^3 + 8\,a^9\,x^3 - 576\,D\,x^3 +
  480\,a\,D\,x^3 + 156\,a^2\,D\,x^3 - 48\,a^3\,D\,x^3\\
  & \quad - 12\,a^4\,D\,x^3 + 288\,F\,x^3 - 168\,a\,F\,x^3 - 12\,a^2\,F\,x^3 +
  12\,a^3\,F\,x^3 - 69120\,x^4 + 192384\,a\,x^4 - 208032\,a^2\,x^4\\
  & \quad + 110504\,a^3\,x^4 - 27048\,a^4\,x^4 - 42\,a^5\,x^4 +
  1752\,a^6\,x^4 - 444\,a^7\,x^4 + 48\,a^8\,x^4 - 2\,a^9\,x^4 - 288\,D\,x^4\\
  & \quad + 384\,a\,D\,x^4 - 78\,a^2\,D\,x^4 -
  24\,a^3\,D\,x^4 + 6\,a^4\,D\,x^4 + 120\,F\,x^4 - 130\,a\,F\,x^4 + 45\,a^2\,F\,x^4 - 5\,a^3\,F\,x^4  ,\\
 & a_3(x)  = \frac{X_3}{96 - 176\,a + 90\,a^2 - 5\,a^3 - 6\,a^4 + a^5} \mbox{ with }   \\
 & \quad X_3 =  96\,D + 128\,a\,D - 114\,a^2\,D - 91\,a^3\,D - 18\,a^4\,D - a^5\,D - 144\,F + 4\,a\,F + 10\,a^2\,F + 8\,a^3\,F\\
 & \quad + 2\,a^4\,F +
  288\,D\,x - 144\,a\,D\,x - 210\,a^2\,D\,x + 33\,a^3\,D\,x + 30\,a^4\,D\,x + 3\,a^5\,D\,x - 192\,F\,x + 112\,a\,F\,x\\
  & \quad +
  16\,a^2\,F\,x - 4\,a^3\,F\,x - 4\,a^4\,F\,x + 288\,D\,x^2 - 384\,a\,D\,x^2 + 42\,a^2\,D\,x^2 + 63\,a^3\,D\,x^2 -
  6\,a^4\,D\,x^2\\
  & \quad - 3\,a^5\,D\,x^2 - 48\,F\,x^2 + 76\,a\,F\,x^2 - 26\,a^2\,F\,x^2 - 4\,a^3\,F\,x^2 + 2\,a^4\,F\,x^2 +
  96\,D\,x^3 - 176\,a\,D\,x^3\\
  & \quad + 90\,a^2\,D\,x^3 - 5\,a^3\,D\,x^3 - 6\,a^4\,D\,x^3 + a^5\,D\,x^3 ,\\
   \end{split} \]
 \[
\begin{split}
 & a_2(x)  = \frac{X_2}{2\,a\,\left( -4 + 3\,a + a^2 \right) } \mbox{ with } \\
& \quad X_2 =  -14\,a\,F - 9\,a^2\,F - a^3\,F - 8\,G + 6\,a\,G +
2\,a^2\,G - 8\,a\,F\,x + 6\,a^2\,F\,x + 2\,a^3\,F\,x - 2\,a\,F\,x^2\\
& \quad
+
  3\,a^2\,F\,x^2 - a^3\,F\,x^2 + 8\,G\,x^2 - 6\,a\,G\,x^2 - 2\,a^2\,G\,x^2   ,\\
 & a_1(x)  = \frac{G\,\left( -4 + a\,\left( -1 + x \right)  \right) }{a}
\end{split} \]
\end{footnotesize}
with eigenvalues
$$
\lambda_n = \frac{n (-a+n-1)}{2 (a-4) (a-3) (a-2) (a-1) a (a+4)} \Lambda_n
$$
where
\begin{small}
\begin{eqnarray*}
\Lambda_n &=& -2 n^3 a^9+12 n^2 a^9-22 n a^9+12 a^9+6 n^4 a^8-18 n^3 a^8-42 n^2 a^8+162 n a^8-108 a^8\\
&& -6 n^5 a^7-12 n^4 a^7+192 n^3 a^7-312 n^2 a^7-42 n a^7+180 a^7+2 n^6 a^6+30 n^5 a^6\\
&& -184
   n^4 a^6-12 n^3 a^6+1370 n^2 a^6-2322 n a^6+1116 a^6-12 n^6 a^5+66 n^5 a^5\\
   && +480 n^4 a^5-2538 n^3 a^5-2 D n^2 a^5+900 n^2 a^5-4 D a^5-F a^5-2 G a^5+6 D n a^5\\
   && +F n a^5+6576 n
   a^5-5472 a^5-10 n^6 a^4-510 n^5 a^4+1154 n^4 a^4+4 D n^3 a^4\\
   &&+5118 n^3 a^4-16 D n^2 a^4+F n^2 a^4-10648 n^2 a^4-8 D a^4+13 F a^4+12 G a^4\\
   && +20 D n a^4-14 F n a^4-4752 n
   a^4+9648 a^4+180 n^6 a^3+516 n^5 a^3-2 D n^4 a^3\\
   && -4548 n^4 a^3+16 D n^3 a^3-4 F n^3 a^3-532 n^3 a^3-20 D n^2 a^3+15 F n^2 a^3+16272 n^2 a^3\\
   && +28 D a^3-38 F a^3+10 G a^3-22 D n
   a^3+27 F n a^3-4208 n a^3-7680 a^3-352 n^6 a^2\\
   && +480 n^5 a^2-6 D n^4 a^2+2 F n^4 a^2+4064 n^4 a^2-4 D n^3 a^2+4 F n^3 a^2-5088 n^3 a^2\\
   && +46 D n^2 a^2-48 F n^2 a^2-8320 n^2
   a^2-16 D a^2+32 F a^2-180 G a^2-20 D n a^2\\
   &&+10 F n a^2+6912 n a^2+2304 a^2+192 n^6 a-576 n^5 a+8 D n^4 a-4 F n^4 a-960 n^4 a\\
   &&-16 D n^3 a+8 F n^3 a+2880 n^3 a-8 D n^2 a+28 F
   n^2 a+768 n^2 a+352 G a+16 D n a\\
   &&-32 F n a-2304 n a-192 G
\end{eqnarray*}
\end{small}
\item $J^{8.III.c}:$ \ $(b-a-2)>6,\ b+a+2=-4$
\begin{footnotesize}
\[
\begin{split}
 & p(x)  = {\left( 1 - x \right) }^2\,{\left( 1 + x \right) }^{-4 - a} \\
 & a_8 (x) ={\left( 1 - x^2 \right) }^4 ,\\
 & a_7 (x)  = -4\,{\left( -1 + x^2 \right) }^3\,\left( a\,\left( -1 + x \right)  - 6\,\left( 1 + x \right)  \right) ,\\
 & a_6 (x) = \frac{{\left( -1 + x^2 \right) }^2}{3\,\left( -36 + 60\,a - 25\,a^2 + a^4 \right) } X_6 \mbox{ with } \\
& \quad X_6 =   -19440 + 23976\,a - 108\,a^2 - 4770\,a^3 + 90\,a^4 +
234\,a^5 + 18\,a^6 + D - a\,D - F \\
& \quad - 38880\,x + 66096\,a\,x -
  27864\,a^2\,x - 1260\,a^3\,x + 1980\,a^4\,x - 36\,a^5\,x - 36\,a^6\,x \\
  & \quad - 19440\,x^2 + 39528\,a\,x^2 - 26028\,a^2\,x^2 +
  6030\,a^3\,x^2 + 90\,a^4\,x^2 - 198\,a^5\,x^2 + 18\,a^6\,x^2 \\
  & \quad - D\,x^2 + a\,D\,x^2 + F\,x^2   ,\\
 & a_5 (x)  = -\left( \frac{-1 + x^2}{-36 + 60\,a - 25\,a^2 + a^4} \right) X_5 \mbox{ with } \\
 & \quad X_5 =   17280 - 12960\,a - 11376\,a^2 + 6104\,a^3 + 1380\,a^4 - 340\,a^5 - 84\,a^6 - 4\,a^7 - 6\,D \\
 & \quad + 5\,a\,D + a^2\,D + 6\,F + a\,F +
  51840\,x - 76896\,a\,x + 16272\,a^2\,x + 12648\,a^3\,x \\
  & \quad - 3420\,a^4\,x - 564\,a^5\,x + 108\,a^6\,x + 12\,a^7\,x - 4\,D\,x +
  5\,a\,D\,x - a^2\,D\,x + 4\,F\,x - a\,F\,x \\
  & \quad + 51840\,x^2 - 101088\,a\,x^2 + 59184\,a^2\,x^2 - 7608\,a^3\,x^2 - 3060\,a^4\,x^2 +
  708\,a^5\,x^2 + 36\,a^6\,x^2 \\
  & \quad - 12\,a^7\,x^2 + 6\,D\,x^2 - 5\,a\,D\,x^2 - a^2\,D\,x^2 - 6\,F\,x^2 - a\,F\,x^2 + 17280\,x^3 -
  39456\,a\,x^3 \\
  & \quad + 31920\,a^2\,x^3 - 11144\,a^3\,x^3 + 1260\,a^4\,x^3 + 196\,a^5\,x^3 - 60\,a^6\,x^3 + 4\,a^7\,x^3 + 4\,D\,x^3 \\
  & \quad -
  5\,a\,D\,x^3 + a^2\,D\,x^3 - 4\,F\,x^3 + a\,F\,x^3   ,\\
 & a_4(x)  = \frac{X_4}{2\,\left( -36 + 60\,a - 25\,a^2 + a^4 \right) } \mbox{ with }   \\
 & \quad X_4 =  -25920 - 8208\,a + 48744\,a^2 - 6300\,a^3 - 8902\,a^4 + 48\,a^5 + 476\,a^6 + 60\,a^7 + 2\,a^8 \\
 & \quad + 64\,D - 38\,a\,D - 24\,a^2\,D -
  2\,a^3\,D - 70\,F - 21\,a\,F - 3\,a^2\,F - 103680\,x + 112320\,a\,x \\
  & \quad + 42336\,a^2\,x - 59376\,a^3\,x + 3928\,a^4\,x +
  4800\,a^5\,x - 176\,a^6\,x - 144\,a^7\,x - 8\,a^8\,x \\
  & \quad + 72\,D\,x - 84\,a\,D\,x + 8\,a^2\,D\,x + 4\,a^3\,D\,x - 72\,F\,x +
  12\,a\,F\,x + 4\,a^2\,F\,x - 155520\,x^2 \\
  & \quad + 282528\,a\,x^2 - 125712\,a^2\,x^2 - 21672\,a^3\,x^2 + 22908\,a^4\,x^2 -
  1728\,a^5\,x^2 - 888\,a^6\,x^2 \\
  & \quad + 72\,a^7\,x^2 + 12\,a^8\,x^2 - 40\,D\,x^2 + 40\,a^2\,D\,x^2 + 52\,F\,x^2 + 30\,a\,F\,x^2 +
  2\,a^2\,F\,x^2 \\
  & \quad - 103680\,x^3 + 236736\,a\,x^3 - 185760\,a^2\,x^3 + 54672\,a^3\,x^3 + 1048\,a^4\,x^3 - 3456\,a^5\,x^3 \\
  & \quad +
  400\,a^6\,x^3 + 48\,a^7\,x^3 - 8\,a^8\,x^3 - 72\,D\,x^3 + 84\,a\,D\,x^3 - 8\,a^2\,D\,x^3 - 4\,a^3\,D\,x^3 + 72\,F\,x^3 \\
  & \quad -
  12\,a\,F\,x^3 - 4\,a^2\,F\,x^3 - 25920\,x^4 + 67824\,a\,x^4 - 67608\,a^2\,x^4 + 32676\,a^3\,x^4 - 7462\,a^4\,x^4 \\
  & \quad +
  336\,a^5\,x^4 + 188\,a^6\,x^4 - 36\,a^7\,x^4 + 2\,a^8\,x^4 - 24\,D\,x^4 + 38\,a\,D\,x^4 - 16\,a^2\,D\,x^4 \\
  & \quad + 2\,a^3\,D\,x^4 +
  18\,F\,x^4 - 9\,a\,F\,x^4 + a^2\,F\,x^4 ,\\
 & a_3(x)  = \frac{X_3}{3\,\left( -36 + 60\,a - 25\,a^2 + a^4 \right) } \mbox{ with }   \\
 & \quad X_3 =  -156\,D + 40\,a\,D + 95\,a^2\,D + 20\,a^3\,D + a^4\,D + 264\,F + 44\,a\,F + 24\,a^2\,F + 4\,a^3\,F \\
 & \quad - 192\,D\,x + 210\,a\,D\,x +
  15\,a^2\,D\,x - 30\,a^3\,D\,x - 3\,a^4\,D\,x + 228\,F\,x - 66\,a\,F\,x \\
  & \quad - 12\,a^2\,F\,x - 6\,a^3\,F\,x - 108\,D\,x^2 +
  180\,a\,D\,x^2 - 75\,a^2\,D\,x^2 + 3\,a^4\,D\,x^2 - 24\,D\,x^3 \\
  & \quad + 50\,a\,D\,x^3 - 35\,a^2\,D\,x^3 + 10\,a^3\,D\,x^3 -
  a^4\,D\,x^3 - 12\,F\,x^3 + 22\,a\,F\,x^3 - 12\,a^2\,F\,x^3 \\
  & \quad + 2\,a^3\,F\,x^3 ,\\
   \end{split} \]
 \[
\begin{split}
 & a_2(x)  = \frac{X_2}{2\,a\,\left( -6 + 5\,a + a^2 \right) } \mbox{ with } \\
& \quad X_2 =   -34\,a\,F - 13\,a^2\,F - a^3\,F - 12\,G + 10\,a\,G +
2\,a^2\,G - 12\,a\,F\,x + 10\,a^2\,F\,x \\
& \quad + 2\,a^3\,F\,x - 2\,a\,F\,x^2 +
  3\,a^2\,F\,x^2 - a^3\,F\,x^2 + 12\,G\,x^2 - 10\,a\,G\,x^2 - 2\,a^2\,G\,x^2  ,\\
 & a_1(x)  = \frac{G\,\left( -6 + a\,\left( -1 + x \right)  \right) }{a}
\end{split} \]
\end{footnotesize}
with eigenvalues
$$
\lambda_n = \frac{n (-a+n-1)}{6 (a-3) (a-2) (a-1) a (a+6)} \Lambda_n
$$
where
\begin{small}
\begin{eqnarray*}
\Lambda_n &=& -6 n^3 a^8+36 n^2 a^8-66 n a^8+36 a^8+18 n^4 a^7-90 n^3 a^7+90 n^2 a^7+90 n a^7-108 a^7\\
&& -18 n^5 a^6+72 n^4 a^6+156 n^3 a^6-1116 n^2 a^6+1734 n a^6-828 a^6+6 n^6 a^5-18 n^5
   a^5\\
   && -480 n^4 a^5+1980 n^3 a^5-66 n^2 a^5-6282 n a^5+4860 a^5+450 n^5 a^4-720 n^4 a^4\\
   && -5334 n^3 a^4+2 D n^2 a^4+9504 n^2 a^4+4 D a^4-3 F a^4-6 G a^4-6 D n a^4+3 F n a^4\\
   && +5676 n
   a^4-9576 a^4-150 n^6 a^3-630 n^5 a^3+4422 n^4 a^3-4 D n^3 a^3+1350 n^3 a^3\\
   && +8 D n^2 a^3-F n^2 a^3-16800 n^2 a^3-8 D a^3+19 F a^3+4 D n a^3-18 F n a^3+3600 n a^3\\
   && +8208 a^3+360
   n^6 a^2-432 n^5 a^2+2 D n^4 a^2-4392 n^4 a^2-4 F n^3 a^2+5184 n^3 a^2\\
   && -12 D n^2 a^2+25 F n^2 a^2+9216 n^2 a^2+4 D a^2-22 F a^2+150 G a^2+6 D n a^2+F n a^2\\
   && -7344 n a^2-2592
   a^2-216 n^6 a+648 n^5 a-2 D n^4 a+2 F n^4 a+1080 n^4 a+4 D n^3 a\\
   &&-4 F n^3 a-3240 n^3 a+2 D n^2 a-20 F n^2 a-864 n^2 a-360 G a-4 D n a+22 F n a\\
   &&+2592 n a+216 G
\end{eqnarray*}
\end{small}
\item $J^{8.III.d}:$ \ $(b-a-2)>6,\ b+a+2=-6$
\begin{footnotesize}
\[
\begin{split}
 & p(x)  = {\left( 1 - x \right) }^3\,{\left( 1 + x \right) }^{-5 - a} \\
 & a_8 (x) ={\left( 1 - x^2 \right) }^4 ,\\
 & a_7 (x)  = -4\,\left( -8 + a\,\left( -1 + x \right)  - 6\,x \right) \,{\left( -1 + x^2 \right) }^3 ,\\
 & a_6 (x) = \frac{-{\left( -1 + x^2 \right) }^2}{3\,\left( -48 + 82\,a - 37\,a^2 + 2\,a^3 + a^4 \right) } X_6 \mbox{ with } \\
& \quad X_6 =  50112 - 70920\,a + 14400\,a^2 + 7758\,a^3 - 990\,a^4
- 342\,a^5 - 18\,a^6 - D + a\,D + F \\
& \quad + 69120\,x - 123264\,a\,x +
  60408\,a^2\,x - 3924\,a^3\,x - 2556\,a^4\,x + 180\,a^5\,x + 36\,a^6\,x \\
  & \quad + 25920\,x^2 - 53784\,a\,x^2 + 37080\,a^2\,x^2 -
  9882\,a^3\,x^2 + 522\,a^4\,x^2 + 162\,a^5\,x^2 \\
  & \quad - 18\,a^6\,x^2 + D\,x^2 - a\,D\,x^2 - F\,x^2    ,\\
 & a_5 (x)  =  -\left( \frac{-1 + x^2}{-48 + 82\,a - 37\,a^2 + 2\,a^3 + a^4} \right)  X_5 \mbox{ with } \\
 & \quad X_5 =    73728 - 87552\,a - 3584\,a^2 + 17864\,a^3 + 532\,a^4 - 868\,a^5 - 116\,a^6 - 4\,a^7 - 8\,D \\
 & \quad + 7\,a\,D + a^2\,D + 8\,F + a\,F +
  133632\,x - 222528\,a\,x + 85680\,a^2\,x + 11088\,a^3\,x \\
  & \quad - 7812\,a^4\,x - 252\,a^5\,x + 180\,a^6\,x + 12\,a^7\,x - 4\,D\,x +
  5\,a\,D\,x - a^2\,D\,x + 4\,F\,x - a\,F\,x \\
  & \quad + 92160\,x^2 - 187392\,a\,x^2 + 121632\,a^2\,x^2 - 25368\,a^3\,x^2 -
  2100\,a^4\,x^2 + 1092\,a^5\,x^2 \\
  & \quad - 12\,a^6\,x^2 - 12\,a^7\,x^2 + 8\,D\,x^2 - 7\,a\,D\,x^2 - a^2\,D\,x^2 - 8\,F\,x^2 -
  a\,F\,x^2 + 23040\,x^3 \\
  & \quad - 53568\,a\,x^3 + 44912\,a^2\,x^3 - 17024\,a^3\,x^3 + 2660\,a^4\,x^3 + 28\,a^5\,x^3 - 52\,a^6\,x^3 +
  4\,a^7\,x^3 \\
  & \quad + 4\,D\,x^3 - 5\,a\,D\,x^3 + a^2\,D\,x^3 - 4\,F\,x^3 + a\,F\,x^3  ,\\
 & a_4(x)  = \frac{X_4}{2\,\left( -48 + 82\,a - 37\,a^2 + 2\,a^3 + a^4 \right) } \mbox{ with }   \\
 & \quad X_4 =  -228096 + 202080\,a + 100952\,a^2 - 64120\,a^3 - 14966\,a^4 + 3080\,a^5 + 988\,a^6 \\
 & \quad + 80\,a^7 + 2\,a^8 + 120\,D - 86\,a\,D -
  32\,a^2\,D - 2\,a^3\,D - 126\,F - 29\,a\,F - 3\,a^2\,F \\
  & \quad - 442368\,x + 672768\,a\,x - 153600\,a^2\,x - 114352\,a^3\,x +
  32536\,a^4\,x + 6272\,a^5\,x \\
  & \quad - 1040\,a^6\,x - 208\,a^7\,x - 8\,a^8\,x + 96\,D\,x - 116\,a\,D\,x + 16\,a^2\,D\,x +
  4\,a^3\,D\,x - 96\,F\,x \\
  & \quad + 20\,a\,F\,x + 4\,a^2\,F\,x - 400896\,x^2 + 801216\,a\,x^2 - 479568\,a^2\,x^2 + 52416\,a^3\,x^2 \\
  & \quad +
  34524\,a^4\,x^2 - 7056\,a^5\,x^2 - 792\,a^6\,x^2 + 144\,a^7\,x^2 + 12\,a^8\,x^2 - 96\,D\,x^2 + 48\,a\,D\,x^2 \\
  & \quad +
  48\,a^2\,D\,x^2 + 108\,F\,x^2 + 38\,a\,F\,x^2 + 2\,a^2\,F\,x^2 - 184320\,x^3 + 436224\,a\,x^3 \\
  & \quad - 368192\,a^2\,x^3 +
  131824\,a^3\,x^3 - 12712\,a^4\,x^3 - 3584\,a^5\,x^3 + 752\,a^6\,x^3 + 16\,a^7\,x^3 \\
  & \quad - 8\,a^8\,x^3 - 96\,D\,x^3 +
  116\,a\,D\,x^3 - 16\,a^2\,D\,x^3 - 4\,a^3\,D\,x^3 + 96\,F\,x^3 - 20\,a\,F\,x^3 \\
  & \quad - 4\,a^2\,F\,x^3 - 34560\,x^4 + 91872\,a\,x^4 -
  94152\,a^2\,x^4 + 47992\,a^3\,x^4 - 12502\,a^4\,x^4 \\
  & \quad + 1288\,a^5\,x^4 + 92\,a^6\,x^4 - 32\,a^7\,x^4 + 2\,a^8\,x^4 -
  24\,D\,x^4 + 38\,a\,D\,x^4 - 16\,a^2\,D\,x^4 \\
  & \quad + 2\,a^3\,D\,x^4 + 18\,F\,x^4 - 9\,a\,F\,x^4 + a^2\,F\,x^4 ,\\
 & a_3(x)  = \frac{X_3}{3\,\left( -48 + 82\,a - 37\,a^2 + 2\,a^3 + a^4 \right) } \mbox{ with }   \\
 & \quad X_3 =  -432\,D + 226\,a\,D + 179\,a^2\,D + 26\,a^3\,D + a^4\,D + 576\,F + 104\,a\,F + 36\,a^2\,F \\
 & \quad + 4\,a^3\,F - 360\,D\,x +
  438\,a\,D\,x - 33\,a^2\,D\,x - 42\,a^3\,D\,x - 3\,a^4\,D\,x + 396\,F\,x \\
  & \quad - 126\,a\,F\,x - 24\,a^2\,F\,x - 6\,a^3\,F\,x -
  144\,D\,x^2 + 246\,a\,D\,x^2 - 111\,a^2\,D\,x^2 + 6\,a^3\,D\,x^2 \\
  & \quad + 3\,a^4\,D\,x^2 - 24\,D\,x^3 + 50\,a\,D\,x^3 -
  35\,a^2\,D\,x^3 + 10\,a^3\,D\,x^3 - a^4\,D\,x^3 - 12\,F\,x^3 \\
  & \quad + 22\,a\,F\,x^3 - 12\,a^2\,F\,x^3 + 2\,a^3\,F\,x^3 ,\\
   \end{split} \]
 \[
\begin{split}
 & a_2(x)  = \frac{X_2}{2\,a\,\left( -8 + 7\,a + a^2 \right) } \mbox{ with } \\
& \quad X_2 =  -62\,a\,F - 17\,a^2\,F - a^3\,F - 16\,G + 14\,a\,G +
2\,a^2\,G - 16\,a\,F\,x + 14\,a^2\,F\,x \\
& \quad + 2\,a^3\,F\,x - 2\,a\,F\,x^2 +
  3\,a^2\,F\,x^2 - a^3\,F\,x^2 + 16\,G\,x^2 - 14\,a\,G\,x^2 - 2\,a^2\,G\,x^2   ,\\
 & a_1(x)  = \frac{G\,\left( -8 + a\,\left( -1 + x \right)  \right) }{a}
\end{split} \]
\end{footnotesize}
with eigenvalues
$$
\lambda_n = \frac{n (-a+n-1)}{6 (a-3) (a-2) (a-1) a (a+8)}\Lambda_n
$$
where
\begin{small}
\begin{eqnarray*}
\Lambda_n &=& -6 n^3 a^8+36 n^2 a^8-66 n a^8+36 a^8+18 n^4 a^7-102 n^3 a^7+162 n^2 a^7-42 n a^7-36 a^7\\
&& -18 n^5 a^6+108 n^4 a^6+48 n^3 a^6-1368 n^2 a^6+2706 n a^6-1476 a^6+6 n^6 a^5-54 n^5
   a^5\\
   &&-552 n^4 a^5+2940 n^3 a^5-786 n^2 a^5-8646 n a^5+7092 a^5+12 n^6 a^4+630 n^5 a^4\\
   &&-1248 n^4 a^4-7134 n^3 a^4+2 D n^2 a^4+13692 n^2 a^4+4 D a^4-3 F a^4-6 G a^4-6 D n a^4\\
   &&+3 F
   n a^4+7296 n a^4-13248 a^4-222 n^6 a^3-810 n^5 a^3+6150 n^4 a^3-4 D n^3 a^3\\
   &&+1482 n^3 a^3+8 D n^2 a^3-F n^2 a^3-22920 n^2 a^3-8 D a^3+19 F a^3-12 G a^3+4 D n a^3\\
   &&-18 F n
   a^3+5232 n a^3+11088 a^3+492 n^6 a^2-612 n^5 a^2+2 D n^4 a^2-5916 n^4 a^2\\
   &&-4 F n^3 a^2+7092 n^3 a^2-12 D n^2 a^2+25 F n^2 a^2+12336 n^2 a^2+4 D a^2-22 F a^2\\
   &&+222 G a^2+6 D n
   a^2+F n a^2-9936 n a^2-3456 a^2-288 n^6 a+864 n^5 a-2 D n^4 a\\
   &&+2 F n^4 a+1440 n^4 a+4 D n^3 a-4 F n^3 a-4320 n^3 a+2 D n^2 a-20 F n^2 a-1152 n^2 a\\
   &&-492 G a-4 D n a+22 F n
   a+3456 n a+288 G
\end{eqnarray*}
\end{small}
\end{itemize}
\newpage
\noindent {\bf 4.4.3.4 The case $J^{8.IV}:$\quad $(b-a-2)>6,\
(b+a+2)<-6 $}

%
\begin{footnotesize}
\[
\begin{split}
 & p(x)   = \fr{(1+x)^{\fr{b-a-2}{2}}}{(1-x)^{\fr{b+a+2}{2}}} \\
 & a_8 (x) = {\left( 1 - x^2 \right) }^4 ,\\
 & a_7 (x)  = -4\,\left( b + \left( -6 + a \right) \,x \right) \,{\left( -1 + x^2 \right) }^3 ,\\
 & a_6 (x) = \frac{{\left( -1 + x^2 \right) }^2}{a\,\left( 24 - 50\,a + 35\,a^2 - 10\,a^3 + a^4 \right) } X_6 \mbox{ with } \\
& \quad X_6 =  -864\,a + 1944\,a^2 - 1560\,a^3 + 570\,a^4 - 96\,a^5
+
  6\,a^6 + 144\,a\,b^2 - 300\,a^2\,b^2 +
  210\,a^3\,b^2 \\
  & \quad - 60\,a^4\,b^2 + 6\,a^5\,b^2 - a\,D +
  a^2\,D + 2\,a\,F + 2\,G - 1440\,a\,b\,x +
  3288\,a^2\,b\,x - 2700\,a^3\,b\,x \\
  & \quad +
  1020\,a^4\,b\,x - 180\,a^5\,b\,x + 12\,a^6\,b\,x +
  4320\,a\,x^2 - 10584\,a^2\,x^2 + 9744\,a^3\,x^2 -
  4410\,a^4\,x^2 \\
  & \quad + 1050\,a^5\,x^2 - 126\,a^6\,x^2 +
  6\,a^7\,x^2 + a\,D\,x^2 - a^2\,D\,x^2 -
  2\,a\,F\,x^2 - 2\,G\,x^2  ,\\
 & a_5 (x)  = -\left( \frac{-1 + x^2}{a\,\left( 24 - 50\,a + 35\,a^2 - 10\,a^3 + a^4 \right) } \right) X_5 \mbox{ with }\\
& \quad X_5 =  -1536\,a\,b + 3488\,a^2\,b - 2840\,a^3\,b +
  1060\,a^4\,b - 184\,a^5\,b + 12\,a^6\,b +
  96\,a\,b^3 - 200\,a^2\,b^3 \\
  & \quad + 140\,a^3\,b^3 -
  40\,a^4\,b^3 + 4\,a^5\,b^3 - 3\,a\,b\,D +
  3\,a^2\,b\,D + 6\,a\,b\,F + 6\,b\,G + 6912\,a\,x -
  17280\,a^2\,x \\
  & \quad + 16368\,a^3\,x - 7680\,a^4\,x +
  1908\,a^5\,x - 240\,a^6\,x + 12\,a^7\,x -
  1152\,a\,b^2\,x + 2688\,a^2\,b^2\,x \\
  & \quad -
  2280\,a^3\,b^2\,x + 900\,a^4\,b^2\,x -
  168\,a^5\,b^2\,x + 12\,a^6\,b^2\,x + 12\,a\,D\,x -
  15\,a^2\,D\,x + 3\,a^3\,D\,x \\
  & \quad - 24\,a\,F\,x +
  6\,a^2\,F\,x - 24\,G\,x + 6\,a\,G\,x +
  5760\,a\,b\,x^2 - 14592\,a^2\,b\,x^2 +
  14088\,a^3\,b\,x^2 \\
  & \quad - 6780\,a^4\,b\,x^2 +
  1740\,a^5\,b\,x^2 - 228\,a^6\,b\,x^2 +
  12\,a^7\,b\,x^2 + 3\,a\,b\,D\,x^2 -
  3\,a^2\,b\,D\,x^2 - 6\,a\,b\,F\,x^2 \\
  & \quad - 6\,b\,G\,x^2 -
  11520\,a\,x^3 + 31104\,a^2\,x^3 - 33040\,a^3\,x^3 +
  18256\,a^4\,x^3 - 5740\,a^5\,x^3 + 1036\,a^6\,x^3 \\
  & \quad -
  100\,a^7\,x^3 + 4\,a^8\,x^3 - 12\,a\,D\,x^3 +
  15\,a^2\,D\,x^3 - 3\,a^3\,D\,x^3 + 24\,a\,F\,x^3 -
  6\,a^2\,F\,x^3 \\
  & \quad + 24\,G\,x^3 - 6\,a\,G\,x^3 ,\\
 & a_4(x)  = \frac{X_4}{a\,\left( 24 - 50\,a + 35\,a^2 - 10\,a^3 + a^4 \right) } \mbox{ with } \\
& \quad X_4 =  1728\,a - 4320\,a^2 + 4092\,a^3 - 1920\,a^4 +
  477\,a^5 - 60\,a^6 + 3\,a^7 - 672\,a\,b^2 +
  1544\,a^2\,b^2 \\
  & \quad - 1280\,a^3\,b^2 + 490\,a^4\,b^2 -
  88\,a^5\,b^2 + 6\,a^6\,b^2 + 24\,a\,b^4 -
  50\,a^2\,b^4 + 35\,a^3\,b^4 - 10\,a^4\,b^4 +
  a^5\,b^4 \\
  & \quad + 12\,a\,D - 15\,a^2\,D + 3\,a^3\,D -
  3\,a\,b^2\,D + 3\,a^2\,b^2\,D - 12\,a\,F - a^2\,F +
  a^3\,F + 6\,a\,b^2\,F - 12\,G \\
  & \quad - a\,G + a^2\,G +
  6\,b^2\,G + 4608\,a\,b\,x - 12000\,a^2\,b\,x +
  12008\,a^3\,b\,x - 6020\,a^4\,b\,x +
  1612\,a^5\,b\,x \\
  & \quad - 220\,a^6\,b\,x + 12\,a^7\,b\,x -
  288\,a\,b^3\,x + 696\,a^2\,b^3\,x -
  620\,a^3\,b^3\,x + 260\,a^4\,b^3\,x -
  52\,a^5\,b^3\,x \\
  & \quad + 4\,a^6\,b^3\,x + 18\,a\,b\,D\,x -
  24\,a^2\,b\,D\,x + 6\,a^3\,b\,D\,x -
  36\,a\,b\,F\,x + 12\,a^2\,b\,F\,x - 36\,b\,G\,x \\
  & \quad +
  12\,a\,b\,G\,x - 10368\,a\,x^2 + 29376\,a^2\,x^2 -
  33192\,a^3\,x^2 + 19704\,a^4\,x^2 - 6702\,a^5\,x^2 +
  1314\,a^6\,x^2 \\
  & \quad - 138\,a^7\,x^2 + 6\,a^8\,x^2 +
  1728\,a\,b^2\,x^2 - 4608\,a^2\,b^2\,x^2 +
  4764\,a^3\,b^2\,x^2 - 2490\,a^4\,b^2\,x^2 +
  702\,a^5\,b^2\,x^2 \\
  & \quad - 102\,a^6\,b^2\,x^2 +
  6\,a^7\,b^2\,x^2 - 48\,a\,D\,x^2 + 72\,a^2\,D\,x^2 -
  27\,a^3\,D\,x^2 + 3\,a^4\,D\,x^2 +
  3\,a\,b^2\,D\,x^2 \\
  & \quad - 3\,a^2\,b^2\,D\,x^2 +
  72\,a\,F\,x^2 - 34\,a^2\,F\,x^2 + 4\,a^3\,F\,x^2 -
  6\,a\,b^2\,F\,x^2 + 72\,G\,x^2 - 34\,a\,G\,x^2 \\
  & \quad +
  4\,a^2\,G\,x^2 - 6\,b^2\,G\,x^2 - 5760\,a\,b\,x^3 +
  16512\,a^2\,b\,x^3 - 18952\,a^3\,b\,x^3 +
  11476\,a^4\,b\,x^3 \\
  & \quad - 4000\,a^5\,b\,x^3 +
  808\,a^6\,b\,x^3 - 88\,a^7\,b\,x^3 +
  4\,a^8\,b\,x^3 - 18\,a\,b\,D\,x^3 +
  24\,a^2\,b\,D\,x^3 - 6\,a^3\,b\,D\,x^3 \\
  & \quad +
  36\,a\,b\,F\,x^3 - 12\,a^2\,b\,F\,x^3 +
  36\,b\,G\,x^3 - 12\,a\,b\,G\,x^3 + 8640\,a\,x^4 -
  26208\,a^2\,x^4 + 32556\,a^3\,x^4 \\
  & \quad -
  21952\,a^4\,x^4 + 8869\,a^5\,x^4 - 2212\,a^6\,x^4 +
  334\,a^7\,x^4 - 28\,a^8\,x^4 + a^9\,x^4 +
  36\,a\,D\,x^4 - 57\,a^2\,D\,x^4 \\
  & \quad + 24\,a^3\,D\,x^4 -
  3\,a^4\,D\,x^4 - 60\,a\,F\,x^4 + 35\,a^2\,F\,x^4 -
  5\,a^3\,F\,x^4 - 60\,G\,x^4 + 35\,a\,G\,x^4 -
  5\,a^2\,G\,x^4 ,\\
   \end{split} \]
 \[
\begin{split}
 & a_3(x)  = \frac{X_3}{a\,\left( 24 - 50\,a + 35\,a^2 - 10\,a^3 + a^4 \right) } \mbox{ with }  \\
 & \quad X_3 = 10\,a\,b\,D - 13\,a^2\,b\,D + 3\,a^3\,b\,D -
  a\,b^3\,D + a^2\,b^3\,D + 4\,a\,b\,F -
  8\,a^2\,b\,F + 2\,a^3\,b\,F \\
  & \quad + 2\,a\,b^3\,F +
  4\,b\,G - 8\,a\,b\,G + 2\,a^2\,b\,G + 2\,b^3\,G -
  24\,a\,D\,x + 42\,a^2\,D\,x - 21\,a^3\,D\,x +
  3\,a^4\,D\,x \\
  & \quad + 6\,a\,b^2\,D\,x - 9\,a^2\,b^2\,D\,x +
  3\,a^3\,b^2\,D\,x + 16\,a^2\,F\,x - 12\,a^3\,F\,x +
  2\,a^4\,F\,x - 12\,a\,b^2\,F\,x \\
  & \quad +
  6\,a^2\,b^2\,F\,x + 16\,a\,G\,x - 12\,a^2\,G\,x +
  2\,a^3\,G\,x - 12\,b^2\,G\,x + 6\,a\,b^2\,G\,x -
  18\,a\,b\,D\,x^2 \\
  & \quad + 33\,a^2\,b\,D\,x^2 -
  18\,a^3\,b\,D\,x^2 + 3\,a^4\,b\,D\,x^2 +
  12\,a\,b\,F\,x^2 - 16\,a^2\,b\,F\,x^2 +
  4\,a^3\,b\,F\,x^2 + 12\,b\,G\,x^2 \\
  & \quad -
  16\,a\,b\,G\,x^2 + 4\,a^2\,b\,G\,x^2 +
  24\,a\,D\,x^3 - 50\,a^2\,D\,x^3 + 35\,a^3\,D\,x^3 -
  10\,a^4\,D\,x^3 + a^5\,D\,x^3 ,\\
 & a_2(x)  = \frac{X_2}{a\,\left( 2 - 3\,a + a^2 \right) } \mbox{ with } \\
& \quad X_2 = -2\,a\,F + a^2\,F + a\,b^2\,F - 2\,a\,G + a^2\,G +
  b^2\,G - 2\,a\,b\,F\,x + 2\,a^2\,b\,F\,x \\
  & \quad -
  2\,b\,G\,x + 2\,a\,b\,G\,x + 2\,a\,F\,x^2 -
  3\,a^2\,F\,x^2 + a^3\,F\,x^2 ,\\
 & a_1(x)  = \frac{G\,\left( b + a\,x \right) }{a}
\end{split} \]
\end{footnotesize}
with eigenvalues
$$
\lambda_n = \frac{n (-a+n-1)}{(a-4) (a-3) (a-2) (a-1) a} \Lambda_n
$$
where
\begin{small}
\begin{eqnarray*}
\Lambda_n &=& -n^3 a^8+6 n^2 a^8-11 n a^8+6 a^8+3 n^4 a^7-5 n^3 a^7-45 n^2 a^7+125 n a^7-78 a^7-3 n^5 a^6\\
&& -18 n^4 a^6+116 n^3 a^6+24 n^2 a^6-521 n a^6+402 a^6+n^6 a^5+27 n^5 a^5-20 n^4
   a^5-470 n^3 a^5\\
   &&+589 n^2 a^5+923 n a^5-1050 a^5-10 n^6 a^4-75 n^5 a^4+320 n^4 a^4+611 n^3 a^4-D n^2 a^4\\
   &&-1906 n^2 a^4-2 D a^4+F a^4-G a^4+3 D n a^4-F n a^4-404 n a^4+1464
   a^4+35 n^6 a^3\\
   &&+45 n^5 a^3-703 n^4 a^3+2 D n^3 a^3+115 n^3 a^3-4 D n^2 a^3-F n^2 a^3+2300 n^2 a^3+4 D a^3\\
   &&-11 F a^3+11 G a^3-2 D n a^3+12 F n a^3-G n a^3-760 n a^3-1032 a^3-50
   n^6 a^2+78 n^5 a^2\\
   &&-D n^4 a^2+538 n^4 a^2+4 F n^3 a^2-726 n^3 a^2+6 D n^2 a^2-17 F n^2 a^2-G n^2 a^2-1064 n^2 a^2\\
   &&-2 D a^2+16 F a^2-46 G a^2-3 D n a^2-3 F n a^2+12 G n a^2+936
   n a^2+288 a^2+24 n^6 a\\
   &&-72 n^5 a+D n^4 a-2 F n^4 a-120 n^4 a-2 D n^3 a+4 F n^3 a+4 G n^3 a+360 n^3 a-D n^2 a\\
   &&+14 F n^2 a-17 G n^2 a+96 n^2 a+66 G a+2 D n a-16 F n a-3 G n
   a-288 n a-2 G n^4\\
   &&+4 G n^3+14 G n^2-24 G-16 G n
\end{eqnarray*}
\end{small}
%
The zeros of denominator terms in above expressions consist of $a=0,
a=1$, $a=2$, $a=3$ and $a=4$. Since $(b-a-2)>6$ and $(b+a+2)<-6 $
imply $a<-8$, the denominator terms in above expressions never
vanish.\\

\subsection*{Acknowledgement}
We thank Professor Mourad Ismail for his valuable comments and suggestions during the preparation of this paper.

\end{document}